## Cographs[i]

Robert Haas

### 1. Introduction

This article introduces the theory of cographs, combinatorial structures that both dualize and generalize ordinary graphs; they promise as rich a range as graph theory itself of applications inside and outside mathematics. Their double origin produces two equivalent definitions. For the first, note that a graph may be defined abstractly as a function (a primitive boundary operator) $\mathcal{G}: \mathcal{E} \rightarrow \mathcal{P}^{(2)}$ that associates to each element of an abstract set $\mathcal{E}$ ("edges") an (unordered) pair of distinct elements from the abstract set $\mathcal{P}$ ("points") called the element's "endpoints." Dualizing:

**Definition 1.1.** A *cograph* is a set function $\mathcal{C}: \mathcal{P}^{(2)} \rightarrow \mathcal{E}$, where elements of $\mathcal{P}$ are "points," elements of $\mathcal{E}$ are "edges," and $\mathcal{P}^{(2)}$ is (unordered) pairs of distinct points of $\mathcal{P}$.

Representing the cograph $\mathcal{C}$ by ordinary points, and lines between each pair of points, yields the equivalent second definition:

**Definition 1.2.** A *cograph* is a complete graph with colored lines.

Here the *complete graph* $K_n$ is the graph on $n$ points where each pair is joined by a line. In keeping with the desired abstract combinatorial nature of cographs, the second definition actually involves *classes* of graphs *equivalent* under position or color labeling: In the cographs with three points, for example, no distinction will be made between the triangle having a red lower line and blue right and left sides, and the same triangle rotated to put its red line at the right; or from the triangle with one yellow and two green lines; or from the one with one blue and two red lines.

**Remark.** Note that it is not required that lines incident at a point have different colors. This type of constraint, common in graph coloring theory, can introduce considerable combinatorial complexity (e.g. Haas [10]), and examples below show it does not hold in some quite natural classes of cographs.

Figure 1.1 shows the three cographs on 3 points, representing the colors by patternings--solid, dashed, or dotted--of the edges; "Type" classifies the cographs by the number of copies of each edge.

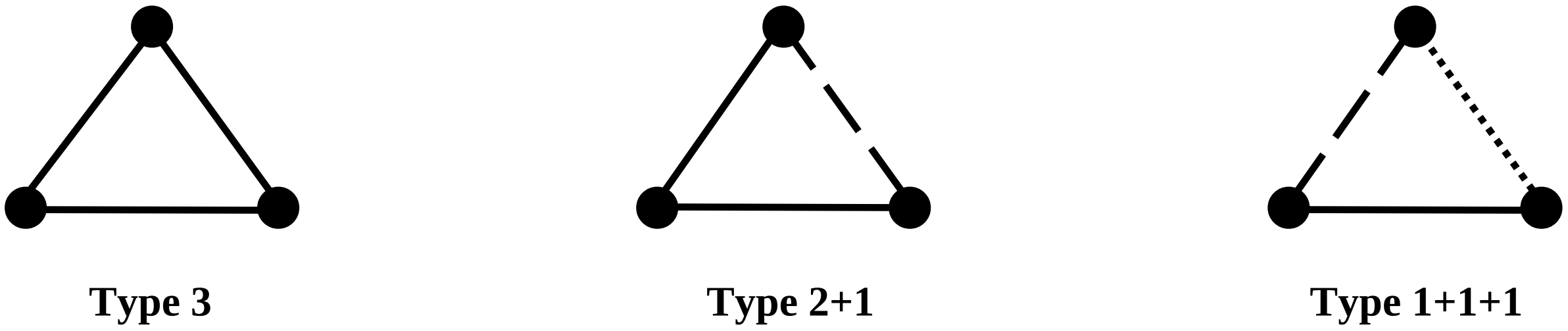

**Fig. 1.1.** The three cographs on three points.

Figure 1.2 gives the cographs on 2, 3, or 4 points (for clarity showing only multiple-copy edges). The problem how to enumerate cographs (i.e. complete graphs with colored lines) was first solved by R. C. Read, and Table 1.1 lists the number of cographs on $n \leq 9$ points calculated by the computer algebra program Maple from the enumeration theorem of Polya (Harary & Palmer [12] 35-6, 41-3, 82-8, 135-7, 142-5, 251-2).

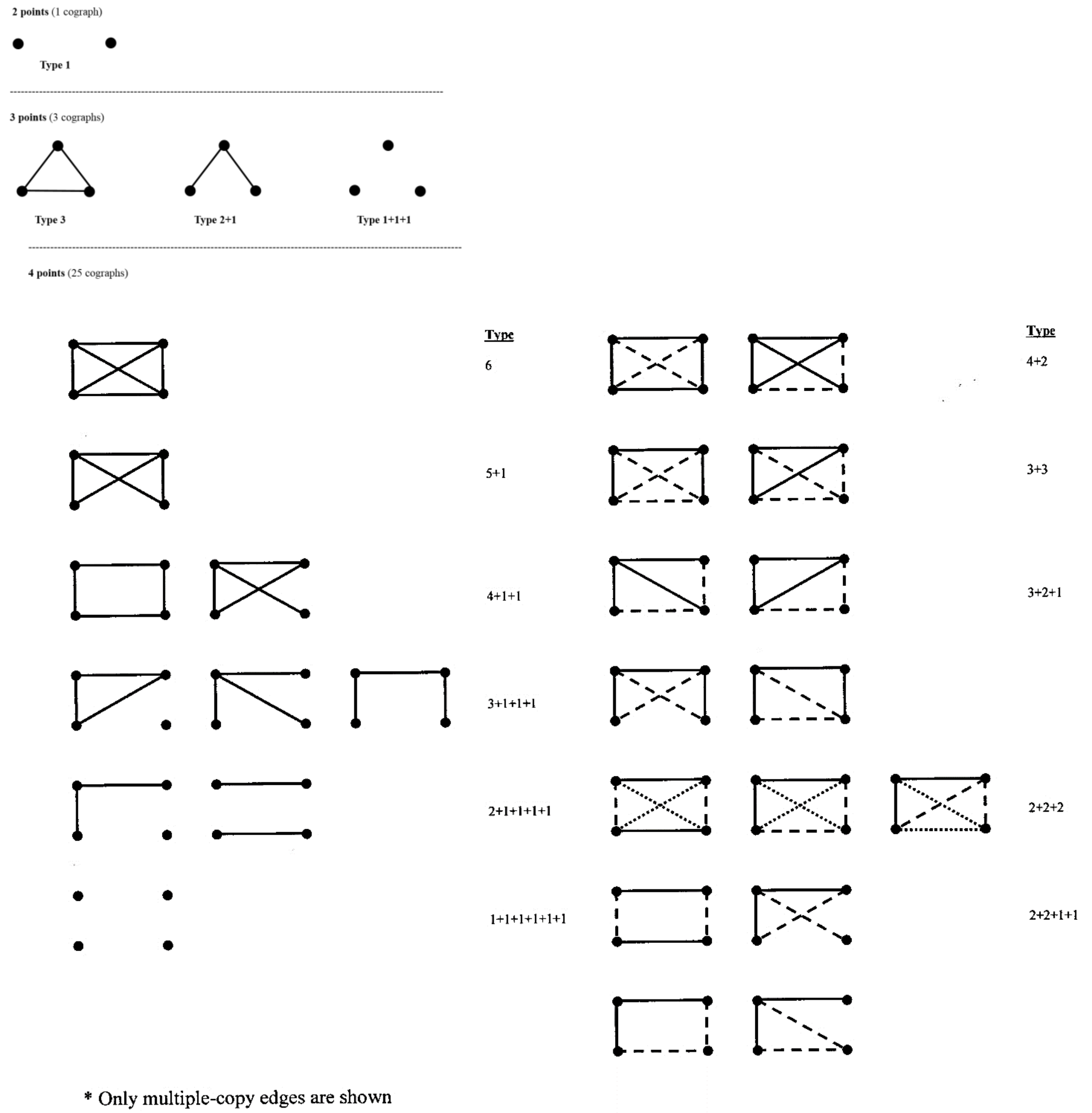

**Fig. 1.2.** Catalogue of cographs with 2, 3, or 4 points.

**Table 1.1. Number of Cographs on *n* Points**

| Points ($n$) | Cographs |
|---|---|
| 2 | 1 |
| 3 | 3 |
| 4 | 25 |
| 5 | 1299 |
| 6 | $1.97 \times 10^{6}$ |
| 7 | $9.43 \times 10^{10}$ |
| 8 | $1.53 \times 10^{17}$ |
| 9 | $1.05 \times 10^{25}$ |

Appreciation how cographs generalize other combinatorial structures comes from considering the ones with just a few edges. Cographs with a single edge coincide with complete graphs $K_n$ (or equivalently with their complements $\overline{K_n}$ having $n$ points and no lines). Viewing a cograph with two edges as being a graph with black and (invisible) white lines shows how it coincides with an ordinary graph (or, again, the complement). Any cograph thus naturally breaks up into a line-disjoint union of graphs, one for each color, on the same set of points. And, as Figure 1.3 illustrates, some cographs with at most three edges (the union of two of which must be planar) correspond to knots.

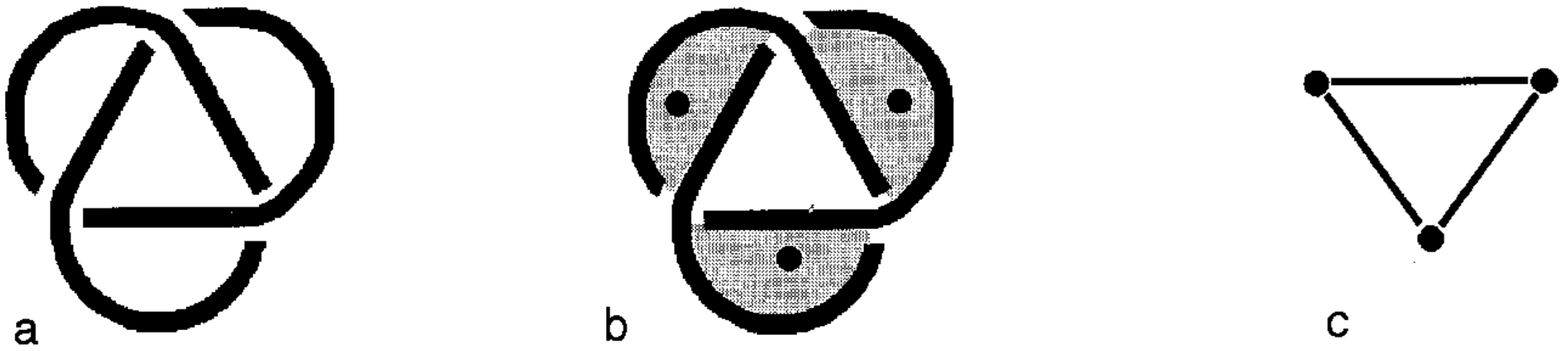

**Fig 1.3.** Converting a trefoil knot into a one-edge cograph:
a) knot; b) shading; c) line-colored graph = cograph.

A knot (a) is a circle mapped into the plane with over- and under-crossings specified. Shade in alternate regions and place graph vertices at their centers (b). Encode crossing information by coloring the lines joining the vertices: blue, if rotating the upper strand of the knot a quarter turn counterclockwise sweeps over two shaded regions; red if over two unshaded regions (Adams [1] 51-5; Kauffman [14] 223-4); and green if the vertices are not adjacent (c: three blue lines).

Figure 1.4 illustrates how families of cographs arise naturally in several areas of mathematics:

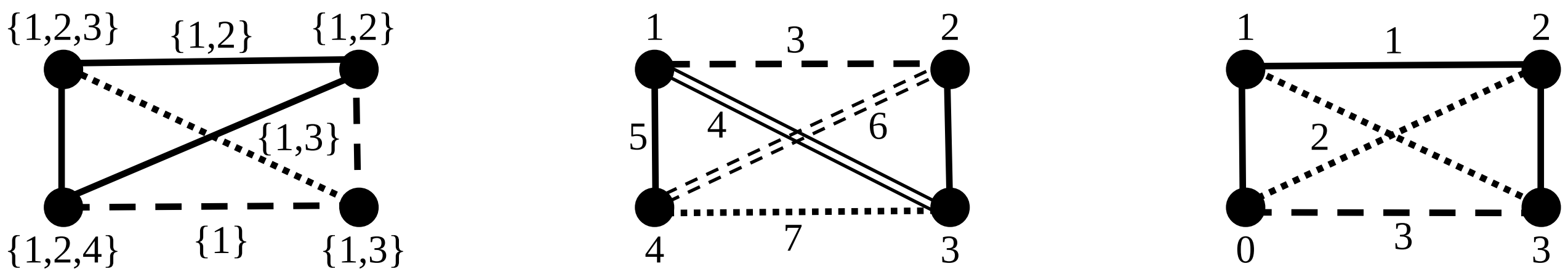

a) Intersection cograph b) Sum cograph c) Distance cograph

**Fig. 1.4.** Some simple mathematical cographs.

In an intersection cograph (a) the points are sets, and the edge $\mathcal{C}(P,Q)$ joining points P and Q is their intersection $P \cap Q$. In algebraic cographs the points are elements of an algebraic structure, and edges the corresponding algebraic product; thus in sum cograph (b) the points are integers, with $\mathcal{C}(P,Q) = P + Q$. Metric structures are yet another source of cographs: in the distance cograph (c) the points are real numbers, with $\mathcal{C}(P,Q) = |P-Q|$.

One theme in this study will be to investigate which cographs can arise in each of these various ways. Figure 1.5 illustrates some configurations that are forbidden to occur:

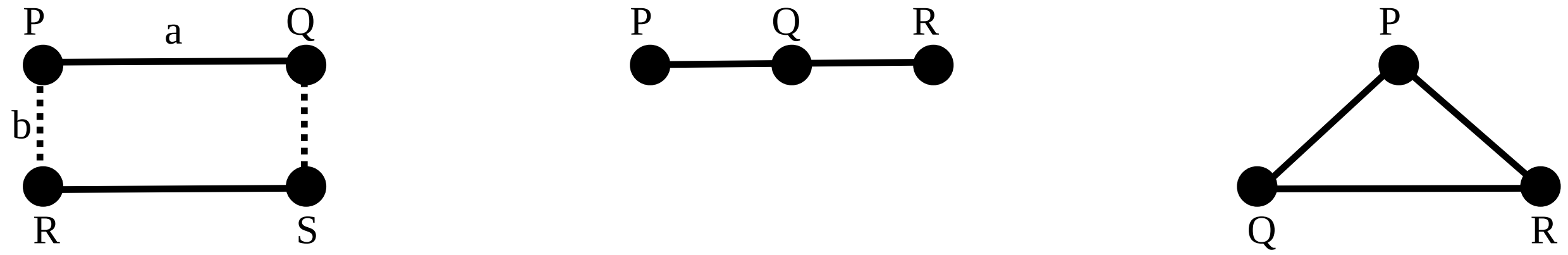

a) In an intersection cograph b) In a sum cograph over an abelian group c) In a distance cograph over the real numbers

**Fig. 1.5**. Forbidden configurations.

For in (a), $a \neq b$, yet $a = a \cap a = (P \cap Q) \cap (R \cap S) = (P \cap R) \cap (Q \cap S) = b \cap b = b$; in (b), $P \neq R$, yet $P+Q = \mathcal{C}(P,Q) = \mathcal{C}(Q,R) = Q+R$ implies $P = R$; and in (c), relabeling points, if necessary, to make $P < Q < R$, $\mathcal{C}(P,R) = R-P = (R-Q) + (Q-P) = \mathcal{C}(Q,R) + \mathcal{C}(P,Q) = 2\,\mathcal{C}(P,R)$, a contradiction.

**Remarks**. 1.1) *Other variants*: One might also define a "union cograph" on sets, with $\mathcal{C}(P,Q) = P \cup Q$; but taking complements and applying DeMorgan's laws, the same cograph may be obtained as an intersection cograph. "Number theoretic cographs," with $\mathcal{C}(P,Q)$ being the gcd or lcm of numbers P and Q, are likewise just variants of intersection cographs.

1.2) *Intersection cographs as algebraic*: An intersection cograph is equivalent to an algebraic cograph of the form $Z_2 \oplus \ldots \oplus Z_2$ indexed by the union of the sets of the cograph. For simply replace each point or edge S by the "characteristic function" $\chi_S$ that is 1 exclusively on coordinates in S, and set intersections then correspond to products in the ring.

1.3) *Potential applications*: This initial article must focus on the mathematics of cographs; but intersection cographs, especially, promise to be useful for modeling real-world phenomena. One might model, for instance, the people participating in social or political interactions as each one equivalent to the set of his own "interests," and describe a process like negotiation as a search through cographs for maximal overlap. Table 1.1, showing the vast numbers of possible cographs (i.e. configurations, alignments, coalitions) on even a few points, might give insight why reaching consensus can sometimes be so difficult. A similar viewpoint appears in branches of philosophy that regard every object as underlying "substance" bearing the set of its "accidents." Another class of applications might let cograph edges represent forces between the points: a crystal, atoms held together by electrostatic forces between each pair; a galaxy, stars held by pairwise gravitational attraction; a poem, words each linked to many others by sound or sense. (See §4.2 for further applications to aesthetics.)

As a first theorem about cographs, one proves that they may be represented in several ways:

**Theorem 1.1**. *Every cograph may be represented as 1) a "fat" intersection cograph. Every finite cograph may also be represented as 2a) an inner product cograph, 2b) a polynomial cograph, and 3) a sufficiently high-dimensional geometric cograph.*

**Proof.** 1) Let $\mathcal{P}$ and $\mathbb{S}$ be collections of sets, where $\mathbb{S}$ is closed under intersections and contains $U = \bigcup_{P \in \mathcal{P}} P$. Define a *fat intersection cograph* $\mathcal{F}$ to have points $\mathcal{P}$, where the edge $\mathcal{F}(P,Q)$ is the smallest set in $\mathbb{S}$ containing $P \cap Q$.

**Fig. 1.6.** Example of a fat intersection cograph. Let $\mathbb{S} = \{\emptyset, \{a,b\}, \{c,d\}, \{a,b,c,d\}\}$. (Note that this contains the configuration of Figure 1.5a, hence is not an intersection cograph.)

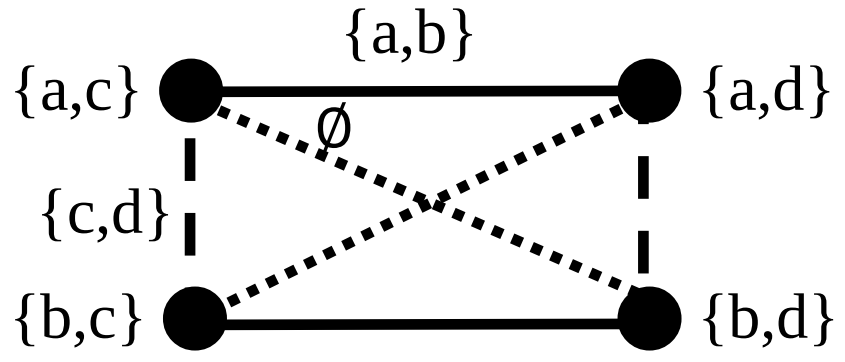

Given then an abstract cograph $\mathcal{C}$, first make the copies of each edge distinguishable by relabeling $\mathcal{C}(P,Q) = e$ as $\mathcal{C}'(P,Q) = e_{\{P,Q\}}$, and then label each point P with the set $P^* = P_o \cup \bigcup_{Q \in \mathcal{C}} \{e_{\{P,Q\}}\}$. Then, for all P and Q, $P^* \cap Q^*$ is precisely $\mathcal{C}'(P,Q)$. Define $U = \bigcup_{P \in \mathcal{C}} P^*$, and define the set $\mathbb{S}$ as:

$$\mathbb{S} = \{\emptyset, U\} \cup \bigcup_{e \in \mathcal{C}} \{e_{\{P,Q\}} : e_{\{P,Q\}} \text{ is a relabeled copy of } e\}$$

Edges then belong to the same set in $\mathbb{S} \sim \{\emptyset, U\}$ if and only if they originated from multiple copies in $\mathcal{C}$, and thus $\mathcal{C}$ is equivalent to the fat intersection cograph with points $P^*$ and list $\mathbb{S}$.

2a) The $n$-point abstract cograph $\mathcal{C}$ will be represented by an $(n-1)$-dimensional *inner product cograph* $\mathcal{I}$--its points from a real inner product space X, the edge between two points defined as their inner (dot) product, $\mathcal{I}(P,Q) = (P,Q)$--making an *arbitrary* choice of distinct real numbers for the edges. Thus the goal is to represent each point $P_i$ ($1 \le i \le n$) of $\mathcal{C}$ by a point $(p_{i1}, p_{i2}, \ldots, p_{i(n-1)})$ in $\mathcal{I}$. As an initial labeling, let $P_1$ be $(1,0,\ldots,0)$, and let $p_{i1} = \mathcal{C}(P_1,P_i)$ for $i = 2,3,\ldots,n$. Proceeding inductively, once $p_{ij}$ has been defined for $i,j < k$, let $p_{kk} = 1$, $p_{ki} = 0$ for $i>k$, and for $i<k$ let $p_{ik} = \mathcal{C}(P_i,P_k) - \sum_{j=1}^{k-1} p_{ij}p_{kj}$. From this construction, the edge label $\mathcal{C}(P_i,P_j)$ coincides with the inner product. $P_1,P_2,\ldots,P_{n-1}$ have distinct labels (since $p_{ij} = 1$ when $i = j$ and 0 when $i<j$), but not necessarily $P_n$. To guarantee it distinct the initial labeling must be modified. At the cost of increasing the dimension, one can simply add an $n$th coordinate, 0 for $P_1,P_2,\ldots,P_{n-1}$ and 1 for $P_n$.

Otherwise, one modifies the initial labeling stepwise: find the lowest index $i$ where $P_i$ coincides with $P_n$, and change $p_{ii}$ from 1 to 2, and $p_{ik}$ to half its value for $k>i$; then repeat inductively to eliminate any introduced new coincidences.

**Remark**. The bound on the dimension is the best possible: It can be shown that the cograph with $n$ points ($n \geq 3$) and a single edge cannot be represented by an inner product cograph of dimension $n$-2.

2b) The $n$-point abstract cograph $\mathcal{C}$ will be represented by a *polynomial cograph* $\mathcal{P}$--its points from Z (or $Z_m$), the edge $\mathcal{P}$(P,Q) defined to equal f(P,Q) where f is a symmetric polynomial in two variables--labeling the points of $\mathcal{C}$ by the integers 1,2,...,$n$, and its edges by *arbitrarily* chosen distinct positive integers. Define first a family of elementary symmetric polynomials $f_{ij}$, $1 \leq i < j \leq n$, by $f_{ij}(x,y) = \prod_{\substack{1 \leq s < t \leq n \\ \{s,t\} \neq \{i,j\}}} [(x-s)^2 + (y-t)^2] \cdot [(x-t)^2 + (y-s)^2]$; thus $f_{ij}(i,j) = f_{ij}(j,i) \neq 0$, but $f_{ij}(x,y) = 0$ otherwise. A symmetric interpolatory polynomial f satisfying f(P,Q) = $\mathcal{C}$(P,Q) for all the chosen values of $\mathcal{C}$(P,Q) will then be an appropriate linear combination, with rational coefficients, of the $f_{ij}$'s. Multiplying f by a constant to clear the denominator will not change the cograph, which is thus over Z (or, since finite, in $Z_m$ for sufficiently large $m$).

3) The goal is to represent the $n$-point abstract cograph $\mathcal{C}$ by a *geometric cograph* $\mathcal{G}$--points from $R^{n-1}$, Euclidean distances for edges. The cograph with just one edge (of length 1) is given by the points $(1/\sqrt{2}, 0,0,...,0)$, $(0, 1/\sqrt{2},0,...,0)$,...,$(0,0,0,...,0, 1/\sqrt{2})$, and (c,c,c,...,c), where $c = (1+ 1/\sqrt{(n+1)})/(n/\sqrt{2})$. To obtain a cograph with more edges one adjusts edge lengths in this figure: any side *independently* can be expanded slightly by moving an endpoint in a way that changes the angles but not the lengths of the other sides. The sole "geometric" constrant is that subpolyhedra remain nondegenerate, i.e. keep nonzero volume, and these volumes can be expressed by determinants (Borsuk [5] 116-120) which are continuous functions of the coordinates, hence remain bounded away from zero within small neighborhoods. Hence for a sufficiently large integer $K$ one may convert the starting cograph to an arbitrary one with $m$ edges, where edge $i$ has length $1+i/K$, for $i$ = 1,2,...,$m$. Since scaling the figure up by a factor of $K$ does not change the cograph, this proof also shows that one can realize an arbitrary finite cograph by a geometric cograph with edges of integer lengths.

**Remarks**. The constraint on dimensions is the best possible: one cannot realize a three-point cograph with one edge--an equilateral triangle--in $R^1$. Note also that geometric representation is not always possible with rational or integer *points*--for instance, there is no equilateral triangle with all three vertices rational in $R^2$. ■

Theorem 1.1 seeks concrete terms like sets for both points and edges to produce a given abstract cograph--the cograph representation problem, that in general seems quite difficult. To conclude this section, here are some results for the easier situation in which the cograph edges are "prelabeled," and the problem is only to find the points.

**Theorem 1.2.** *1) A cograph $\mathcal{C}$ with edges prelabeled as sets may be represented as an intersection cograph if and only if every triangle abc satisfies $a \cap b = b \cap c = a \cap c$ and every quadrilateral abcd satisfies $a \cap c = b \cap d$ (see Figure 1.7).*

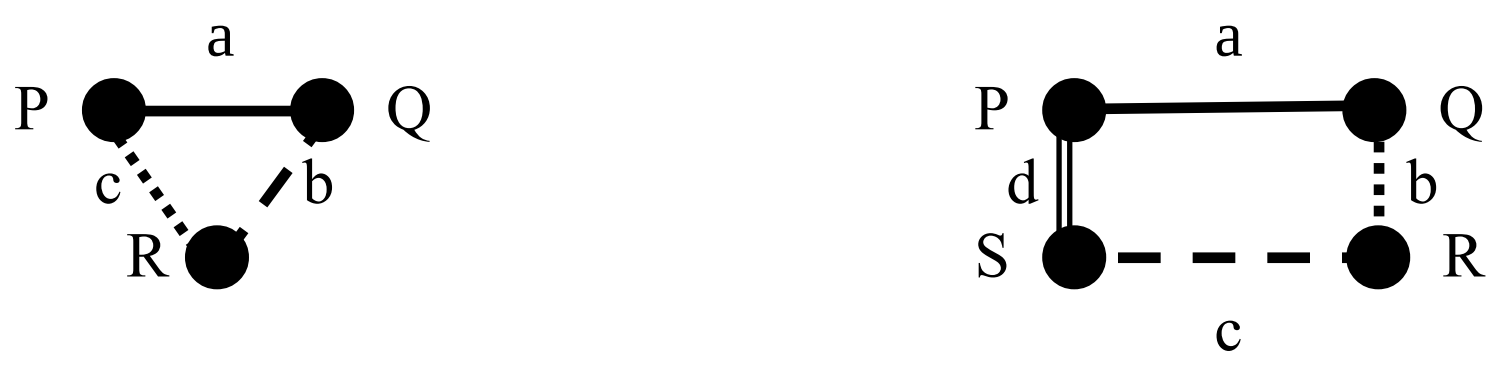

**Fig. 1.7.** A triangle and a quadrilateral.

*Each point P may then be labeled with the set* $P' = \{P_o\} \cup \bigcup_{Q \in \mathcal{C}} \mathcal{C}(P,Q)$ *, where* $P_o$ *is a unique "point name" (i.e.* $P' \cap Q' = \mathcal{C}(P, Q)$ *for all P and Q in* $\mathcal{C}$*).*

*2) A cograph* $\mathcal{C}$ *with edges prelabeled by elements of an abelian group can have point values assigned to make it a sum cograph if and only if (a) all edges incident to each point are distinct, (b) for each triangle of edges abc, a+b-c is divisible by 2, and (c) for each quadrilateral of edges (in clockwise order) abcd, a+c = b+d. If these conditions hold, assign one point P the label L(P), where* $2L(P) = \mathcal{C}(P,Q) + \mathcal{C}(P,R) - \mathcal{C}(Q,R)$*, for any points Q and R, and then label each* $X \neq P$ *by* $L(X) = \mathcal{C}(P,X) - L(P)$*.*

*3) A cograph* $\mathcal{C}$ *with edges prelabeled with a closed and bounded set of positive numbers is a distance cograph on the real line if and only if (a) in every triangle, one edge is the sum of the other two, and (b) the subcograph of Figure 1.8 does not occur. If these conditions hold, choose P and Q with* $\mathcal{C}(P,Q)$ *maximal, assign point P value 0, and point R value* $\mathcal{C}(P,R)$ *for each* $R \neq P$*.*

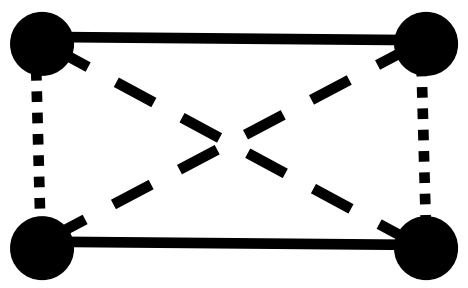

**Fig. 1.8.** The forbidden subcograph.

***Proof.*** For 1) see Proposition 4.1.3 below; the simple proofs of 2) and 3) are omitted. ■

## 2.1. Sum cographs: Elementary properties

A sum cograph, as mentioned above, is a cograph having its points and edges in an abelian group, in which each edge is simply the sum of its endpoints: $\mathcal{C}(P,Q) = P+Q$. In consequence, $Q = \mathcal{C}(P,Q) - P$, and, inductively, the value of each point in any chain of points and edges is determined by the initial point P and the alternating sum of the intervening edges, as follows:

**Proposition 2.1.1.** *In a sum cograph, the points in a chain of points and edges have the form shown in Figure 2.1.1.*

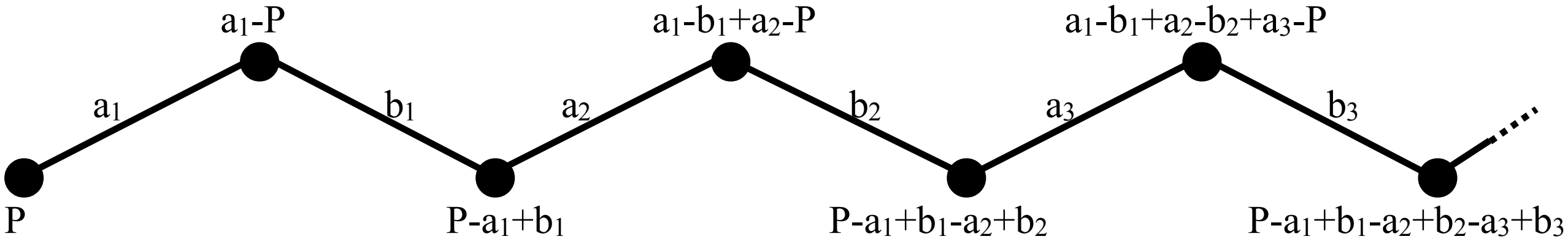

**Fig. 2.1.1**. The points in a chain in a sum cograph.

The elementary properties of sum cographs follow as corollaries:

**Corollary 2.1.1 (V-free)**. *The V-configuration of Figure 2.1.2 is forbidden in a sum cograph: all edges incident at a point must be distinct.*

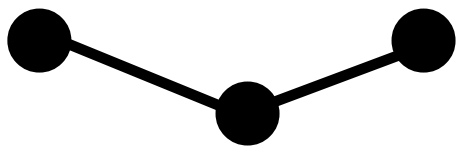

**Fig 2.1.2.** Forbidden configuration.

**Proof**. $a_1 = b_1$ in Proposition 2.1.1 would force the third point to equal P, contrary to the diagram of this Corollary which shows the two points as distinct. ■

**Corollary 2.1.2 (Quadrilateral rule)**. *In any quadrilateral in a sum cograph, (Figure 2.1.3), $a + c = b + d$; thus each edge $d = a - b + c$ is determined by the other three edges.*

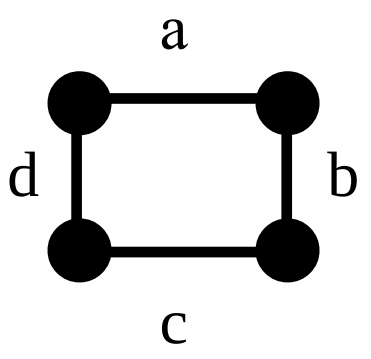

**Fig 2.1.3.** Quadrilateral.

**Proof**: Identifying the fifth point in Proposition 2.1.1 with the first, $-a_1+b_1-a_2+b_2 = 0$. ■

**Corollary 2.1.3 (Alternating pentagon rule).** *In a sum cograph the configuration of Figure 2.1.4 is forbidden.*

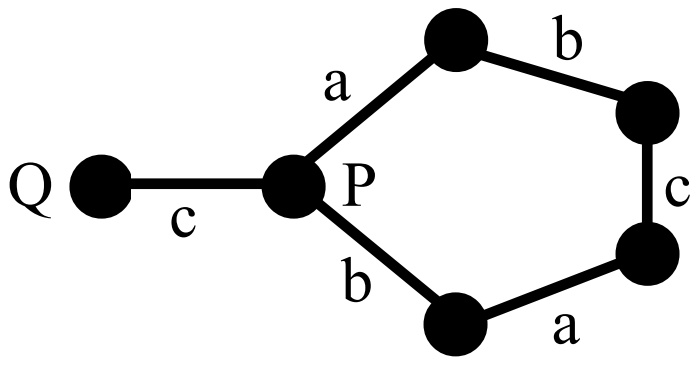

**Fig 2.1.4.** Alternating pentagon.

**Proof**: Identifying the sixth point in Proposition 2.1.1 with the first, it follows that in the diagram of this Corollary, $2P = a - b + c - a + b = c$. But in this sum cograph $P + Q = c$. Hence $2P = c = P+Q$, forcing $P = Q$, which contradicts the diagram showing these points as distinct. ■

The argument obviously generalizes to any alternating 4k+1 – sided polygon (see also Corollary 2.1.6 below).

**Corollary 2.1.4 (Forced hexagon)**. *In a sum cograph, the configuration of Figure 2.1.5, part (a), forces that of part (b).*

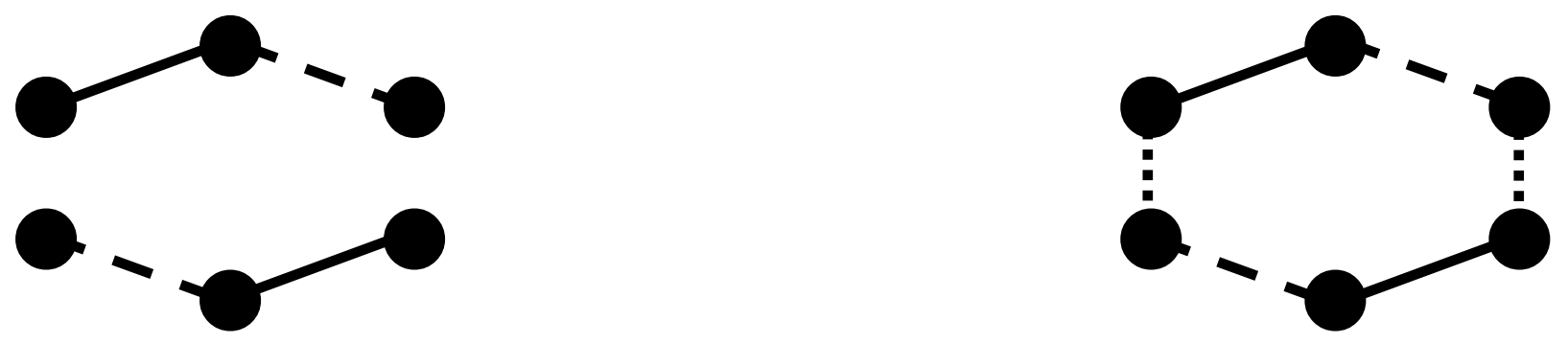

**Fig. 2.1.5.** Configuration a forces configuration b.

**Proof**: The central vertical edge of the hexagon bisects it into two equal three-sided figures, and Corollary 2.1.2, the quadrilateral rule, then forces the two vertical side edges of the hexagon to be equal. ■

Corollary 2.1.4 and Corollary 2.1.1 together yield, incidentally, an alternate proof for Corollary 2.1.3.

**Corollary 2.1.5 (Alternating n-cycle)**. *If a sum cograph contains an alternating cycle of (even) length n, its group must have an element with n/2-torsion.*

**Proof**. In Proposition 2.1.1, denote $a_1 = a_2 = ... = a_{n/2} = a$ and $b_1 = b_2 = ... = b_{n/2} = b$. Then the n+1st point is $P - (n/2)a + (n/2)b$. That the chain closes to a cycle means this point is again P, so $(n/2)(b-a) = 0$, and b-a has n/2-torsion. ■

**Corollary 2.1.6 (2-Torsion)**. *If a sum cograph contains either two identical odd cycles, or the configuration of Figure 2.1.6b, then its group must have 2-torsion.*

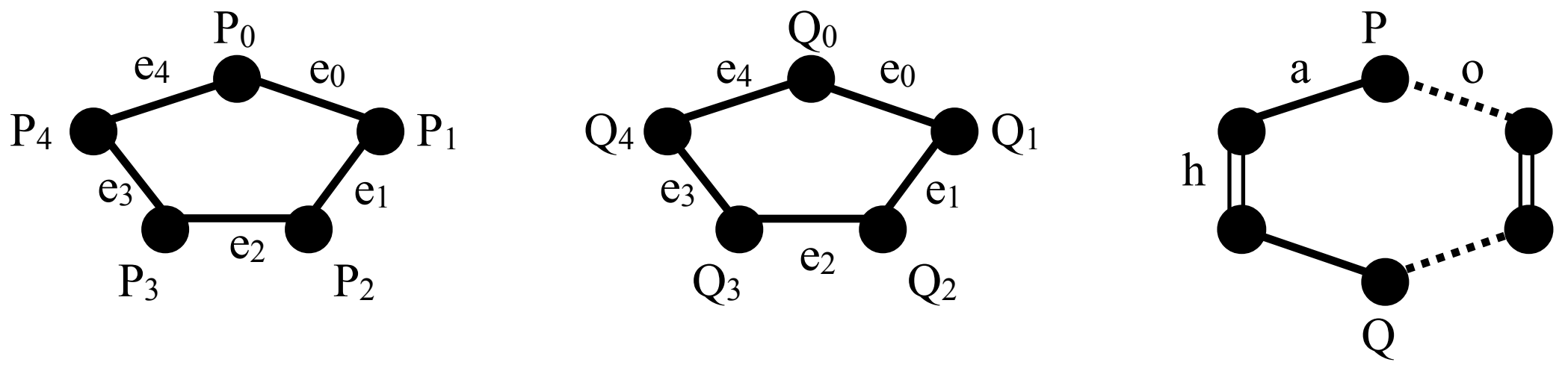

**Fig. 2.1.6.** Configurations forcing 2-torsion

**Proof.** We illustrate the proof in the case that the cograph contains two 5-cycles (Figure 2.1.6(a)). Then $2P_0 = (P_0 + P_1) - (P_1 + P_2) + (P_2 + P_3) - (P_3 + P_4) + (P_4 + P_0) = e_0 - e_1 + e_2 - e_3 + e_4 =$ $(Q_0 + Q_1) - (Q_1 + Q_2) + (Q_2 + Q_3) - (Q_3 + Q_4) + (Q_4 + Q_0) = 2Q_0$, and if the group lacks 2-torsion then $P_0 = Q_0$. Similarly, $P_i = Q_i$ for all i, contradicting the assumption that the two cycles are distinct.

In the case of Figure 2.1.6(b), the edge $b = \mathcal{C}(P,Q) = P + Q$ splits the figure to two quadrilaterals, for which the quadrilateral rule, Corollary 2.1.2, implies $2a - h = b = 2o - h$, forcing the contradiction $a = o$ unless the group has 2-torsion. ■

**Example 2.1 (Shifted n-Diamond).** Consider the configuration of Figure 2.1.7 in a sum cograph:

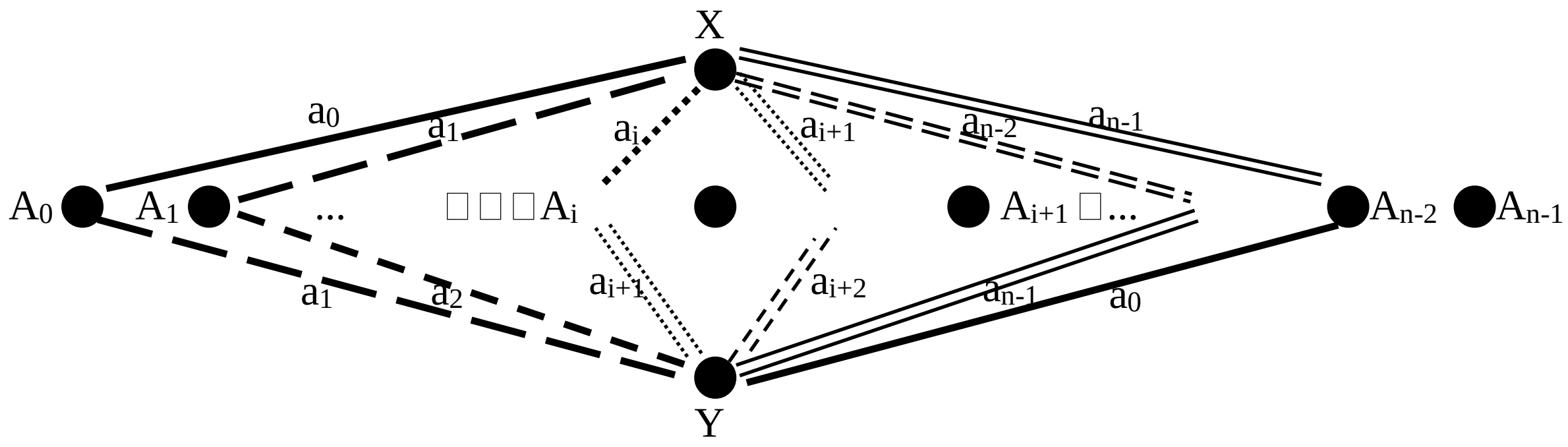

**Fig. 2.1.7.** Shifted n-diamond.

where the edges to the lower point Y are shifted by one as compared to those to the upper point X.

By definition, for all i, $a_i = A_i + X = A_{i-1} + Y$.

By the central quadrilaterals, $a_{i+2} + a_i = 2a_{i+1}$, whence $a_{i+2} - a_{i+1} = a_{i+1} - a_i \equiv \Delta$.

Telescoping, $a_i = (a_i - a_{i-1}) + (a_{i-1} - a_{i-2}) + ... + (a_1 - a_0) + a_0 = i\Delta + a_0$, whence $na_i = na_0$ and $n\Delta = 0$.

Thus all the edges $a_i$ are determined by the choice of $a_0$ and $\Delta$, and similarly for the points,

$A_i = A_0 + i\Delta$ (so $nA_i = nA_0$), $X = a_0 - A_0$, and $Y = X + \Delta$ (so $nX = nY$).

The edges, not shown in the figure, between the $A_i$'s are also determined, and may repeat:

$\mathcal{C}(A_i,A_j) = A_i + A_j = 2A_0 + (i + j)\Delta$ (so $n\mathcal{C}(A_i,A_j) = 2nA_0$),

and so $\mathcal{C}(A_i,A_j) = \mathcal{C}(A_{i'},A_{j'}) \Leftrightarrow (i + j - i' - j')\Delta = 0$.

A cograph's structure is typically determined primarily by its edges: the next proposition explains how a sum cograph can always be renamed to assign an arbitrary value to any given point:

**Proposition 2.1.2.** *Let P be a single arbitrary point in a sum cograph $\mathcal{C}$ and K a single arbitrary value in the underlying abelian group. Then the points and edges of $\mathcal{C}$ may be renamed in a way that assigns P the new name K.*

**Proof.** Assigning to each point X of $\mathcal{C}$ the new name X-P+K, and to each edge e the new name e-2P+2K, produces a sum cograph $\mathcal{C}$' isomorphic to $\mathcal{C}$ with P renamed as K. ■

## 2.2. Sum cographs on six points

The catalogue Figure 2.2.2 below presents the 55 sum cographs on six points (showing only *repeated* edges), and Table 2.1 lists their general and particular numerical representations and other invariants. Since sum cographs on two to five points can naturally be extended to six points by adding isolated points with high numerical values, the figure and table actually present all sum cographs on six *or fewer* points.

The cographs were built up systematically by adding increasing numbers of pairs or triples of equal edges. Hence the finished figure is *complete*. A complementary concern is whether the process might have left duplicates, arriving at the same final structure by distinct paths. But the numerical invariants in the table, counting number of repeated edges and their degrees at each point, suffice to prove that these cographs are all *distinct*.

The construction involved analyzing well over 400 candidate structures. So any tabulation on more than six points will surely require computer assistance.

The cographs are sorted by number of repeated pairs or triples of edges, and whether they require torsion elements. Any (finite) sum cograph can, of course, be realized in a sufficiently large torsion group. But since some simple cograph configurations force torsion (see, e.g., Corollaries 2.1.5 and 2.1.6 and Example 2.1 in the preceding section), torsion-free solutions are likely increasingly rare in large cographs, so seem worthy of special note.

In the table, note that for the point degrees and edge multiplicities only *multiple* edges are counted. In the general P,Q,R,… representations, the forms like P+Q-R or 3Q-2P are explicable by the fact that increasing the value of each variable by 1 yields another valid representation. Thus if $X = \Sigma a_i X_i$, then $X+1 = \Sigma a_i(X_i+1)$, and by subtraction, $\Sigma a_i = 1$. Also useful in finding convenient representations is that fact that if $\{X_i\}$ is a representation, so is $\{M-X_i\}$, for any fixed value M.

### Construction of Figure 2.2.2

A candidate cograph structure yields a system of linear equations, which can have one of four possible outcomes when one tries to solve it:

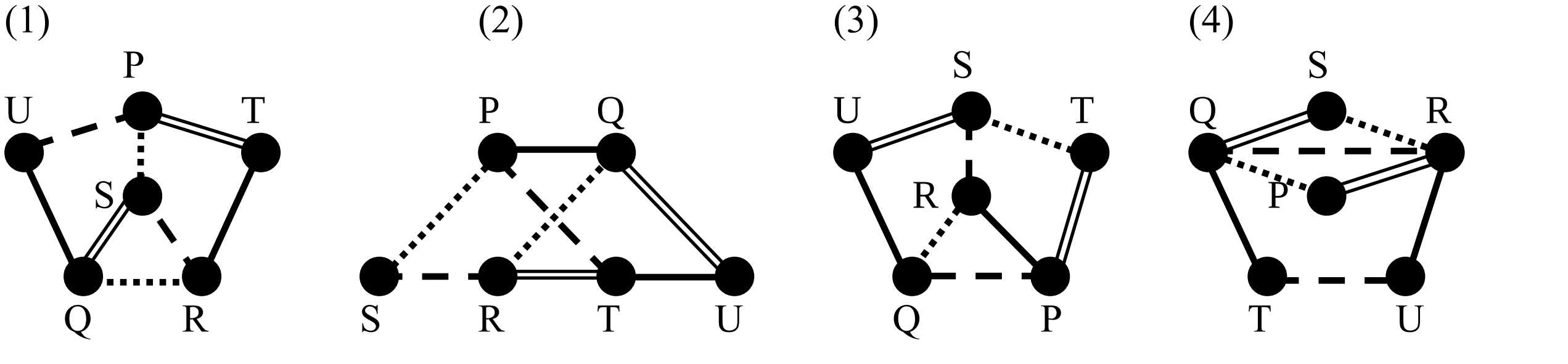

Solvable in integers — Solvable, with torsion — Unsolvable — Forces additional edges

**Figure 2.2.1.** Cographs illustrating the four possible outcomes

(1) Solvable in integers (torsion-free): The candidate cograph (1), for example, yields the system:

| | |
|---|---|
| Q + U = R + T | (from edge ▬▬▬ ) |
| P + U = R + S | (from ▬ ▬ ) |
| P + S = Q + R | (from ▪▪▪▪▪▪) |
| P + T = Q + S | (from ═══) |

which has the general solution S = Q+R-P, T = R+2Q-2P, U = Q+2R-2P for any choice of P, Q, and R. Choosing P = 0, Q = 1, and R = 3 yields a solution (P,Q,R,S,T,U) = (0,1,3,4,5,7) in conveniently small non-negative integers, and one may check that no extraneous further equalities arise from or among the other pairwise sums P+Q, P+R, Q+T, R+U, S+T, S+U, or T+U. So this is torsion-free sum cograph M of the figure.

(2) Solvable, but with torsion: Cograph (2) yields the system:

[P + Q = T + U, P + T = R + S, P + S = Q + R, Q + U = R + T]

with general solution S = Q+R-P, T = Q+2R-2P, U = 3R-2P, with 5P = 5R. This last equation forces torsion, and a solution in $Z_{15}$ is (P,Q,R,S,T,U) = (0,1,3,4,7,9), sum cograph i in the figure.

(3) Unsolvable (torsion = 1): Cograph (3) yields the system:

[P + Q = R + S, P + R = Q + U, P + T = S + U, Q + R = S + T]

which on solving implies P = T (i.e. torsion 1). But this is not allowed, since we insist that the points of a cograph be distinct, hence it is not solvable as a sum cograph. The 400+ systems analyzed for the figure had 46 such failures, 45 of them attributable (as here) to the single alternating pentagon configuration of Corollary 2.1.3 above.

(4) Forces additional edges: Cograph (4) yields the system:

[P + Q = R + S, P + R = Q + S, Q + R = T + U, Q + T = R + U]

with general solution: S = P+R-Q, U = Q+T-R and torsion 2Q = 2R = 2T = 2U, 2P = 2S. But this solution then implies three further equations: P + T = U + S, P + U = S + T, and Q + U = R + T, that is, three new edge pairs will have to be added to the cograph for it to be the sum cograph. The 400+ systems analyzed for the figure had around 130 such failures, all but one of them attributable (as in the first two equations here) to the single "forced hexagon" configuration of Corollary 2.1.4 above.

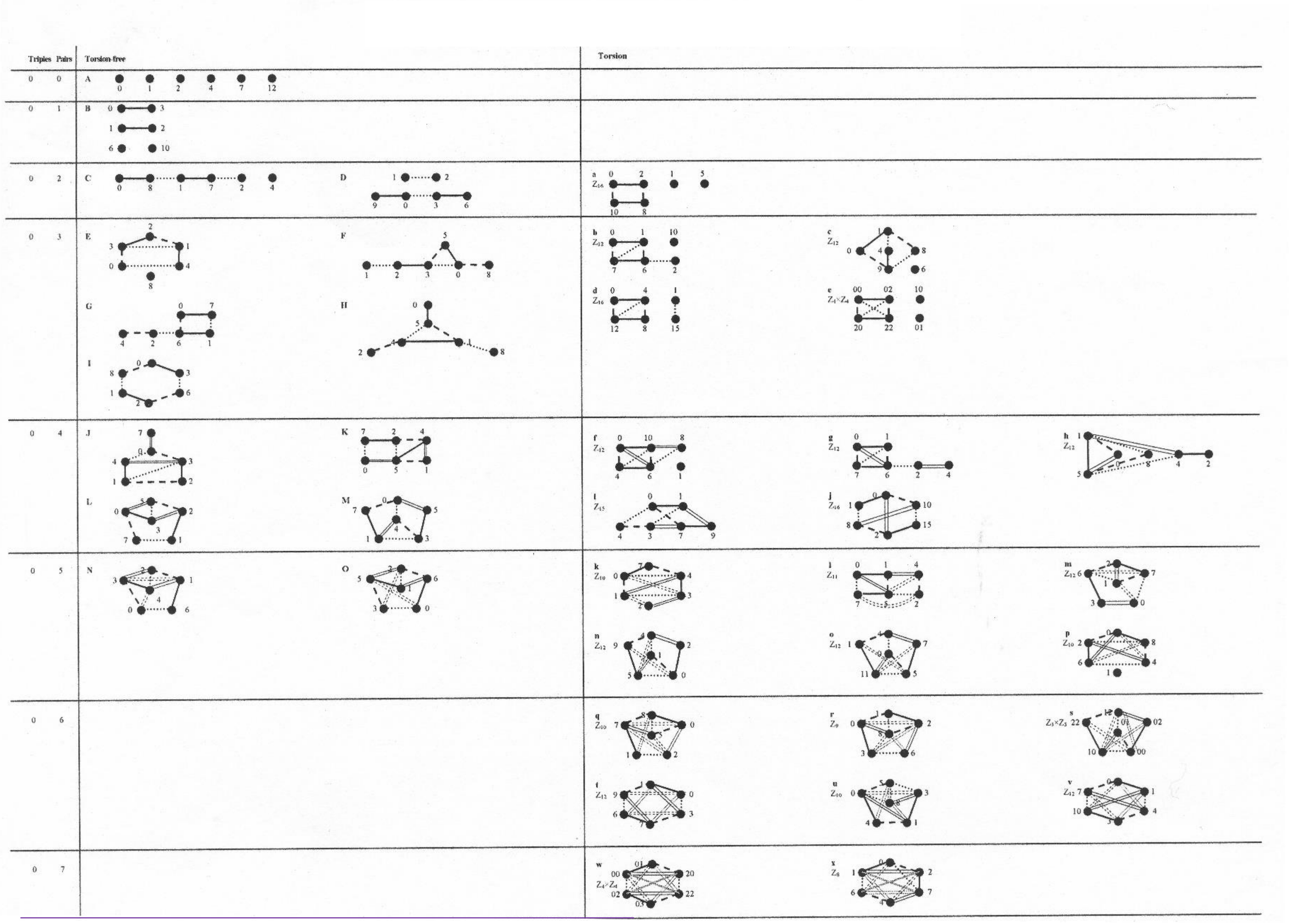

**Fig. 2.2.2.** Catalogue of six-point sum cographs. (continued next page)

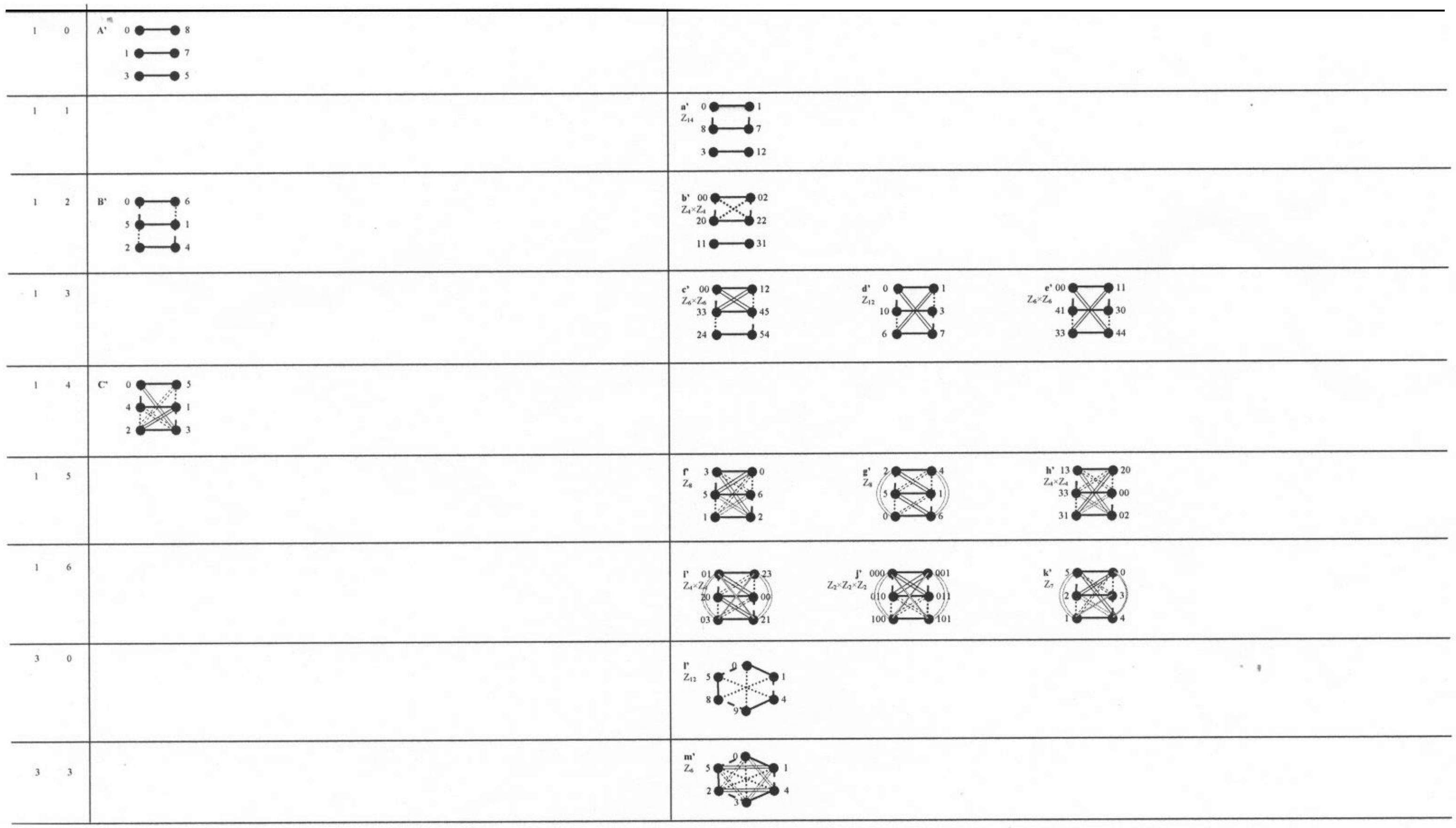

**Fig. 2.2.2.** Catalogue of six-point sum cographs. (continued from previous page)

| Triples | Pairs | Co-graph | Numerical | Group | Point degrees | Edge multi-plicities | P | Q | R | S | T | U | Torsion |
|---|---|---|---|---|---|---|---|---|---|---|---|---|---|
| TORSION-FREE | | | | | | | | | | | | | |
| 0 | 0 | **A** | 0 1 2 4 7 12 | Z | 0 0 0 0 0 0 | 0 | P | Q | R | S | T | U | |
| 0 | 1 | **B** | 0 1 2 3 6 10 | Z | 0 0 1 1 1 1 | 2 | P | Q | R | P+Q-R | T | U | |
| 0 | 2 | **C** | 0 1 2 4 7 8 | Z | 0 1 1 2 2 2 | 2 2 | P | Q | 2Q-P | S | T | Q+T-P | |
| 0 | 2 | **D** | 0 1 2 3 6 9 | Z | 1 1 1 1 2 2 | 2 2 | P | Q | R | Q+R-P | T | Q+R+T-2P | |
| 0 | 3 | **E** | 0 1 2 3 4 8 | Z | 0 2 2 2 3 3 | 2 2 2 | P | Q | 2Q-P | 3Q-2P | 4Q-3P | U | |
| 0 | 3 | **F** | 0 1 2 3 5 8 | Z | 1 1 2 2 3 3 | 2 2 2 | P | Q | R | Q+R-P | Q+2R-2P | 2Q+3R-4P | |
| 0 | 3 | **G** | 0 1 2 4 6 7 | Z | 1 2 2 2 2 3 | 2 2 2 | P | Q | 2Q-P | S | 2Q+S-2P | 3Q+S-3P | |
| 0 | 3 | **H** | 0 1 2 4 5 8 | Z | 1 1 1 3 3 3 | 2 2 2 | P | Q | 2Q-P | S | Q+S-P | 2S-P | |
| 0 | 3 | **I** | 0 1 2 3 6 8 | Z | 2 2 2 2 2 2 | 2 2 2 | P | Q | R | Q+R-P | T | R+T-P | |
| 0 | 4 | **J** | 0 1 2 3 4 7 | Z | 1 2 3 3 3 4 | 2 2 2 2 | P | Q | 2Q-P | 3Q-2P | 4Q-3P | 7Q-6P | |
| 0 | 4 | **K** | 0 1 2 4 5 7 | Z | 2 2 2 3 3 4 | 2 2 2 2 | P | Q | 2Q-P | 4Q-3P | 5Q-4P | 7Q-6P | |
| 0 | 4 | **L** | 0 1 2 3 5 7 | Z | 2 2 3 3 3 3 | 2 2 2 2 | P | Q | R | Q+R-P | Q+2R-2P | Q+3R-3P | |
| 0 | 4 | **M** | 0 1 3 4 5 7 | Z | 2 2 3 3 3 3 | 2 2 2 2 | P | Q | R | Q+R-P | R+2Q-2P | Q+2R-2P | |
| 0 | 5 | **N** | 0 1 2 3 4 6 | Z | 2 3 3 4 4 4 | 2 2 2 2 2 | P | Q | 2Q-P | 3Q-2P | 4Q-3P | 6Q-5P | |
| 0 | 5 | **O** | 0 1 2 3 5 6 | Z | 3 3 3 3 4 4 | 2 2 2 2 2 | P | Q | 2Q-P | 3Q-2P | 5Q-4P | 6Q-5P | |
| 1 | 0 | **A’** | 0 1 3 5 7 8 | Z | 1 1 1 1 1 1 | 3 | P | Q | R | S | R+S-Q | R+S-P | |
| 1 | 2 | **B’** | 0 1 2 4 5 6 | Z | 2 2 2 2 3 3 | 2 2 3 | P | Q | 2Q-P | S | Q+S-P | S+2Q-2P | |
| 1 | 4 | **C’** | 0 1 2 3 4 5 | Z | 3 3 4 4 4 4 | 2 2 2 2 3 | P | Q | 2Q-P | 3Q-2P | 4Q-3P | 5Q-4P | |

**Table 2.1.** Six-point sum cographs: Representations and invariants.

(continued next page)

| Triples | Pairs | Co-graph | Numerical | Group | Point degrees | Edge multi-plicities | P | Q | R | S | T | U | Torsion |
|---|---|---|---|---|---|---|---|---|---|---|---|---|---|
| TORSION | | | | | | | | | | | | | |
| 0 | 2 | **a** | 0 1 2 5 8 10 | $Z_{16}$ | 0 0 2 2 2 2 | 2 2 | P | Q | R | S | T | R+T-P | 2P=2T |
| 0 | 3 | **b** | 0 1 2 6 7 10 | $Z_{12}$ | 0 1 2 3 3 3 | 2 2 2 | P | Q | 2Q-P | S | Q+S-P | U | 2P=2S |
| 0 | 3 | **c** | 0 1 4 6 8 9 | $Z_{12}$ | 0 2 2 2 3 3 | 2 2 2 | P | Q | 2T-P | S | T | Q+T-P | 3P=3T |
| 0 | 3 | **d** | 0 1 4 8 12 15 | $Z_{16}$ | 1 1 2 2 3 3 | 2 2 2 | P | Q | R | S | R+S-P | 2R+S-P-Q | 2P=2S |
| 0 | 3 | **e** | 00 01 02 10 20 22 | $Z_4$x$Z_4$ | 0 0 3 3 3 3 | 2 2 2 | P | Q | R | S | T | R+T-P | 2P=2R=2T |
| 0 | 4 | **f** | 0 1 4 6 8 10 | $Z_{12}$ | 0 2 3 3 4 4 | 2 2 2 2 | 3T-2S | Q | 2S-T | S | T | 2T-S | 6S=6T |
| 0 | 4 | **g** | 0 1 2 4 6 7 | $Z_{12}$ | 1 2 3 3 3 4 | 2 2 2 2 | P | Q | 2Q-P | S | 2Q+S-2P | 3Q+S-3P | 6P=4Q+2S |
| 0 | 4 | **h** | 0 1 2 4 5 8 | $Z_{12}$ | 1 2 2 2 4 4 | 2 2 2 2 | P | Q | 2Q-P | S | Q+S-P | 2S-P | 3P=3S |
| 0 | 4 | **i** | 0 1 3 4 7 9 | $Z_{15}$ | 2 2 3 3 3 3 | 2 2 2 2 | P | Q | R | Q+R-P | Q+2R-2P | 3R-2P | 5P=5R |
| 0 | 4 | **j** | 0 1 2 8 10 15 | $Z_{16}$ | 2 2 3 3 3 3 | 2 2 2 2 | P | Q | R | S | P+R-S | Q+2S-P-R | 2P=2S |
| 0 | 5 | **k** | 0 1 2 3 4 7 | $Z_{10}$ | 2 2 4 4 4 4 | 2 2 2 2 2 | P | Q | 2Q-P | 3Q-2P | 4Q-3P | 4P-3Q | 10P=10Q |
| 0 | 5 | **l** | 0 1 2 4 5 7 | $Z_{11}$ | 3 3 3 3 3 5 | 2 2 2 2 2 | P | Q | 2Q-P | 4Q-3P | 5Q-4P | 7Q-6P | 11P=11Q |
| 0 | 5 | **m** | 0 1 2 3 6 7 | $Z_{12}$ | 2 3 3 4 4 4 | 2 2 2 2 2 | P | Q | 2Q-P | 3Q-2P | T | P+Q-T | 2P=2T |
| 0 | 5 | **n** | 0 1 2 4 5 9 | $Z_{12}$ | 2 3 3 4 4 4 | 2 2 2 2 2 | P | Q | 2Q-P | S | Q+S-P | 2S+Q-2P | 3P=3S |
| 0 | 5 | **o** | 0 1 4 5 7 11 | $Z_{12}$ | 3 3 3 3 4 4 | 2 2 2 2 2 | P | Q | 3P+Q-3T | 2P-T | T | 2P+Q-2T | 2Q=2T |
| 0 | 5 | **p** | 0 1 2 4 6 8 | $Z_{10}$ | 0 4 4 4 4 4 | 2 2 2 2 2 | P | Q | 3P-2S | S | 2P-S | 2S-P | 5P=5S |
| 0 | 6 | **q** | 0 1 2 5 7 8 | $Z_{10}$ | 3 4 4 4 4 5 | 2 2 2 2 2 2 | P | Q | 2Q-P | 6P-5Q | 4P-3Q | 3P-2Q | 10P=10Q |
| 0 | 6 | **r** | 0 1 2 3 6 8 | $Z_9$ | 3 3 4 4 5 5 | 2 2 2 2 2 2 | 4S-3Q | Q | 5Q-4S | S | 3Q-2S | 2Q-S | 9Q=9S |
| 0 | 6 | **s** | 00 01 02 10 12 22 | $Z_3$x$Z_3$ | 3 3 3 5 5 5 | 2 2 2 2 2 2 | P | 2R-P | R | S | R+S-P | R+2S-2P | 3P=3R=3S |
| 0 | 6 | **t** | 0 1 3 6 7 9 | $Z_{12}$ | 4 4 4 4 4 4 | 2 2 2 2 2 2 | P | Q | R | S | Q+S-P | R+S-P | 2P=2S |
| 0 | 6 | **u** | 0 1 3 4 5 8 | $Z_{10}$ | 3 4 4 4 4 5 | 2 2 2 2 2 2 | P | Q | 3Q-2P | 7P-6Q | 5Q-4P | 3P-2Q | 10P=10Q |
| 0 | 6 | **v** | 0 1 3 4 7 10 | $Z_{12}$ | 4 4 4 4 4 4 | 2 2 2 2 2 2 | P | Q | R | Q+R-P | Q+2R-2P | Q+3R-3P | 4P=4R |
| 0 | 7 | **w** | 00 01 02 03 20 22 | $Z_4$x$Z_4$ | 4 4 5 5 5 5 | 2 2 2 2 2 2 2 | P | Q | R | Q+R-P | T | R+T-P | 2P=2R=2T |
| 0 | 7 | **x** | 0 1 2 4 6 7 | $Z_8$ | 4 4 5 5 5 5 | 2 2 2 2 2 2 2 | P | Q | R | 2R-P | 3R-2P | P+Q-R | 4P=4R |
| 1 | 1 | **a'** | 0 1 3 7 8 12 | $Z_{14}$ | 1 1 2 2 2 2 | 2 3 | P | Q | R | S | Q+S-P | P+Q-R | 2P=2S |
| 1 | 2 | **b'** | 00 02 11 20 22 31 | $Z_4$x$Z_4$ | 1 1 3 3 3 3 | 2 2 3 | P | Q | R | S | Q+S-P | P+Q-R | 2P=2Q=2S |
| 1 | 3 | **c'** | 00 12 24 33 45 54 | $Z_6$x$Z_6$ | 2 2 3 3 4 4 | 2 2 2 3 | P | Q | 2Q-P | S | Q+S-P | 2P-Q | 2P=2S |
| 1 | 3 | **d'** | 0 1 3 6 7 10 | $Z_{12}$ | 3 3 3 3 3 3 | 2 2 2 3 | P | Q | R | 2R-P | Q+2R-2P | P+Q-R | 2P=2S |
| 1 | 3 | **e'** | 00 11 30 33 41 44 | $Z_6$x$Z_6$ | 3 3 3 3 3 3 | 2 2 2 3 | P | Q | R | S | Q+R-P | Q+S-P | 2P=2R=2S |
| 1 | 5 | **f'** | 0 1 2 3 5 6 | $Z_8$ | 4 4 4 4 5 5 | 2 2 2 2 2 3 | P | Q | 2Q-P | 3Q-2P | 4P-3Q | 3P-2Q | 8P=8Q |
| 1 | 5 | **g'** | 0 1 2 4 5 6 | $Z_8$ | 4 4 4 4 5 5 | 2 2 2 2 2 3 | P | Q | P+2Q | S | P+Q-S | 2Q-S | 2P=2S |
| 1 | 5 | **h'** | 00 02 13 20 31 33 | $Z_4$x$Z_4$ | 4 4 4 4 5 5 | 2 2 2 2 2 3 | P | 3S-2R | R | S | P+R-S | 2S-R | 2P=2S,4R=4S |
| 1 | 6 | **i'** | 00 01 03 20 21 23 | $Z_4$x$Z_4$ | 5 5 5 5 5 5 | 2 2 2 2 2 2 3 | P | Q | 3Q-2P | P+Q-T | T | 2P+T-2Q | 4P=4Q,2Q=2T |
| 1 | 6 | **j'** | 000 001 010 011 100 101 | $Z_2$x$Z_2$x$Z_2$ | 5 5 5 5 5 5 | 2 2 2 2 2 2 3 | P | Q | R | Q+R-P | T | Q+T-P | 2P=2Q=2R=2T |
| 1 | 6 | **k'** | 0 1 2 3 4 5 | $Z_7$ | 5 5 5 5 5 5 | 2 2 2 2 2 2 3 | P | 3P-2S | 5P-4S | S | 2P-S | 4P-3S | 7P=7S |
| 3 | 0 | **l'** | 0 1 4 5 8 9 | $Z_{12}$ | 3 3 3 3 3 3 | 3 3 3 | P | Q | R | 2P+Q-2R | 2P-R | P+Q-R | 3P=3R |
| 3 | 3 | **m'** | 0 1 2 3 4 5 | $Z_6$ | 5 5 5 5 5 5 | 2 2 2 3 3 3 | P | Q | 5P-4Q | 4P-3Q | 3P-2Q | 2P-Q | 6P=6Q |

**Table 2.1.** Six-point sum cographs: Representations and invariants.

(continued from previous page)

## 2.3. Fibonacci wheels

An alternating cycle of n points in a sum cograph--points $P_0,P_1,...,P_{n-1}$, edges $\mathcal{C}(P_0,P_1) = \mathcal{C}(P_2,P_3) = ... = a$, $\mathcal{C}(P_1,P_2) = \mathcal{C}(P_3,P_4) = ... = \mathcal{C}(P_{n-1},P_0) = b$, n even--forces n/2 torsion. Fibonacci wheels are a class of sum cograph configurations that can compel far higher torsion.

**Definition**: A *Fibonacci wheel* in a sum cograph consists of a central hub point O = 0, and n points $P_0 = a$, $P_1 = b$ satisfying $P_{i+2} = P_i + P_{i+1}$ for all $i \geq 0$. (The issue is thus how the wheel *closes up,* so that $P_n = P_0 = a$, $P_{n+1} = P_1 = b$, and thereafter $P_{i+n} = P_i$ for all $i \geq n$.) The n edges $\mathcal{C}(O,P_i) = 0 + P_i = P_i$, the *spokes*, repeat the points. By the Fibonacci property $P_{i+2} = P_i + P_{i+1}$, the wheel's *rim* terms $\mathcal{C}(P_i,P_{i+1}) = P_i + P_{i+1} = P_{i+2}$ also match the spokes, and also close up.

**Example**: Figure 2.3.1 is a Fibonacci wheel with seven spokes.

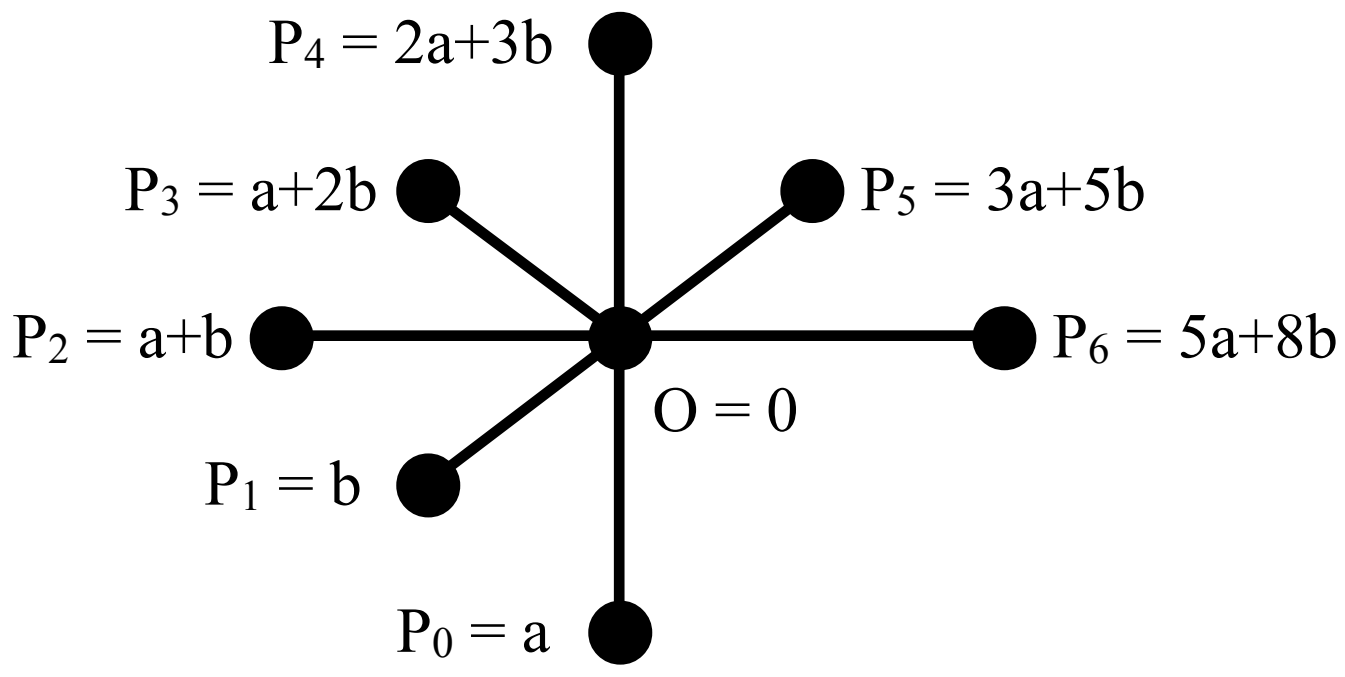

**Fig. 2.3.1.** A Fibonacci wheel with seven spokes.

For this wheel to close up one must have $P_5 + P_6 = P_0$ (namely 8a + 13b = a) and $P_6 + P_0 = P_1$ (namely 6a + 8b = b); in summary, one seeks an abelian group containing elements a and b nontrivially satisfying the pair of equations 7a + 13b = 0 and 6a +7b = 0. Subtracting the second equation from the first yields a = -6b, which, substituted into the first, implies 29b = 0. Similarly, subtracting the first from twice the second yields b = -5a, which, substituted into the second, implies 29a = 0. Thus, working in $Z_{29}$, one may choose a = -1, whence b = 5, and by Fibonacci additivity the wheel terms are $P_0 = -1$, $P_1 = 5$, $P_2 = 4$, $P_3 = 9$, $P_4 = 13$, $P_5 = 22$, $P_6 = 35 = 6$, then closing up to $P_7 = 28 = -1 = P_0$ and $P_8 = 5 = P_1$. Since O = 0, the spokes $\mathcal{C}(P_i,O)$ trivially repeat the points $P_i$, and the rim terms $\mathcal{C}(P_i,P_{i+1}) = P_i + P_{i+1} = P_{i+2}$ also match the spokes, and also close up.

To generalize, let {Fi} be the Fibonacci integer sequence 0,1,1,2,3,5,8,13,21,... defined by $F_0 = 0$, $F_1 = 1$, and $F_{i+2} = F_i + F_{i+1}$ for all $i \geq 0$. Then a Fibonacci wheel with n spokes contains points O = 0, $P_0 = a$, $P_1 = b$, and $P_i = F_{i-1}\, a + F_i\, b$ for $i \geq 1$, and the equations for it to close up are $(F_{n-1} -1)a +F_n\, b = 0$ and $(F_{n-2} +1)a + (F_{n-1} -1)b = 0$. Let $d = \gcd (F_{n-2} +1, F_{n-1} -1)$ be the greatest common divisor of the coefficients of a and b in the second equation. Then, by Fibonacci sequence additivity, it is also the gcd of the coefficients of the first equation.

Factoring out d, the closing-up equations thus become $\mu da + (\mu+\nu)db = 0$ and $\nu da + \mu db = 0$, where the integers $\mu = (F_{n-1} -1)/d$ and $\nu = (F_{n-2} +1)/d$ are relatively prime. Subtracting μ times the second equation from ν times the first, one finds tb = 0, where $t = d(\nu^2 + \mu\nu - \mu^2)$; similarly, ta = 0. Thus the group elements a and b have torsion t.

Applying Cassini's identity $F_nF_{n-2} - (F_{n-1})^2 = (-1)^{n-1}$, the integer $dt = (\nu d)^2 + \mu d \cdot \nu d - (\mu d)^2$ (the negative determinant of the coefficients of the system of equations) reduces to $dt = F_{n-1} + F_{n+1} - 1 - (-1)^n$, so $t = [F_{n-1} + F_{n+1} -1 - (-1)^n]/d$.

Since $\mu$ and $\nu$ are relatively prime, $\mu + \nu = F_n/d$ and $\nu$ are also relatively prime, and the Euclidean algorithm guarantees that integers r and s exist (not uniquely determined) so that $r\mu + s(\mu + \nu) = 1$. Adding s times the first closing-up equation to r times the second then gives $0 = (s\mu + r\nu)da + [r\mu + s(\mu + \nu)]db = (s\mu + r\nu)da + db$, yielding the relation $db = -(s\mu + r\nu)da$, that is, $db = ha$ where $h = -(s\mu + r\nu)d$.

To evaluate d one needs an index-shifting lemma:

**Lemma 2.3.1**: *$\gcd(F_{i-1} + F_{j+1}, F_i - F_j) = \gcd(F_{i-3} + F_{j+3}, F_{i-2} - F_{j+2})$.*

**Proof**: $\gcd(F_{i-1} + F_{j+1}, F_i - F_j) = \gcd(F_{i-1} + F_{j+1}, F_{i-2} + F_{i-1} + F_{j+1} - F_{j+2}) = \gcd(F_{i-1} + F_{j+1}, F_{i-2} - F_{j+2}) = \gcd(F_{i-3} + F_{i-2} - F_{j+2} + F_{j+3}, F_{i-2} - F_{j+2}) = \gcd(F_{i-3} + F_{j+3}, F_{i-2} - F_{j+2})$, where the first and third equalities follow from the inductive definition of the Fibonacci sequence, and the second and fourth from gcd properties, e.g. that $\gcd(x,y) = \gcd(x,x+y)$. ■

Evaluating d, t, r, s, and h then splits into four cases, as n (even) satisfies $n \equiv 0$ or 2 mod 4, or n (odd) satisfies either $n \equiv 1, 5, 7$, or 11, or else $n \equiv 3$ or 9 mod 12:

**Proposition 2.3.1**: *The values of d, t, r, s, and h are as in the table:*

| **n** | **d** $\gcd(F_{n-2}+1, F_{n-1}-1)$ | **t** $ta = tb = 0$ $[L_n-1-(-1)^n]/d$ | **r** $r\mu+s(\mu+\nu)$ | **s** $=1$ | **h** $db = ha$ $-(s\mu+r\nu)d$ |
|---|---|---|---|---|---|
| **4k** | $F_{n/2}$ | $5F_{n/2}$ | $F_{n/2-1}$ | $-F_{n/2-2}$ | $-2F_{n/2}$ |
| **4k+2** | $L_{n/2}$ | $L_{n/2}$ | $F_{n/2-2}$ | $-F_{n/2-3}$ | $-L_{n/2}$ |
| **12k+c c = 1,5,7, or 11** | 1 | $L_n$ | $(F_{n-2} -1) /2$ * | $(1-F_{n-3}) /2$ | $(1-L_{n-2})/2$ * |
| **12k+c c=3 or 9** | 2 | $L_n/2$ | $F_{n-2} -1$ | $1-F_{n-3}$ | $1-L_{n-2}$ |

**Table 2.2.** Values of d, t, r, s, and h. (* = half an integer for c = 5 or 11)

Here $\mu = (F_{n-1} -1)/d$, $\nu = (F_{n-2} +1)/d$, so $\mu + \nu = F_n/d$, $t = [F_{n-1} + F_{n+1} -1 - (-1)^n]/d = [L_n -1 - (-1)^n]/d$, and $h = -(s\mu + r\nu)d$ , where the $L_n$'s are the Lucas numbers 1,3,4,7,11,18,... defined $L_1 = 1$, $L_2 = 3$, and $L_{i+2} = L_i + L_{i+1}$ for all $i \geq 1$, which serve as a convenient abbreviation $L_n = F_{n-1} + F_{n+1}$ (Koshy [15] 97 #32). The formulas for t at even n are simplified using the relations $L_{4k} = 5(F_{2k})^2 + 2$ and $L_{4k+2} = (L_{2k+1})^2 + 2$ (Koshy [15] 97 #42 and 41), while $r\mu + s(\mu + \nu) = 1$ uses $F_{n+2} - F_{n-2} = L_n$ (Koshy [15] p. 97 #33).

**Proof**: For d, one in each case uses the lemma inductively to shift the indices close together, then simplifies by Fibonacci and gcd properties. Here are the details for the case n = 4k+2 giving the largest values: $d = \gcd(F_{n-2} + 1, F_{n-1} - 1) = \gcd(F_{4k} + F_2, F_{4k+1} - F_1) = \ldots = \gcd(F_{2k+2} + F_{2k}, F_{2k+3} - F_{2k-1}) = \gcd(F_{2k+2} + F_{2k}, F_{2k+1} + F_{2k+2} - F_{2k-1}) = \gcd(F_{2k+2} + F_{2k}, F_{2k-1} + F_{2k} + F_{2k+2} - F_{2k-1}) = F_{2k} + F_{2k+2}$.

The values of t follow from $t = [F_{n-1} + F_{n+1} - 1 - (-1)^n]/d$, and one may verify that the given values of r and s yield $r\mu + s(\mu + \nu) = 1$. The proofs for r, s, and h rely on the identity $F_iF_{j-1} - F_{i-1}F_j = (-1)^{i-1}F_{j-i}$ (Koshy [15] 87 #2, with k = 1). ■

To summarize numerically, here are the formulas and values of d, t, and h for small n; also values for a and b that give full wheels, with n distinct terms, in the cyclic group $Z_t$ (or in $Z_2 \oplus Z_2$, $Z_4 \oplus Z_4$, $Z_3 \oplus Z_3$, and $Z_7 \oplus Z_7$ for n = 3, 6, 8, and 16, respectively):

| n | 3 | 4 | 5 | 6 | 7 | 8 | 9 | 10 | 11 | 12 | 13 | 14 | 15 | 16 | 17 | 18 | 19 | 20 | 21 | 22 | 23 ... |
|---|---|---|---|---|---|---|---|---|---|---|---|---|---|---|---|---|---|---|---|---|---|
| d | 2 | $F_2$ | 1 | $L_3$ | 1 | $F_4$ | 2 | $L_5$ | 1 | $F_6$ | 1 | $L_7$ | 2 | $F_8$ | 1 | $L_9$ | 1 | $F_{10}$ | 2 | $L_{11}$ | 1 |
| | 2 | 1 | 1 | 4 | 1 | 3 | 2 | 11 | 1 | 8 | 1 | 29 | 2 | 21 | 1 | 76 | 1 | 55 | 2 | 199 | 1 ... |
| t | $\frac{L_3}{2}$ | $5F_2$ | $L_5$ | $L_3$ | $L_7$ | $5F_4$ | $\frac{L_9}{2}$ | $L_5$ | $L_{11}$ | $5F_6$ | $L_{13}$ | $L_7$ | $\frac{L_{15}}{2}$ | $5F_8$ | $L_{17}$ | $L_9$ | $L_{19}$ | $5F_{10}$ | $\frac{L_{21}}{2}$ | $L_{11}$ | $L_{23}$ |
| | 2 | 5 | 11 | 4 | 29 | 15 | 38 | 11 | 199 | 40 | 521 | 29 | 682 | 105 | 3571 | 76 | 9349 | 275 | 12238 | 199 | 64079 ... |
| h | $\frac{1-L_1}{2}$ | $-2F_2$ | $\frac{1-L_3}{2}$ | $-L_3$ | $\frac{1-L_5}{2}$ | $-2F_4$ | $1-L_7$ | $-L_5$ | $\frac{1-L_9}{2}$ | $-2F_6$ | $\frac{1-L_{11}}{2}$ | $-L_7$ | $1-L_{13}$ | $-2F_8$ | $\frac{1-L_{15}}{2}$ | $-L_9$ | $\frac{1-L_{17}}{2}$ | $-2F_{10}$ | $1-L_{19}$ | $-L_{11}$ | $\frac{1-L_{21}}{2}$ |
| | 0 | -2 | $\frac{-3}{2}$ | -4 | -5 | -6 | -28 | -11 | $\frac{-75}{2}$ | -16 | -99 | -29 | -520 | -42 | $\frac{-1363}{2}$ | -76 | -1785 | -110 | -9348 | -199 | $\frac{-24475}{2}$ ... |
| a | (1,0) | -1 | -2 | (1,0) | -1 | (1,0) | -2 | -1 | -2 | -1 | -1 | -1 | -2 | (1,0) | -2 | -1 | -1 | -1 | -2 | -1 | -2 ... |
| b | (0,1) | 2 | 3 | (0,1) | 5 | (0,1) | 28 | 1 | 75 | 2 | 99 | 1 | 520 | (0,1) | 1363 | 1 | 1785 | 2 | 9348 | 1 | 24475 ... |

**Table 2.3.** Formulas and values of d, t, h, a, and b for small n.

## 3. Difference cographs

A *difference cograph* is a cograph in which the points are elements in either Z, with $\mathcal{C}(P,Q) = |P-Q|$, or in $Z_n$, with $\mathcal{C}(P,Q) = \min |p-q|$, the minimum being taken over elements p in coset P and q in coset Q. Figure 3.1 shows two examples:

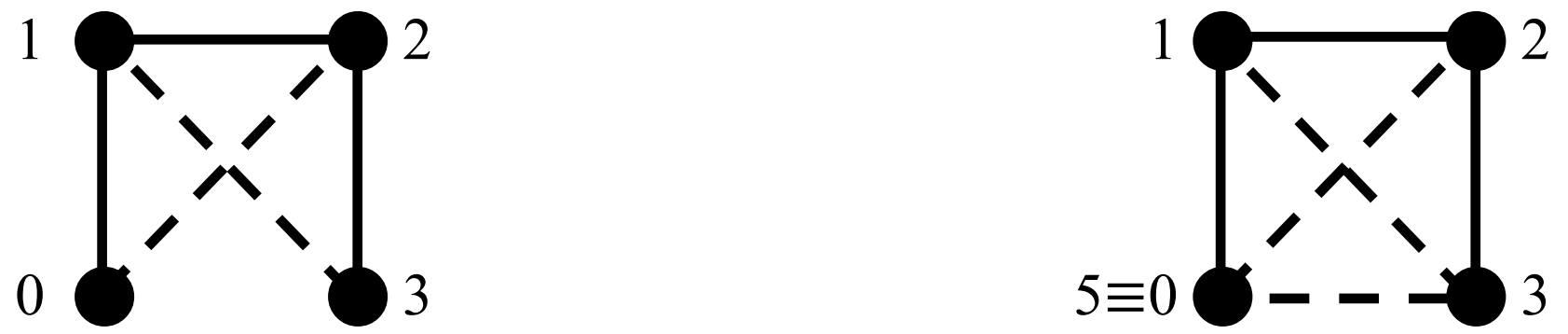

Over Z
(single-copy edge $\mathcal{C}(0,3)$ not shown)

Over $Z_5$

**Fig. 3.1.** Difference cographs on four points.

Difference cographs have some similarities to sum cographs, the points being again basically just numbers, and the edges just differences rather than sums. But the absolute value $|P-Q|$ in difference cograph edges makes the latter class more complicated, and hence also more numerous, because the simple configuration of Figure 3.2a forbidden in sum cographs can easily occur, as both examples illustrate. (Configuration 3.2b continues to be forbidden,though.)

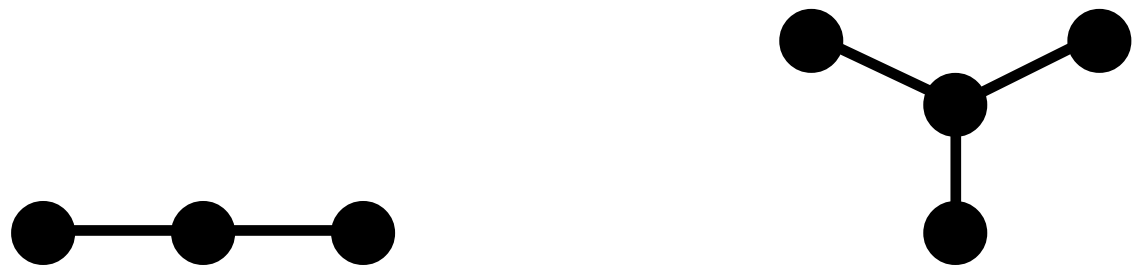

a) Allowed

b) Forbidden

**Fig. 3.2.** Configurations allowed or forbidden in difference cographs.

Difference cographs are therefore better viewed not as algebraic, like sum cographs, but rather as the simplest type of geometric (or "distance") cographs: one-dimensional arrays of points and distances measured along either a line or a circle.

Since any (finite) line segment can be wrapped into a semicircle, any (finite) difference cograph can be considered part of a sufficiently large $Z_n$. But the cographs embeddable in Z are a specialized subclass, and will be distinguished from those (like the second example above) that *require* torsion.

This section considers elementary properties and examples of difference cographs, culminating in the complete listing of the 62 difference cographs on five points.

**Lemma 3.1 ("Labeling")**: *A difference cograph may be relabeled to give any selected single point any selected value in its group Z or $Z_n$.*

**Proof**: If the selected point is X, and the selected value Y, adding Y-X to each point relabels X as Y while leaving all edges (point-to-point distances) unchanged. ■

Chains of equal edges then yield characteristic patterns of point values:

**Lemma 3.2 ("Chain")**. *Consider chains of distinct points and equal edges in a difference cograph:*

1) 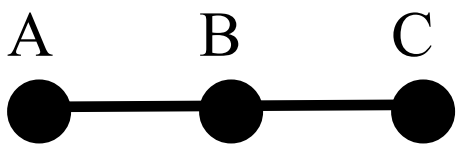

$\Rightarrow$ $A + C = 2B$ *("Midpoint Lemma")*

2) 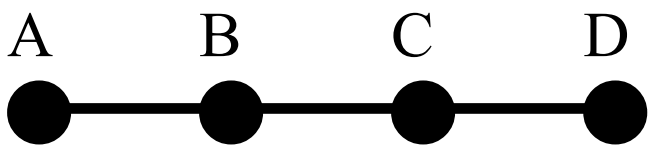

$\Rightarrow$ $A + D = B + C$

3) 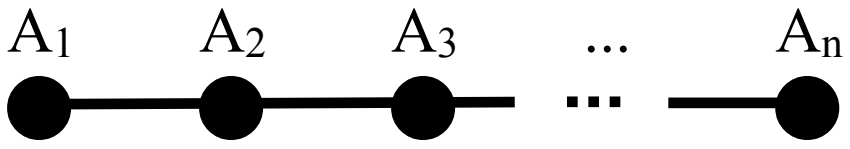

*(n odd)* $\Rightarrow$ $A_1 + A_n = 2A_{(n+1)/2}$

4) 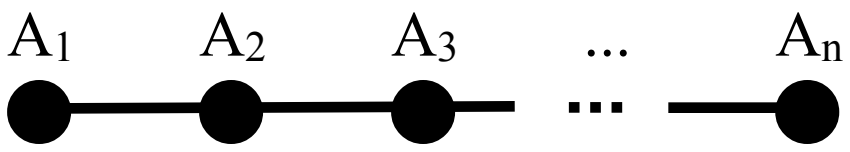

*(n even)* $\Rightarrow$ $A_1 + A_n = A_{n/2} + A_{n/2+1}$

**Proof**: Relabeling to make the first point 0, and calling the second X, then, since the points are distinct, the third point must be 2X, the fourth 3X, and so on. ■

Three characteristic motifs arise from repeated edges:

**Definition**: In a difference cograph, let a "V" be two copies of an edge incident at one point; a "Q" ("quadrilateral") be two disjoint copies of an edge; and a "T" ("triangle") be three copies forming a triangle (Figure 3.3).

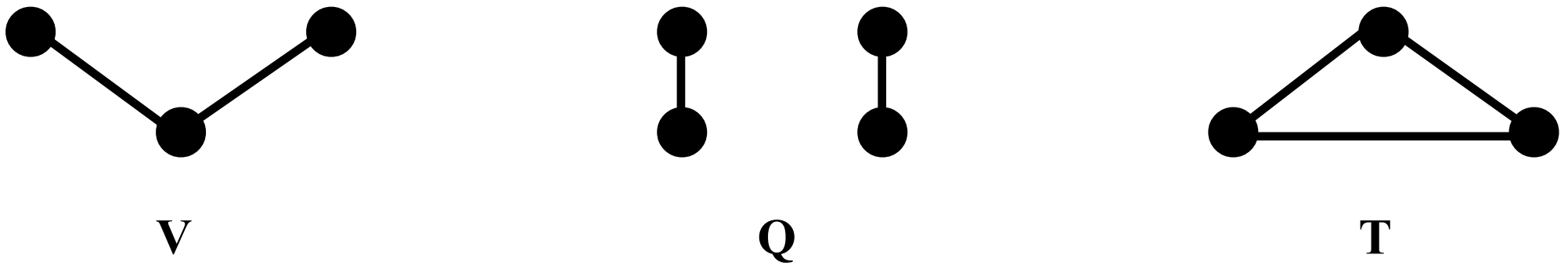

**Fig. 3.3.** Three characteristic motifs: V, Q, T.

**Lemma 3.3 ("Butterfly")**. *Two V's incident at a central point produce a Q (Figure 3.4).*

**Proof**: Label the central point 0. Then if the upper left point is X, the upper right one, forming a V, is -X; similarly, the lower two points are Y and -Y. But then $\mathcal{C}(X,Y) = |X-Y| = \mathcal{C}(-X,-Y)$, forming a Q. ■

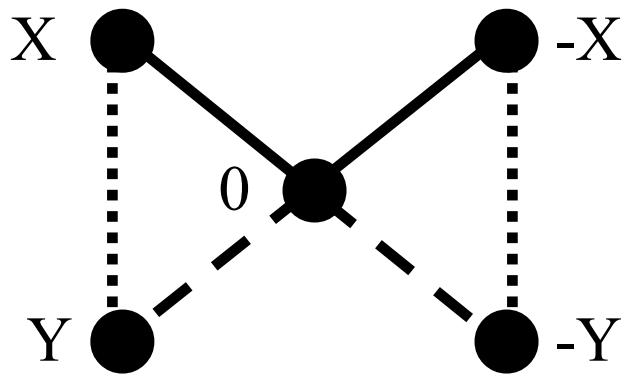

**Fig. 3.4.**
Two V's make a Q.

Q's are actually more complicated: there are five different types, *all* of which have at least a second pair of repeated edges:

**Lemma ("Q")**: *A Q in a difference cograph must be one of the five types in Figure 3.5:*

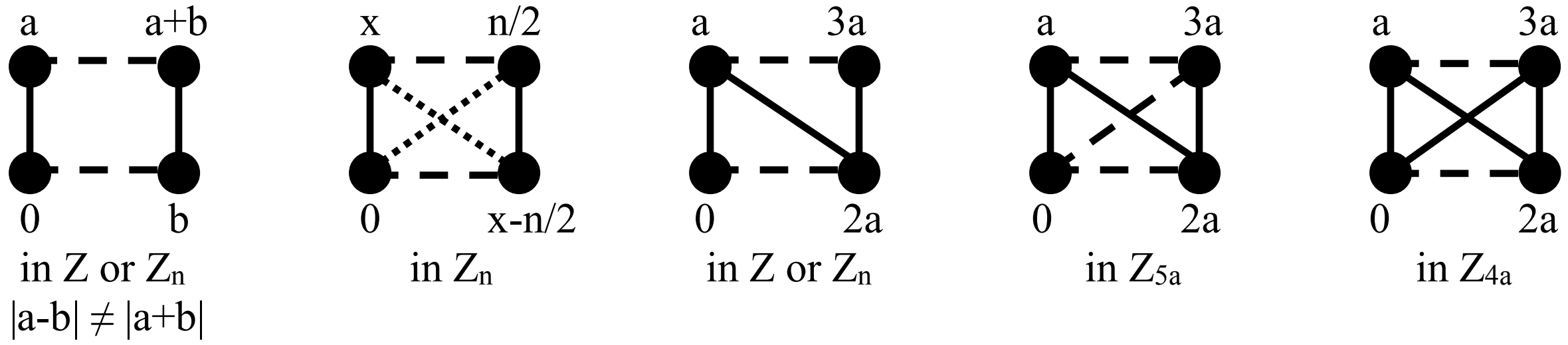

**Fig. 3.5.** The five possibilities for a Q.

**Proof**: Each of the five labelings evidently yields a Q. If, conversely, 0, a, b, x are points with $\mathcal{C}(0,a) = \mathcal{C}(b,x) =$ —— and $\mathcal{C}(0,b) = \mathcal{C}(a,x) =$ – – , then $x = b \pm a$. The alternative $x = b+a$ yields the first type, while $x = b-a$ gives $\mathcal{C}(0,x) = |a-b| = \mathcal{C}(a,b)$, again the first type (or, for either case, the second type if the diagonals are also equal). The third, fourth, and fifth cographs each contain a chain 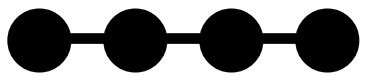 of points that, by the chain lemma, may be labeled 0, a, 2a, 3a, hence satisfy $\mathcal{C}(0,2a) = \mathcal{C}(a,3a) = 2a =$ – – $\neq$ —— $= a = \mathcal{C}(0,a)$, the choice among the three types then depending on the value of $\mathcal{C}(0,3a)$. ■

As a corollary, observe that if A,B,C,D,E are five points in a difference cograph in which A,B,C,D form a Q, then if $\mathcal{C}(A,E)$ is any of the *repeated* edges of the Q, the cograph must also contain a *second* Q (Figure 3.6). For the Q in A,B,C,D must be one of the five types from the Q lemma, and if it is, say, the first type, and $\mathcal{C}(A,E) =$ —— , then E,A,B,D will be a second Q. Hence in cataloguing the five-point difference cographs, sorted by their number of Vs, Qs, and Ts, requiring there to be just a single Q is a substantial restriction.

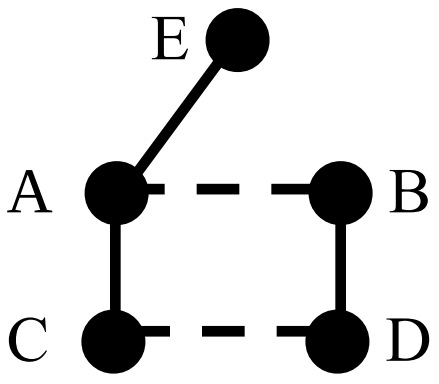

**Fig. 3.6.** Two Q's.

Catalogue Figures 3.7-3.8 below show the difference cographs (*repeated* edges only) on two, three, four, or five points, sorted by whether the cographs require torsion or not, and by the numbers of Qs, Vs, and Ts they contain. The legend gives details about the construction, point labels, cograph names, and uniqueness, as well as a sketch of the proof that the five-point table is complete.

<u>**Legend**</u>:

• Only multiple copy edges are shown.

• Point labels are low, but not rigorously checked to be least.

• Five-point cograph names:

The torsion-free cographs are sequentially named A-V.

Most of the torsion cographs can then be obtained by "reading" torsion-free structures in some $Z_n$ (e.g. $B_{20}$ is B read in $Z_{20}$).

Nine torsion cographs cannot be so obtained: these are named arbitrarily c,f,g,n,o,u,f',g', s', and the appropriate cyclic group $Z_n$ is indicated explicitly by $n \equiv 0$ in each diagram. Each of these nine *can* be obtained by reading an alternate *equivalent form* of some torsion-free cograph in some $Z_n$ (e.g. C = 0 1 2 5 8 could also be named--that is, is the same cograph as--$C_{alt}$ = 0 1 2 7 14, and $C_{alt}$ read in $Z_{16}$ is c).

• The cograph sorting by Qs, Vs, and Ts, together with simple criteria like point "degrees" (= the number of the *multiple copy* edges incident at each point), suffice to show that each cograph in the table is *unique*.

A helpful technique here is "Doubling analysis" of edges, by which the configuration ⋖ implies that edge ···· is twice edge —. In the case of the five five-point torsion cographs with Q = 0, V = 4, for example, one has:

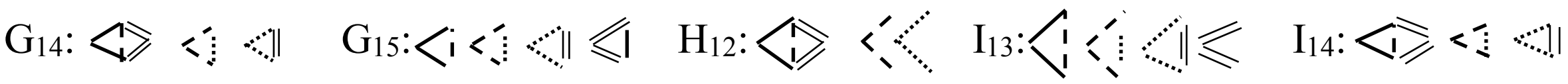

showing at one glance that at least four of the five are distinct.

• The Q-V-T sorting then also gives the framework for a proof that the five-point cograph table is *complete*.

**Sketch of proof**: At Q = 0, the butterfly lemma implies the Vs have distinct vertices, hence $V \leq 5$. The most complicated case here is V = 3: If the five points are a,b,c,d,e, and (without loss of generality) the three Vs have vertices a,b,c, then the V with vertex a--denote it $V_a$--has (without loss of generality) endpoints either b,c, or b,d, or d,e, whereas $V_b$ and $V_c$ each have six possible choices for their pair of endpoints, and two or three possibilities for edges. Two general constraints are:

1) If two Vs overlap (i.e., $V_x$ has endpoints y,z, and $V_y$ has endpoints x and w), they consequently have the same edge, and wyxz consequently would form a Q (contradicting Q = 0) unless w = z, forming a triangle; and

2) If $V_x$ has endpoints y and z, and $V_y$ has endpoints r and s, then the two Vs have different edges (else point y would have the forbidden configuration of three identical edges incident at one point).

One then checks that the individual cases (~84) all fall into the ten cases of the table.

At Q = 1, sorting by the five cases of the Q-lemma severely limits the possibilities, especially when V = 5 or 6, and making V > 6 impossible.

At Q = 2, one finds that the two Q's force there to be always at least one V, but that V > 5 is impossible.

When Q = 3, one finds that the cograph must contain either a triangle and a pair with the same edge, or else a 5-chain. These two constraints then are compatible only with the five cographs in the table.

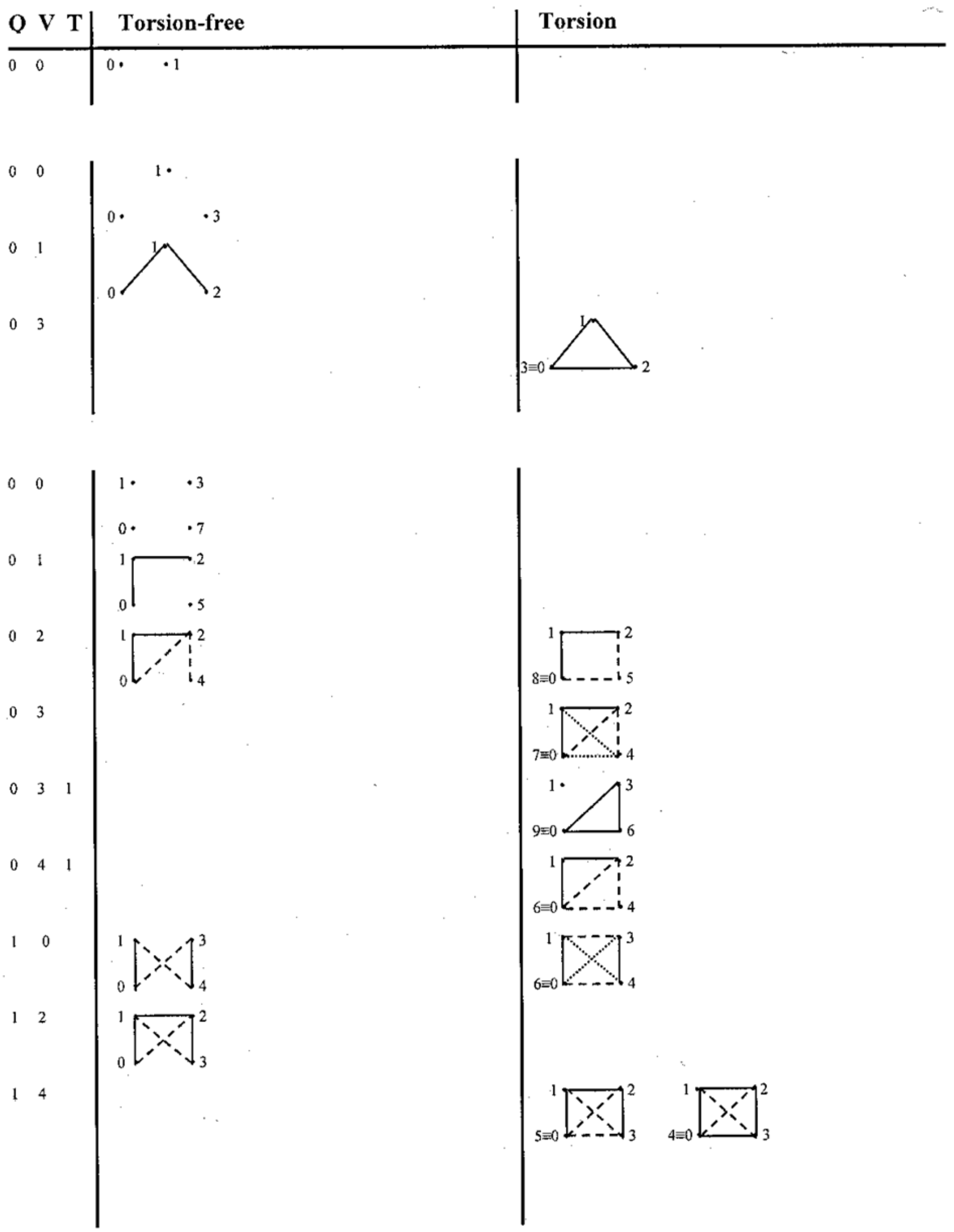

**Fig. 3.7.** Catalogue of two- three-, and four-point difference cographs.

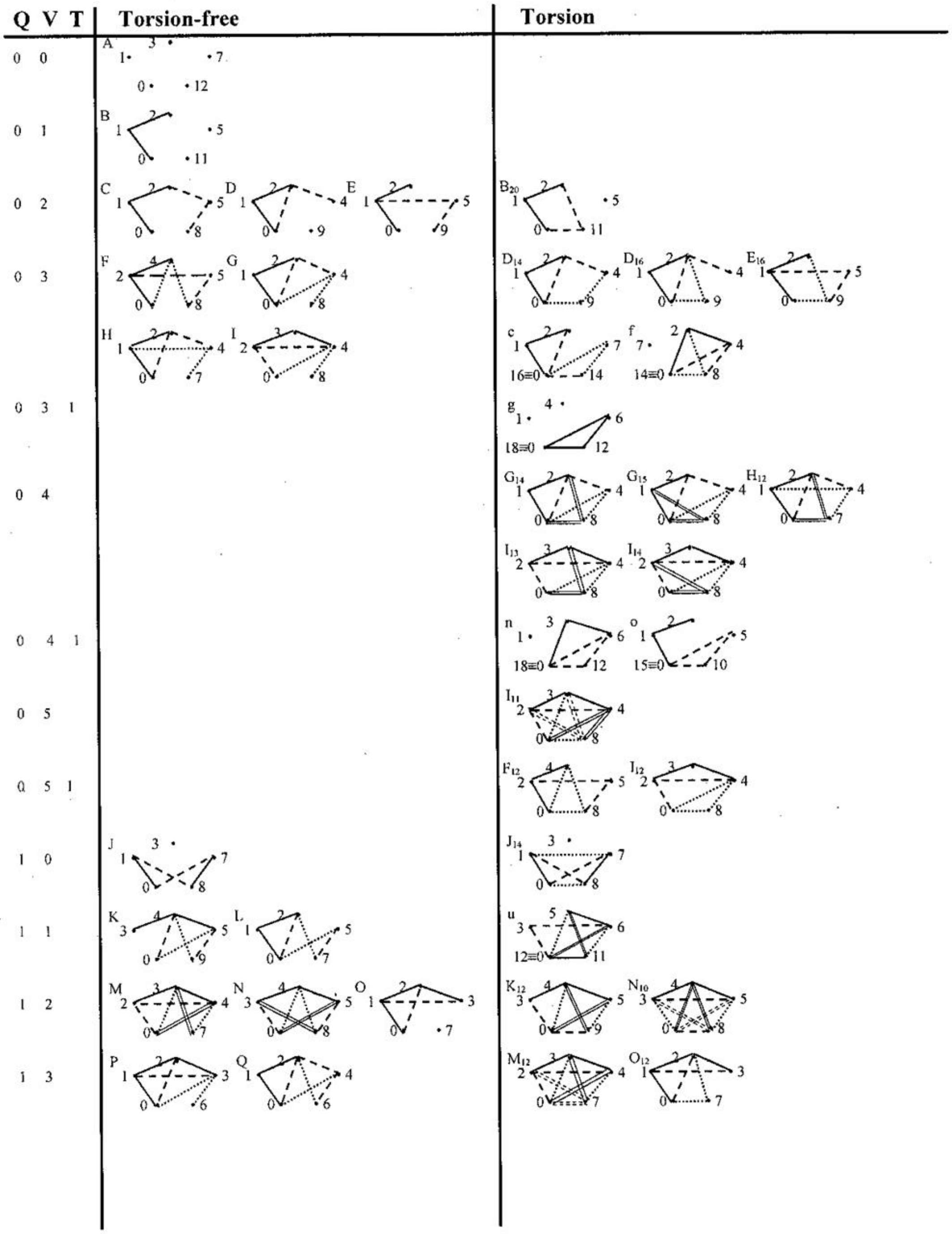

**Fig. 3.8.** Catalogue of five-point difference cographs.

(continued on next page)

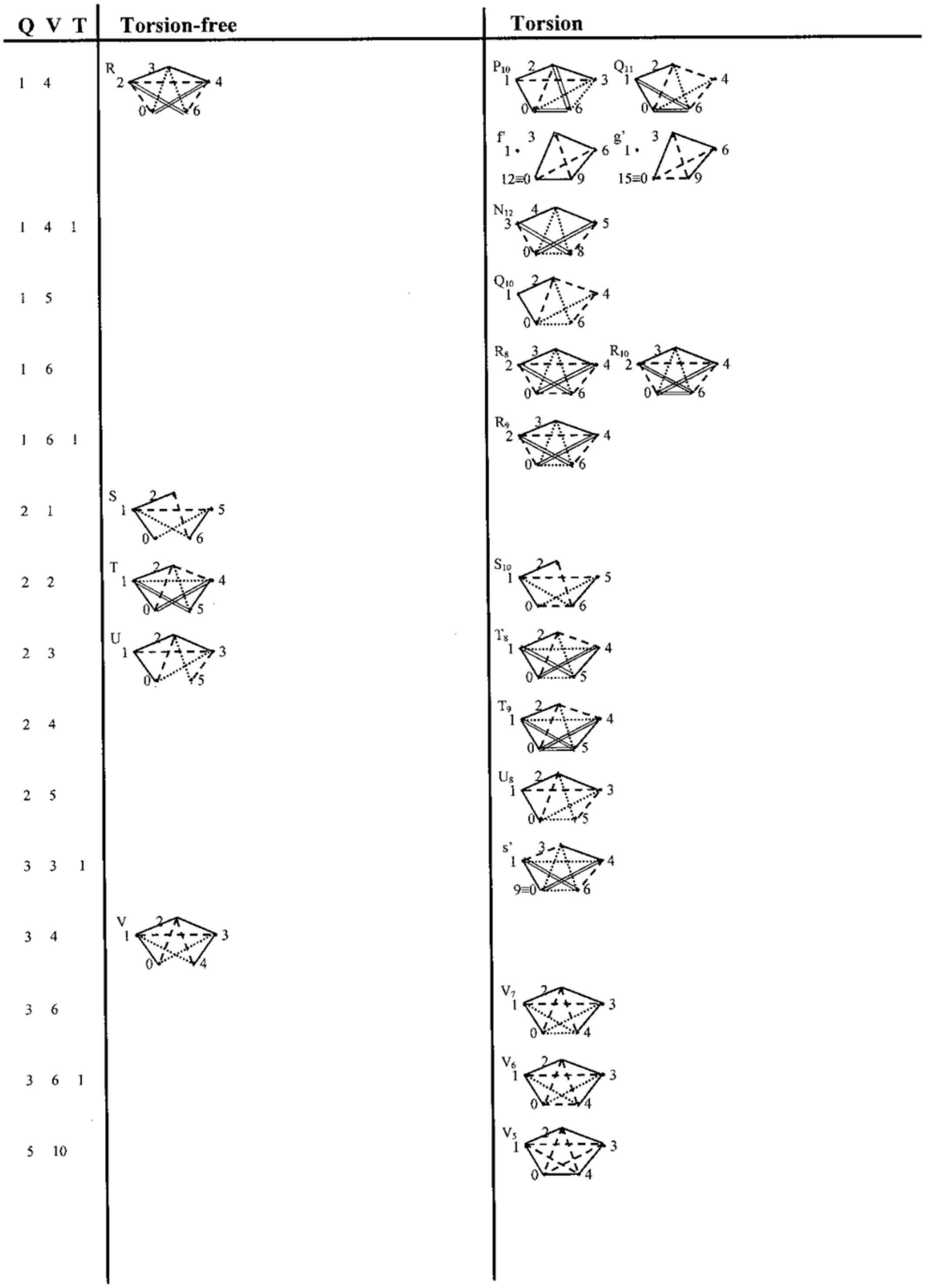

**Fig. 3.8.** Catalogue of five-point difference cographs.

(continued from previous page)

Looking at the torsion cographs in the tables, one can often see simple configurations that force the torsion. The following Proposition 3.1 itemizes some of them:

**Proposition 3.1 ("Torsion")**. *The following edge configurations compel torsion in a difference cograph:*

| | | Examples |
|---|---|---|
| a) **Cycle**: | 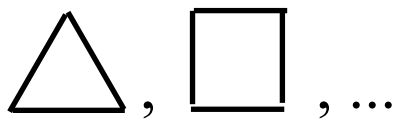 | $F_{12}$, f', $V_5$ |
| b) **Filled quadrangle**: | 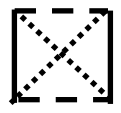 | $J_{14}$ |
| c) **Hemicycle pair**: | 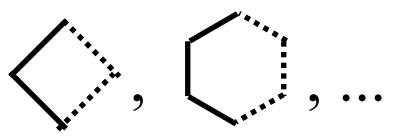 | $B_{20}$ |
| d) **Split hemicycle pair**: | 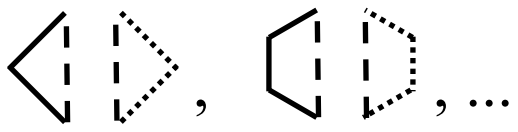 | c |
| e) **(Doubling, etc.) edge cycle**: | 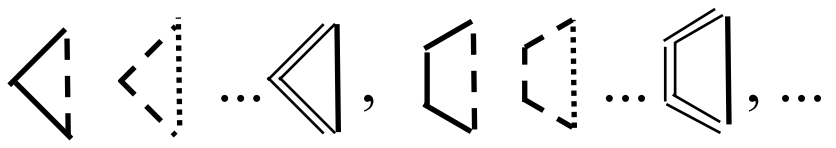 | $G_{14}$, $U_8$ |

**Proof**: a) As in the chain lemma, one may label the points of an n-chain by 0,a,2a, ..., (n-1)a, and closing the n-chain to an n-cycle then forces (n-1)a = 0. Parts (c), (d), and (e) follow similarly, while part (b) is the second type in the Q-lemma. ■

The Proposition accounts for all but four instances of torsion in the table: $I_{13}$, $K_{12}$, $M_{12}$, and $Q_{11}$. For $I_{13}$, edge-doubling analysis (see the table Legend) yields a triangle with sides a, 2a, and 16a, compelling either 13a = 0 or 19a = 0; while for $Q_{11}$, edge doubling and tripling yields 12a = ±a, so that 11a or 13a must be zero. Torsion in $K_{12}$ and $M_{12}$ comes ultimately from quadrangles they contain. It seems likely that there is an infinite expanding family of configurations that compel torsion for the cographs on six, seven, eight, etc., points.

**Remark: 2-Dimensional Four-Point Cographs**

As proved in Theorem 1.1, every finite cograph can be realized as a geometric cograph of sufficiently high dimension. That I have not explored geometric cographs beyond difference cographs, the one-dimensional simplest case, is undoubtedly due to my own limitations (as by nature an algebraist, not a geometer), not the subject's. As a glimpse ahead, it is straightforward to show that only one of the 25 possible four-point cographs, Figure 3.9, cannot be realized geometrically in two dimensions. The 24 realizations are mostly quite immediate, the most challenging being Figure 3.10.

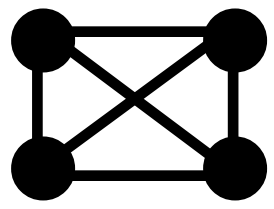

**Fig. 3.9.** Cannot be realized in two dimensions.

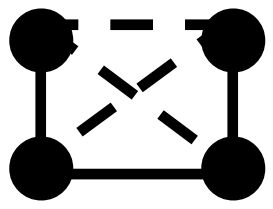

**Fig. 3.10.** Challenging to realize in two dimensions.

A realization of Figure 3.10--the *only* possible one, up to rigid motions, scaling, and switching the two sides x and y--is Figure 3.11:

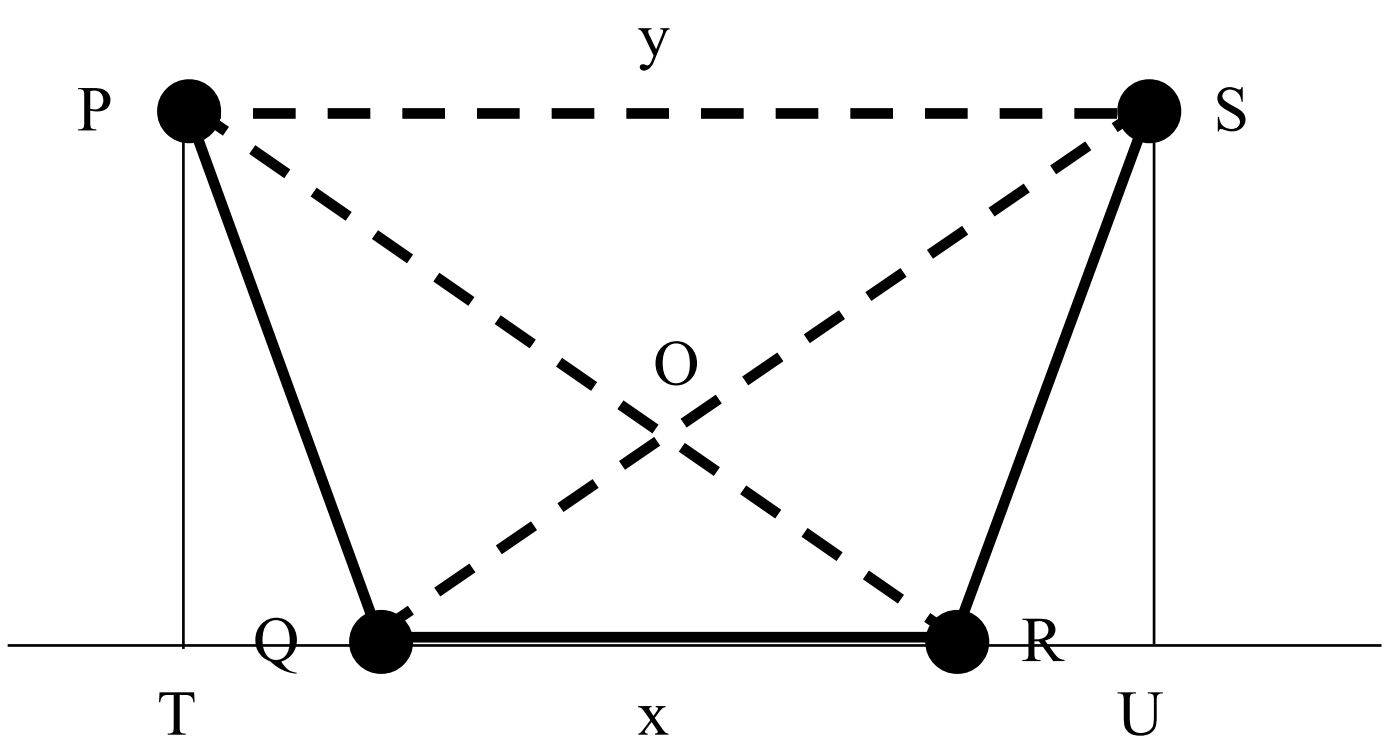

**Fig. 3.11. Geometric realization of the cograph of Fig. 3.10**.

where the angles QPR, RPS, PSQ, QSR, SQR, and PRQ are $\pi/5$; PQT, PQS, POQ, SOR, SRP, and SRU are $2\pi/5$; POS and QOR are $3\pi/5$; and TPQ and RSU are $\pi/10$; and consequently $y = (1+ 2 \sin \pi/10)\, x \approx 1.618\, x$. Drawing the bisector of PRS, to cut PS at a point V, in fact yields an isosceles triangle VRS similar to RPS, proving that PRS is a so-called "golden triangle," with ratio of sides y/x equal to the "golden ratio" $\varphi \approx 1.618$ that is a root $(1+ \sqrt{5})/2$ of the equation $\varphi^2 = \varphi + 1$.

## 4.1. Intersection cographs

An intersection cograph is a cograph in which the points and edges are sets, and each edge is the intersection of its two endpoints. Catalogue Figure 4.1.6 below gives the fifteen four-point intersection cographs (showing only the multiple-copy edges). That these cographs are all that are possible on four points will follow from a simple intersection cograph property:

**Proposition 4.1.1 (Quadrilateral rule)**. *Every quadrilateral abcd in an intersection cograph satisfies $a \cap c = b \cap d$.*

**Proof**: Call the points of the quadrilateral P, Q, R, and S (Fig. 4.1.1). Then $a \cap c = P \cap Q \cap R \cap S = b \cap d$. ■

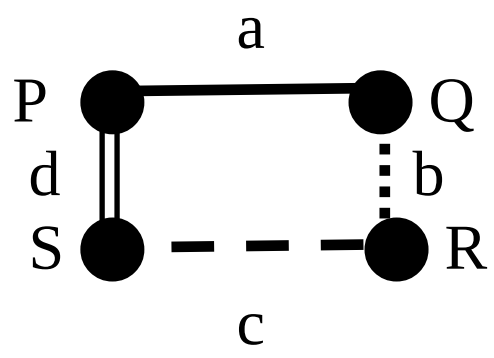

**Fig. 4.1.1.** Quadrilateral abcd.

**Corollary**: *The configurations of Figure 4.1.2 in an intersection cograph each imply $a \supset b$.*

**Proof**: Letting d = b in Proposition 4.1.1, $a \supset a \cap c = b \cap d = b \cap b = b$; the second part is similar. ■

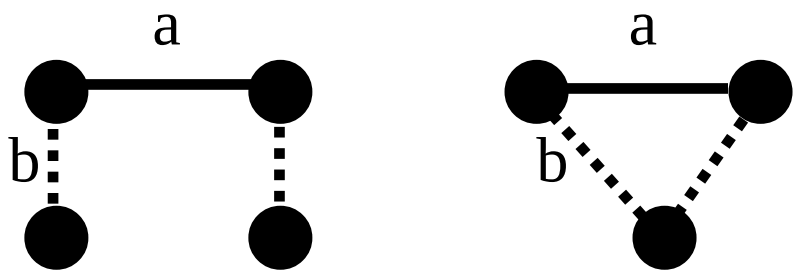

**Fig. 4.1.2.** Configurations that imply $a \supset b$.

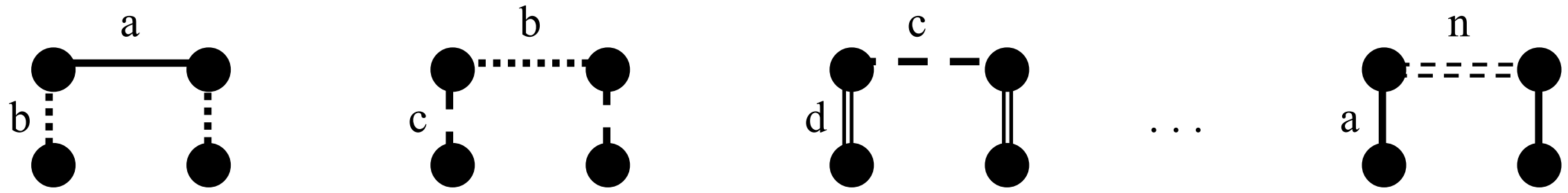

**Fig. 4.1.3.** A forbidden "inclusion cycle."

Any "inclusion cycle" of distinct edges as in Figure 4.1.3 is therefore forbidden in an intersection cograph, since it would entail $a \supset b \supset c \supset d \supset \ldots\ n \supset a$, forcing $a = b = c = d = \ldots = n$, which contradicts the assumption that the edges are distinct. The two simplest such forbidden configurations have cycle lengths two and three (Figure 4.1.4). Checking the Figure 1.2 listing of all 25 four-point cographs, one sees that every one omitted from the intersection cographs contains one of these forbidden configurations.

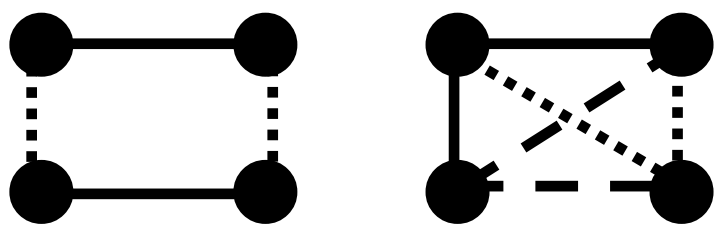

**Fig. 4.1.4** Two simple forbidden configurations.

Inclusion cycles are not the only type of obstruction, though:

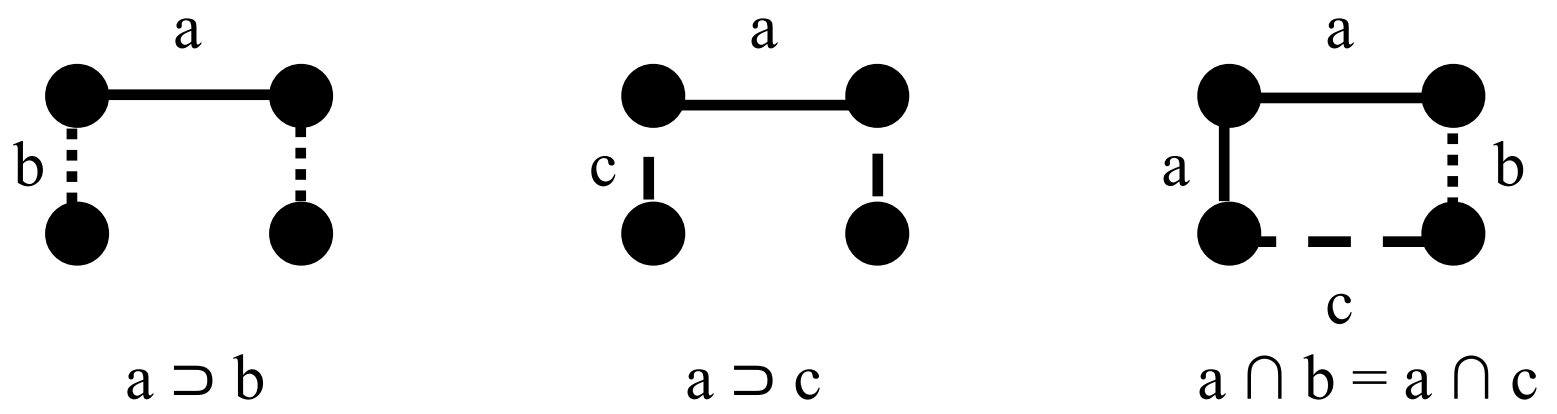

**Fig. 4.1.5.** Another type of obstruction.

The simple 12-point configuration of Figure 4.1.5 is forbidden in an intersection cograph too, since it yields the contradiction $b = a \cap b = a \cap c = c$. And the forbidden configuration can actually be packed into just six points, as in Figure 4.1.7.

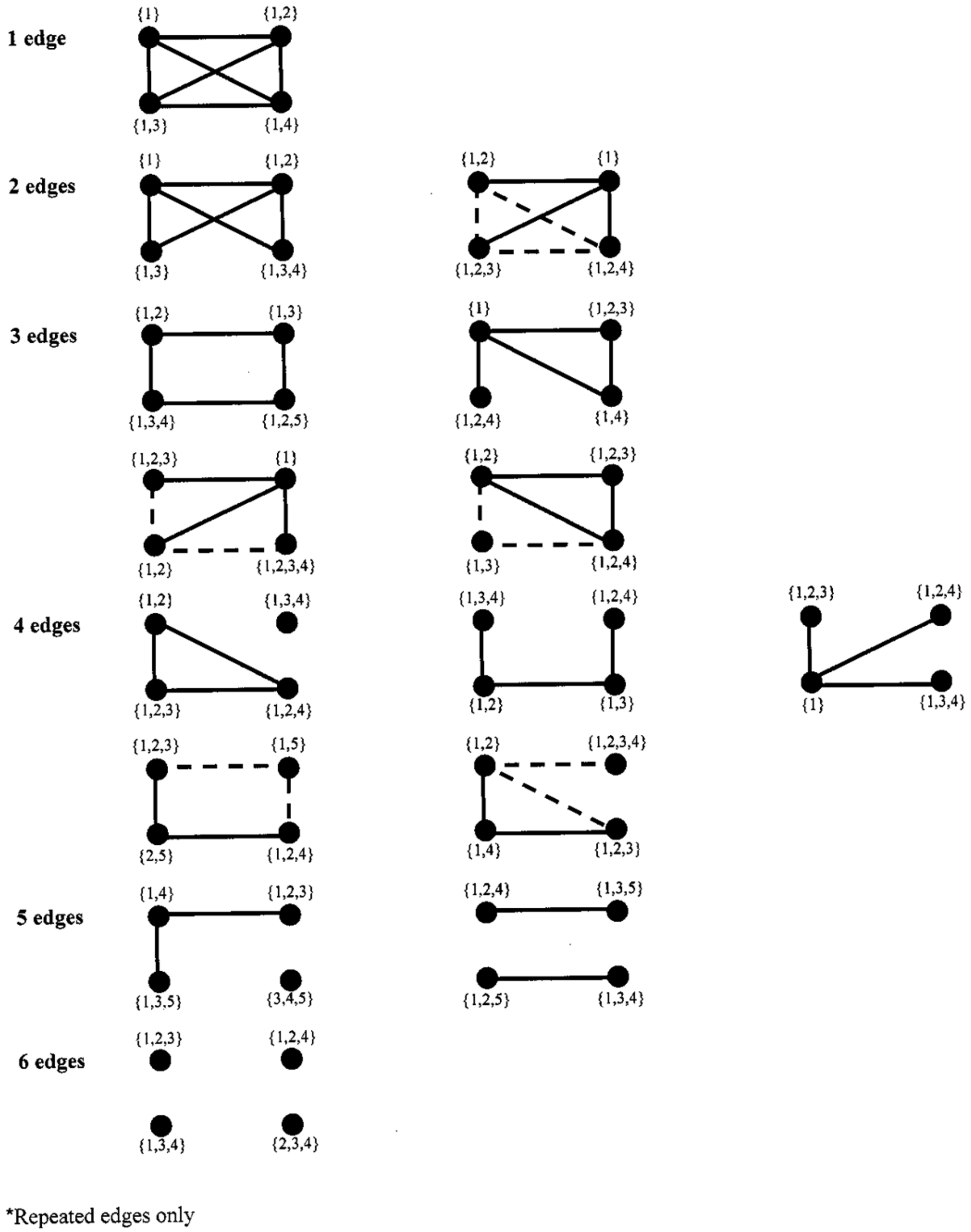

**Fig. 4.1.6.** Catalogue of four-point intersection cographs
(showing repeated edges only)

Proposition 4.1.1 has a parallel concerning triangles which is similarly proved by converting edge intersections to point intersections:

**Proposition 4.1.2 (Triangle rule).** *Any triangle abc in an intersection cograph satisfies* $a \cap b = b \cap c = a \cap c$. ■

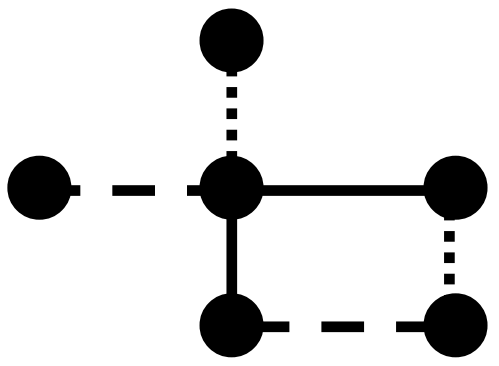

**Fig. 4.1.7.** The forbidden configuration.

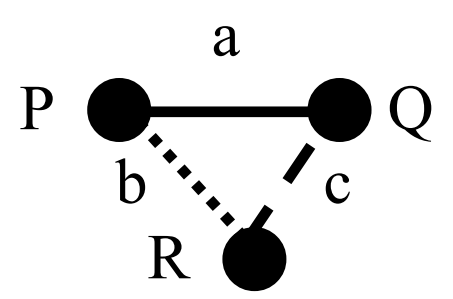

**Fig. 4.1.8.** Triangle rule.

**Remark**: It is easy to verify that the Proposition 4.1.2 conclusion has two further equivalents:

(1) $a \cap b = b \cap c = a \cap c = a \cap b \cap c$
(2) $a \cap b \subset c$, $b \cap c \subset a$, and $a \cap c \subset b$.

Propositions 4.1.1 and 4.1.2 operate in an intersection cograph, in which each edge is the intersection of its pre-existing endpoints. The next proposition, however, provides a UIE ("union of incident edges") construction which, beginning with a cograph in which just the edges are sets satisfying the rules of Propositions 4.1.1 and 4.1.2, generates points to yield a full compatible intersection cograph:

**Proposition 4.1.3 (UIE construction).** *Let* $\mathcal{C}$ *be a cograph in which the edges are sets, every triangle of which satisfies Proposition 4.1.2, and every quadrilateral satisfies Proposition 4.1.1. For each point* $P \in \mathcal{C}$ *define a set*

$$P' = \{P_o\} \cup \bigcup_{Q \in \mathcal{C}} \mathcal{C}(P,Q)$$

*Then the set of points* $\{P' : P \in \mathcal{C}\}$ *yields an intersection cograph having the edges of* $\mathcal{C}$. *(The singleton set* $\{P_o\}$ *contained in each P' is just a label, for the purpose of distinguishing points which may happen to have identical sets of incident edges.)* ■

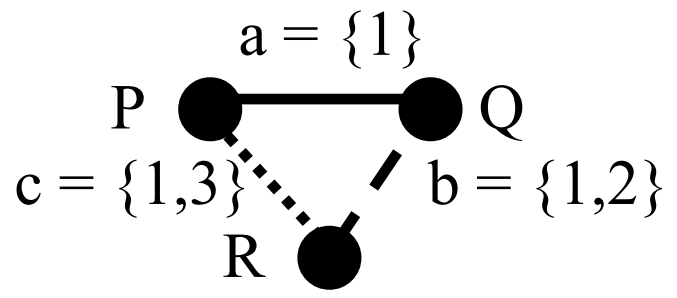

**Fig. 4.1.9.** Example of the UIE construction:
Let $P' = \{P_o\} \cup \{1,3\}$, $Q' = \{Q_o\} \cup \{1,2\}$, $R' = \{R_o\} \cup \{1,2,3\}$.

**Proof**: $P' \cap S' = [\bigcup_{Q \in \mathcal{C}} \mathcal{C}(P,Q)] \cap [\bigcup_{R \in \mathcal{C}} \mathcal{C}(S,R)] = \bigcup_{Q,\, R \in \mathcal{C}} [\mathcal{C}(P,Q) \cap \mathcal{C}(S,R)]$. The terms in this union fall into five cases:

1) $Q = S, R = P \Rightarrow \mathcal{C}(P,Q) \cap \mathcal{C}(S,R) = \mathcal{C}(P,S) \cap \mathcal{C}(S,P) = \mathcal{C}(P,S)$
2) $Q = S \Rightarrow \mathcal{C}(P,Q) \cap \mathcal{C}(S,R) = \mathcal{C}(P,S) \cap \mathcal{C}(S,R) \subset \mathcal{C}(P,S)$
3) $R = P \Rightarrow \mathcal{C}(P,Q) \cap \mathcal{C}(S,R) = \mathcal{C}(P,Q) \cap \mathcal{C}(S,P) \subset \mathcal{C}(P,S)$
4) $Q \neq S, R \neq P, Q = R \Rightarrow \mathcal{C}(P,Q) \cap \mathcal{C}(S,R) = \mathcal{C}(P,R) \cap \mathcal{C}(S,R) = \mathcal{C}(P,S) \cap \mathcal{C}(S,R) \subset \mathcal{C}(P,S)$
5) $Q \neq S, R \neq P, Q \neq R \Rightarrow \mathcal{C}(P,Q) \cap \mathcal{C}(S,R) = \mathcal{C}(P,S) \cap \mathcal{C}(Q,R) \subset \mathcal{C}(P,S)$

utilizing in case 4 that the edges of triangle PRS satisfy Proposition 4.1.2, and in case 5 that the edges of quadrilateral PQRS satisfy Proposition 4.1.1. The entire union is therefore $P' \cap S' = \mathcal{C}(P,S)$, as was to be proved. ■

**Remark 1**: The UIE-constructed cograph is not unique, but it is minimal: If {P”} also has intersection cograph $\mathcal{C}$, then for each P, $P'' \supset P' - \{P_o\}$, and the elements of P” - P’ are not contained in any other Q” - Q’, else they would appear in the intersection $P'' \cap Q''$.

**Remark 2**: An example illustrates the issue of representing an abstract cograph as an intersection cograph. Consider the abstract cograph PQRS of Fig. 4.1.10. Its edges form triangles aaa (from PQS), abb (from PQR), and abc (from SPR and SQR), and quadrilaterals abca (from PQRS), bbaa (from PRQS), and bcaa (from PRSQ). While triangle aaa yields no information, abb implies $a \cap b = b \cap b = b$, forcing $b \subset a$ (see the Corollary to Proposition 4.1.1), and abc then yields (see Proposition 4.1.2) that $a \cap c = a \cap b = b$. Since the three qadrilaterals turn out to contribute no further information, the general solution is to select a and c arbitrarily and let $b = a \cap c$. A satisfactory set representation is thus a = {1,2}, c = {1,3}, b = {1}; and UIE-constructed points P = {1,2}, Q = {1,2,4} (the element 4 added to distinguish it from P), R = {1,3}, and S = {1,2,3}.

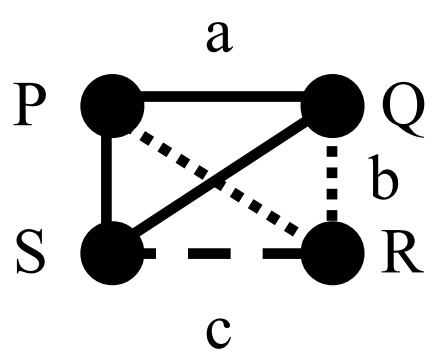

**Fig. 4.1.10.** Example of representation

## 4.2. Intersection cographs and aesthetics

Mathematics is such a beautiful subject it is not surprising that many mathematicians have been strongly moved by beauty, and a number of them have devised mathematical formulations of aesthetics. For instance, H. Weyl's book *Symmetry* traces the role of group theoretical symmetry in the visual arts [19], while Birkhoff offers a definition of beauty through a concept of "aesthetic measure" [4]. To a degree, of course, "beauty is in the eye of the beholder"; that is, the beauty must arise not solely from the beautiful object itself, but rather in the interaction of that object with the percipient, in the act of perception. This section will suggest how intersection cographs might offer a mathematical model for that interaction.

We begin by reformulating intersection cographs into a product binary form, less compact but more transparent for generalization, as illustrated by the following example:

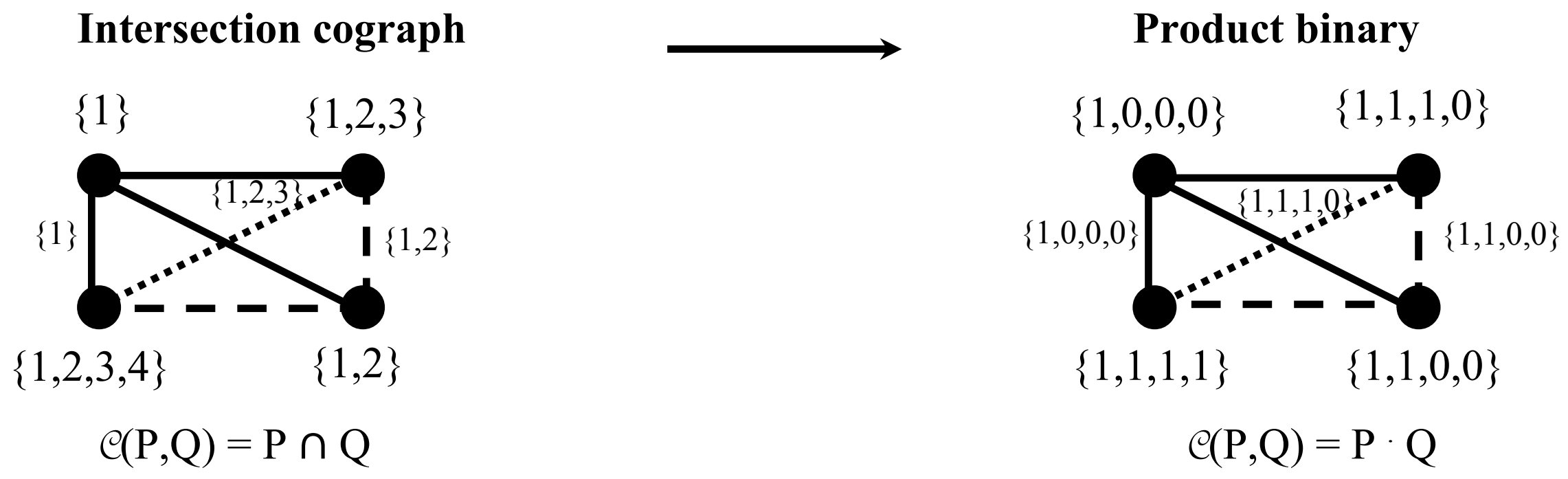

**Fig. 4.2.1.** Intersection cograph reformulated into product binary form.

Here the cograph on the left is the ordinary intersection cograph on the four points (sets) {1}, {1,2}, {1,2,3}, and {1,2,3,4}, with its three edges {1}, {1,2}, and {1,2,3}, represented respectively by solid, dashed, and dotted lines, given by the rule $\mathcal{C}(P,Q) = P \cap Q$. On the right the same abstract cograph is represented by four elements in the algebraic product $Z_2 \times Z_2 \times Z_2 \times Z_2$, where $Z_2$ is the "binary" two-element ring of integers mod two; each set of the intersection cograph is represented by its "characteristic function," and the cograph rule is $\mathcal{C}(P,Q) = P \cdot Q$. By this formulation it is evident how the notion might be generalized to products having many more, or even an infinite number, of "dimensions."

The approach here now is to view the perceived world abstractly as such an intersection cograph. Perceived "objects" are sets--the sets of perceptions (or "properties," or, in philosophy, "accidents") that an observer can ascertain from each of them. The focus of interest, for instance, in making the judgment "beauty," is to compare these sets among themselves, that is, to contemplate their intersections. Here, first, are three or four specific examples from the arts of painting, poetry, and music. They have been chosen for their extreme simplicity (and therefore rather unadorned, abstract character), to highlight the remarkable richness and complexity inherent in the judgment "beautiful".

The first example is the famous painting "Six Persimmons" by the thirteenth century Chinese painter Mu Qi (Figure 4.2.2) [16]. This picture is art of the utmost simplicity: six stylized pieces of fruit painted in black ink, without color, background, shadows, pictorial details, or

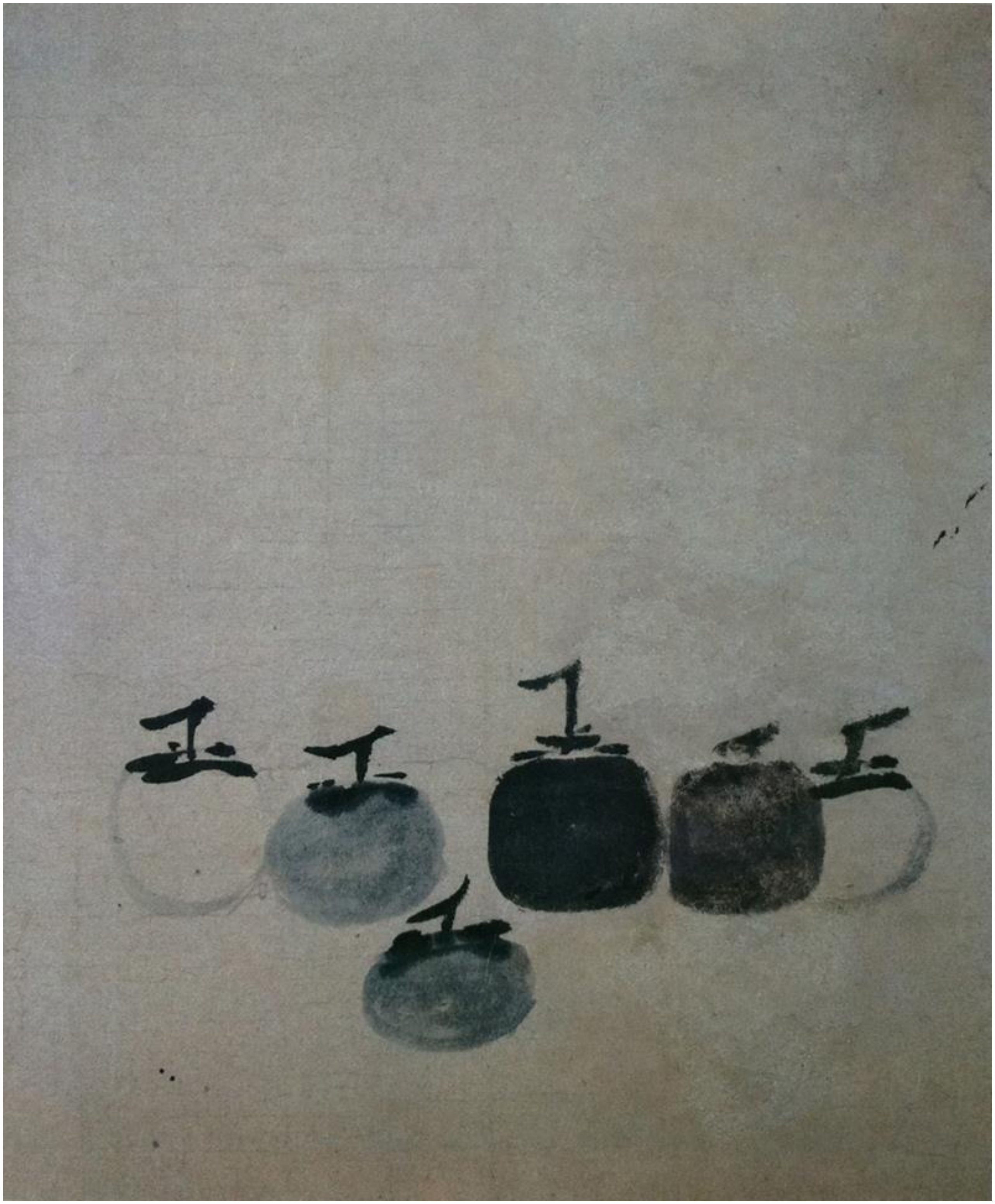

**Fig. 4.2.2.** Mu Qi: "Six Persimmons," 13th c. Chinese.

dramatic perspective. It can, nevertheless, arouse a powerful emotional response in a sensitive viewer: "passion ... congealed into a stupendous calm," is the reaction of one critic [18]. Intersection cographs can describe the process of perception as the appreciating eye plays over this picture:

| Aspect \ Persimmon | | 1 | 2 | 3 | 4 | 5 | 6 |
|---|---|---|---|---|---|---|---|
| a | Frontality | 0 | 0 | 1 | 0 | 0 | 0 |
| b | Frontality’ | 0 | 1 | 1 | 0 | 0 | 0 |
| c | Frontality” | 0 | 1 | 2 | 1 | 1 | 0 |
| d | Color | 0 | 1 | 1 | 3 | 2 | 0 |
| e | Size | 1 | 1 | 0 | 2 | 1 | 1 |
| f | Shape | 0 | 1 | 1 | 2 | 2 | 1 |
| g | Stem | 2 | 2 | 1 | 3 | 1 | 2 |

**Fig. 4.2.3.** Aspects of Mu Qi’s “Six Persimmons.”

The first line of the figure shows the result of the first glance (Aspect a): Of the six persimmons in the picture, numbers 1, 2, and 4, 5, 6 share the characteristic of being in one row in back, while number 3 is further to the front. Aspect b gives the second closer glance: fruit 2 is slightly ahead of the rest of the back row, and thus is united in similarity with number 3. Aspect c gives the most detailed look: the persimmons at both ends are subtly overlapped, hence behind the adjacent fruits. As one studies the painting, first one fruit and then another catches the eye, gaining prominence not only from position, but also by size, shape, or shading. Aspect d shows the cograph for color: the two fruits at the ends share the palest color, the fourth is darkest, the others are intermediate. Aspect e shows size; aspect f shape--oval, round, or squarish. Aspect g compares the lengths of the stems of the fruits. All of us become amateur artists when we make photos, and triumph when we center our friends in a snapshot. In his picture Mu-Qi has "balanced" all seven different aspects of position, shape, color, ... (and there are more) described in the cographs. The composition would be destroyed by omitting any one of the six fruits, and ruined by so much as shifting any position, size, shape, or color. It is this perfect equipoise of strong figural forces that produces the feeling of calm and passion noted by the critic; it is this exquisite balance that makes the painting a great work of art.

A similar balance characterizes the beauty of a poem. The example here is a verse from a lyric by Emily Dickinson (the first poem she valued highly enough to send to a critic [7]) describing the noble repose of the redeemed dead awaiting their resurrection on Judgment Day:

> Safe in their Alabaster Chambers--
> Untouched by Morning--
> And untouched by Noon--
> Sleep the meek members of the Resurrection,
> Rafter of Satin--and Roof of Stone--

The first impression one receives in reading this poem is perhaps the rhythm (Figure 4.2.4, Aspect a). The idiomatically mixed dactylic (¯ ˘ ˘) and trochaic (¯ ˘) pulse carries the words along to make them "verse" rather "prose", while placing special emphasis on emotionally important words like

Safe in their Alabaster Chambers--
Untouched by Morning--
And untouched by Noon--
Sleep the meek members of the Resurrection,
Rafter of Satin--and Roof of Stone--

**Aspect a: Rhythm**

Safe in their Alabaster Chambers--
Untouched by Morning--
And untouched by Noon--
Sleep the meek members of the Resurrection,
Rafter of Satin--and Roof of Stone--

**Aspect b: Vowel assonances**

Safe in their Alabaster Chambers--
Untouched by Morning--
And untouched by Noon--
Sleep the meek members of the Resurrection,
Rafter of Satin--and Roof of Stone--

**Aspect c: Consonant alliteration**

Safe in their Alabaster Chambers--
Untouched by Morning--
And untouched by Noon--
Sleep the meek members of the Resurrection,
Rafter of Satin--and Roof of Stone--

**Aspect d: Sense**

**Fig. 4.2.4.** Aspects of the Dickinson poem “Safe in their Alabaster Chambers.”

"safe," "untouched," and "sleep." The rhythm alone creates the powerful effect in the last line, where the drumbeat of the dactylic "Rafter of Satin" (¯ ˘ ˘ ¯ ˘) slows to the iambics "and Roof of Stone" (˘ ¯ ˘ ¯) to produce a feeling of unshakable solidity and firmness that will outlast the eons.

Reinforcing the poem’s initial rhythmic pattern then is the music of the language itself. Aspect b summarizes the vowel rhymes and assonances; these, for example, link "Safe" to "Chambers" in the first line, "Sleep" to "meek" as an internal rhyme in the fourth, and "Alabaster" in the first to "Rafter" and "Satin" in the last. The consonantal alliterations (Aspect c) provide even more numerous linkages. For example, the "s" sound in the first word "Safe" is echoed in "Alabaster," "Chambers," "Sleep," "members," "Resurrection," "Satin," and the last word "Stone." The "m," "r," and "t" sounds recur similarly. The alliterations also contribute notably to the effect of the last line, where the "r," "f," "s," "t," and "n" of "Rafter of Satin" are echoed exactly by those in "Roof of Stone."

Poetry, finally, requires a harmony of sense mutually reinforcing that of sound. Aspect d indicates some of the sense patterns in this lyric: The central thought is how the physical environment ("Chambers," "Rafter," and "Roof"), charged with emotional connotations of protection and permanence ("Alabaster", "Satin," "Stone"), shields its inhabitants from time ("Morning," "Noon," and "Resurrection"). The great majority of words in the lyric express this protection: "Safe," "Alabaster," "Chambers," "Untouched," "Untouched," "Sleep," "Meek," "Rafter," and "Roof." As with the "Six Persimmons" painting, this lyric is created from only a few ingredients. The exquisite rightness and economy of its crafting, each word linked to the others in a balance of rhythm, sound, and sense, make it, too, a great work of art.

Though this example deals with the minutest elements of sound and meaning, intersection cographs also easily represent much larger literary structures. For example, the cograph of Figure 4.2.5, with solid, dashed, or invisible white line segments, shows the main characters and relationships in Shakespeare’s play *King Lear*:

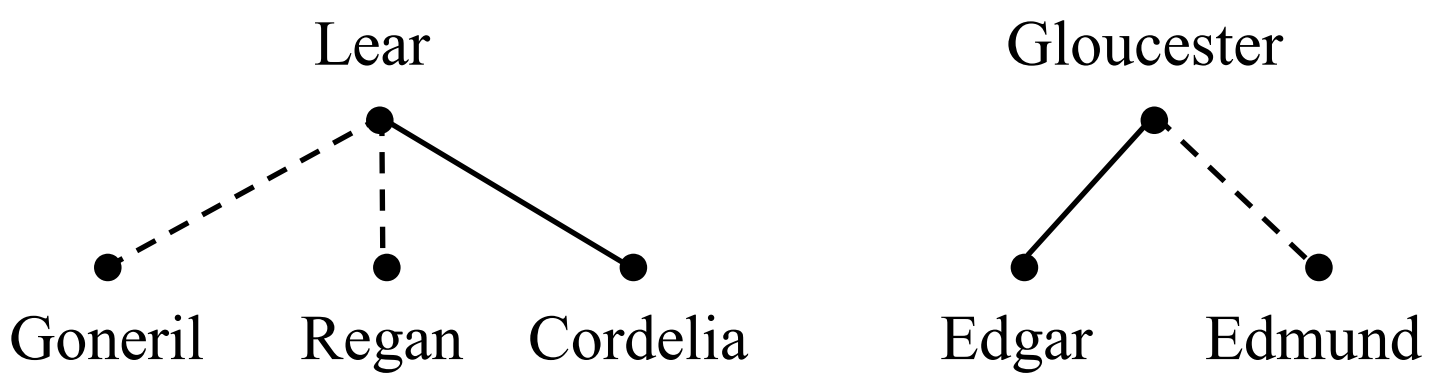

**Fig. 4.2.5.** Cograph of characters in Shakespeare's *King Lear.*

The central issue in this play, announced already in its third line "Is not this your son, my lord?", is the nature of the relationship between parent and child. The cograph schematizes, graphically and instantaneously, the two forms occurring here: the false one, between Lear and Goneril, Lear and Regan, and Gloucester and Edmund, and the true one, between Lear and Cordelia, and Gloucester and Edgar.

The final example is a musical one. It is difficult to find a profound piece of music on as miniature a scale as the Mu Qi painting or the Dickinson lyric, and we will content ourselves with a fragment, the first eight measures of the familiar beginner's minuet in G from the notebook of Anna Magdalena Bach (Fig. 4.2.6) [2]. Music, like poetry but unlike painting, is organized along a strictly linear pattern extended in time. Its first impression is therefore also the underlying rhythmical pattern. The rhythm is stricter for music than poetry, and the first rhythm cograph (not shown) simply records its steady 1-2-3 pattern of beats. This strict foundation, however, then permits the elaboration of more complex hierarchical structures: Fig. 4.2.7a highlights the repeated figure of four eighth notes leading up to a quarter note. Coincident with the rhythmic patterns are melodic and harmonic ones. Musical analysis (pioneered most formally by Schenker [8]) reveals these latter patterns most clearly by "rhythmic reduction" which omits ornamental filigree notes. The underlying pattern then stands out clearly: here, two simple scale passages, ascending, then descending to the tonic note G (Fig. 4.2.7b - circled notes). Other dimensions of musical expression include the shading of dynamics, ranging from soft to loud; progression of the underlying harmonies; small but important adjustments in tempo, such as ritards or accelerandos near musical climaxes; and, in ensemble music, use of the palette of colors of the different instruments. As with the other arts, an aesthetically satisfying musical composition or performance will be one in which the multiple dimensions of structure summarized schematically by the cographs are integrated into a convincing whole.

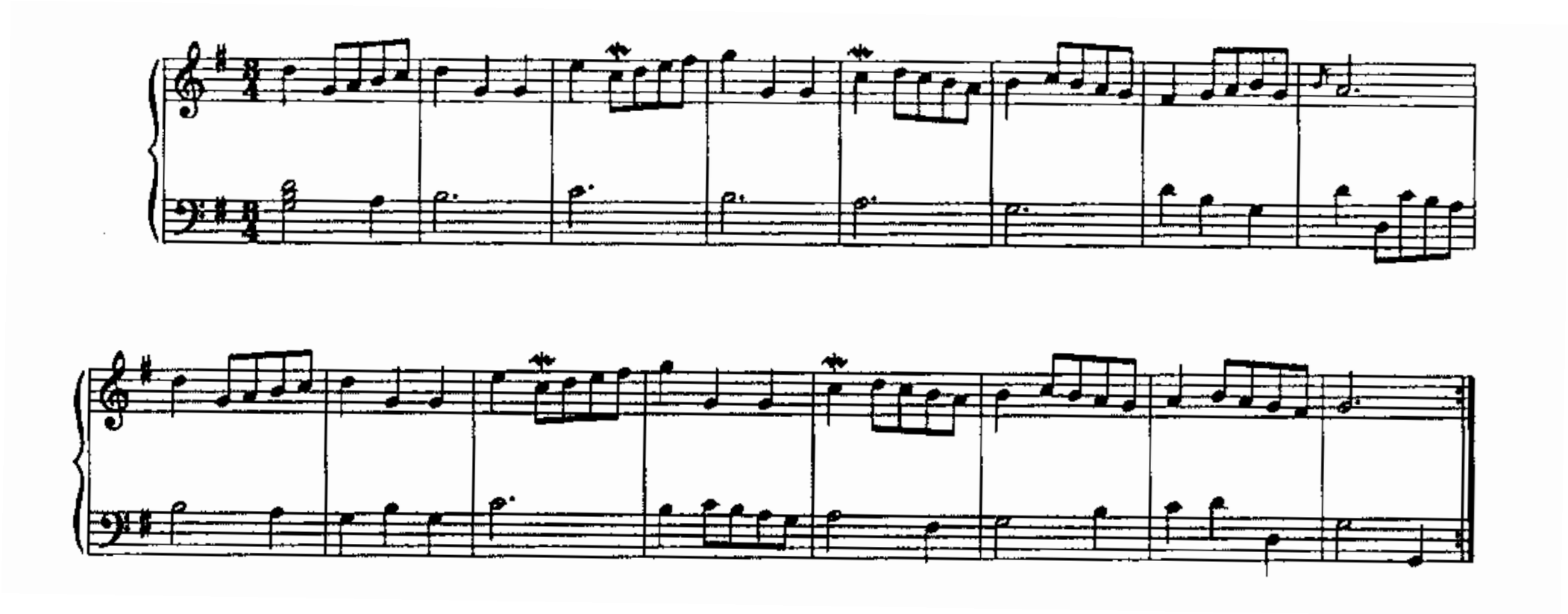

**Fig. 4.2.6.** Beginning of the Minuet in G, from the notebook of Anna Magdalena Bach.

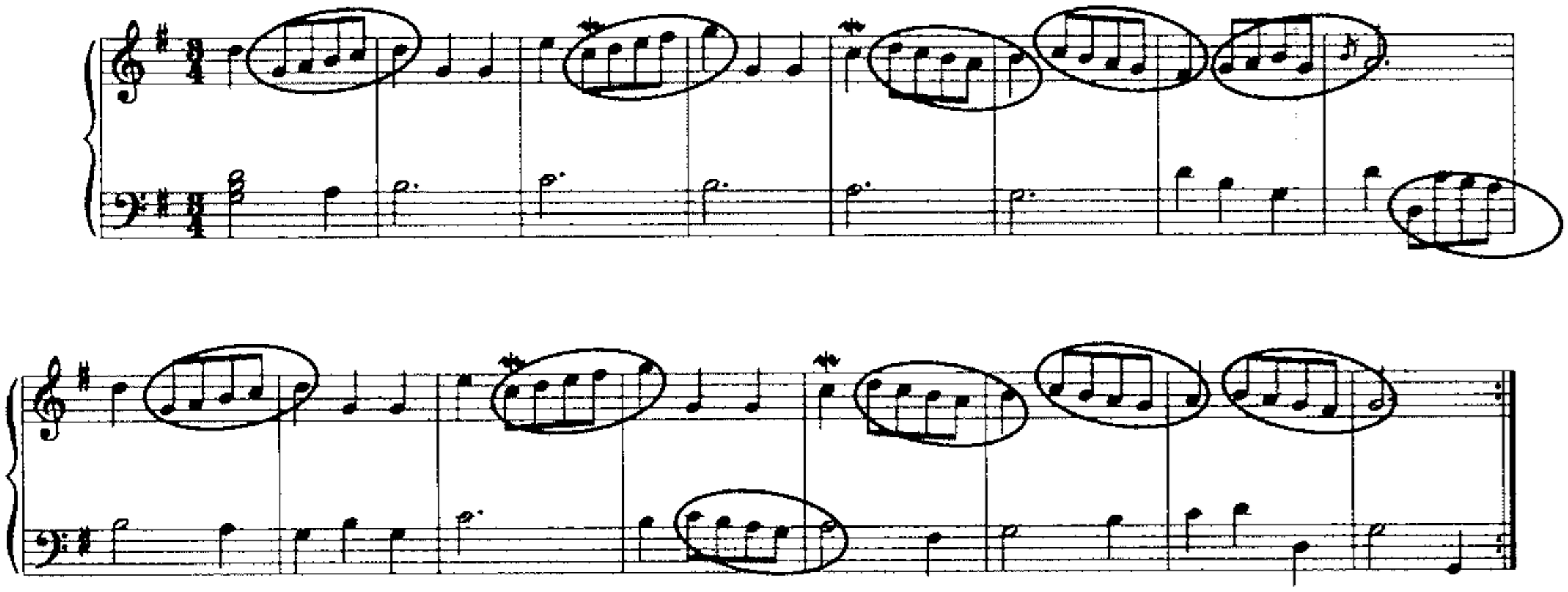

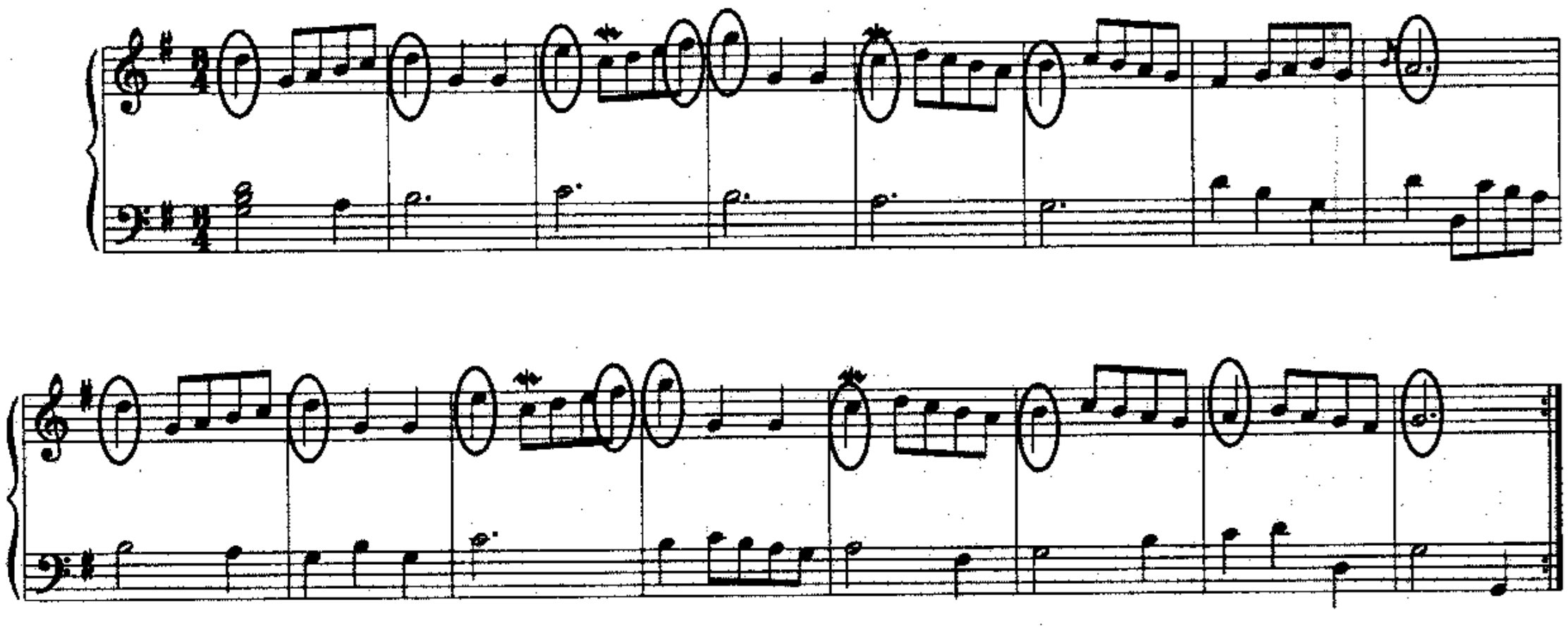

**Fig. 4.2.7.** Minuet in G: a) a repeated figure; b) the underlying scale pattern.

Intersection cographs might serve equally well in analyzing the “balancing of forces” required for stability in a variety of real-world contexts: for instance, the summation of attractive and repulsive electrostatic forces (the “Madelung constant”) in an ionic crystal; the psychological forces of personality, family background, and “chemistry” sustaining a compatible couple in their marriage; or sociological forces like nationality, race, gender, and class needed for a stable society.

## 5. P-L (Point-Line) - cographs

Two points in elementary geometry determine a line. The situation suggests a cograph structure. This section pursues the possibility by a class of PL (Point-Line)-cographs: It formulates equivalent definitions; gives examples; coordinatizes PL-cographs by a method originating in the specialized class of finite projective planes (Hall [11] 353-356); and gives a full tabulation for n = 2 through 7 points. The number of PL-cographs matches (Sloane [17]) the number of finite linear geometries (spaces) (Batten & Beutelspacher [3]), and the section concludes by proving that PL-cographs are equivalent to linear spaces.

How should a PL-cograph be defined? Points on a geometric line obey two special rules that suggest an appropriate definition.

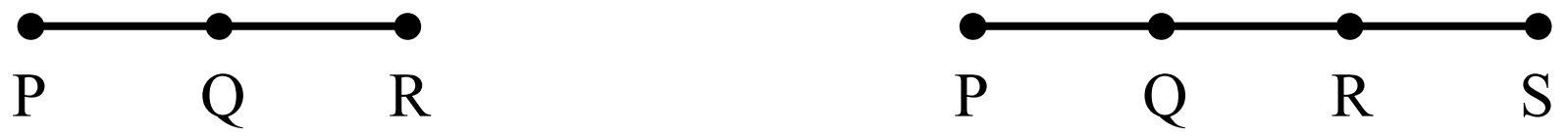

**Fig. 5.1.** Points on a geometric line obey two special rules.

First, if point pairs {P,Q} and {Q,R} lie on the same line, then so too does {P,R}. And secondly, if {P,Q} and {R,S} lie on the same line, then {P,R}, {P,S}, {Q,R}, and {Q,S} do too. Accordingly, one defines:

**Definition**: A *PL (point-line)-cograph* $\mathcal{C}$ is a cograph satisfying:

(1) $\mathcal{C}(P,Q) = \mathcal{C}(Q,R) \;\Rightarrow\; \mathcal{C}(P,Q) = \mathcal{C}(P,R)$

Symbolically,

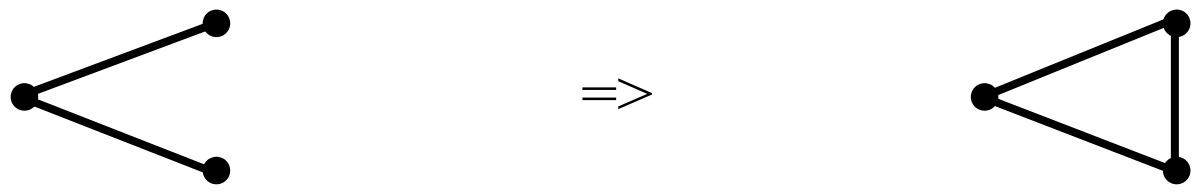

(2) $\mathcal{C}(P,Q) = \mathcal{C}(R,S) \;\Rightarrow\; \mathcal{C}(P,Q) = \mathcal{C}(P,R) = \mathcal{C}(P,S) = \mathcal{C}(Q,R) = \mathcal{C}(Q,S)$

Symbolically,

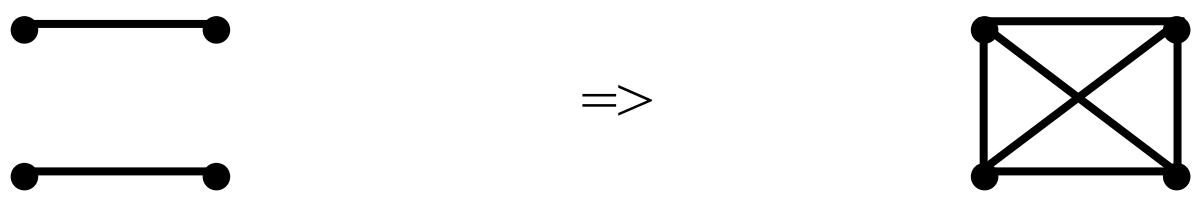

(A) 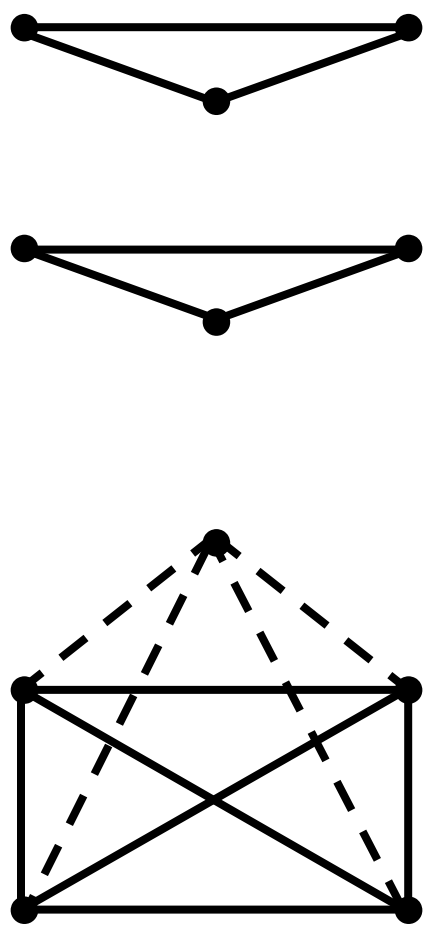

(All unmarked edges are single-copy)

(B)

**Fig. 5.2.** Two simple counterexamples (*not* PL-cographs) show that the defning conditions for PL-cographs are independent: (A) satisfies condition 1 but not 2; (B) satisfies condition 2 but not 1.

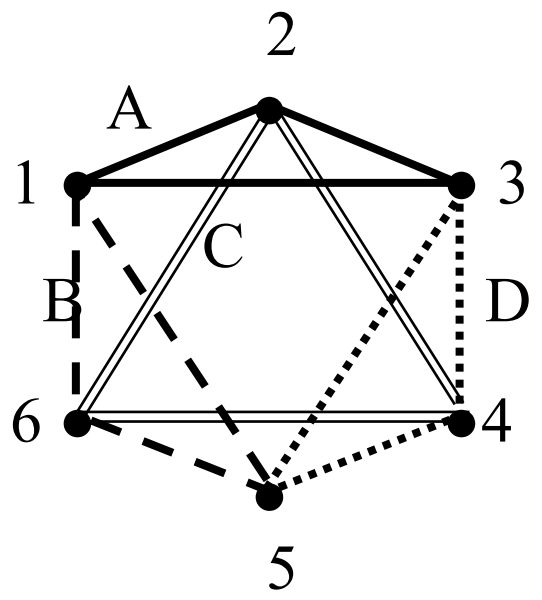

(The three unmarked edges 𝒞(1,4) = E, 𝒞(2,5) = F, and 𝒞(3,6) = G are single-copy.)

**Fig. 5.3.** Example: A six-point seven-edge PL-cograph.
Visualizing the figure as an octahedron makes it evident that the four edges A, B, C, and D are equivalent.

Figure 5.5 further below shows all the PL-cographs on n ≤ 7 points.

In studying PL-cographs it is useful to consider the *graphs*--term them *blocks*--formed by single edges and their endpoints. The definition of a PL-cograph is that every 3- or 4-point block is a complete graph. Here are some equivalent characterizations by these terms:

**Proposition 5.1.** *The following four characterizations of a PL-cograph are equivalent*:

*1)* 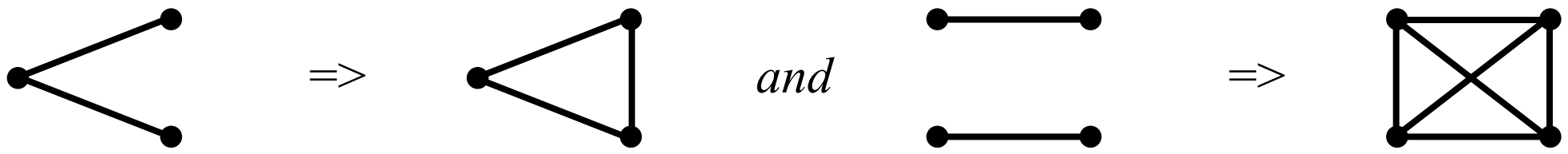

*2) Every 3- or 4-point block is a complete graph.*

*3) Every block is a complete graph.*

*4) Any two distinct blocks intersect in at most one point.*

**Proof**: (1), the definition of PL-cographs, is clearly equivalent to (2), which is a special case of (3). (2) => (3) because if P and Q are in the block with edge e, then ∃ R and S with $\mathcal{C}$(P,R) = $\mathcal{C}$(Q,S) = e. $\mathcal{C}$(P,Q) = e then follows by the first rule in (1) when R = S, and by the second when R ≠ S.

(1) => (4) because if the intersection of blocks with differing edges e and f contained two points P and Q, the rules of (1) would force $\mathcal{C}$(P,Q) to be both e and f, a contradiction. If, conversely, $\mathcal{C}$ contained then the blocks for e and f would intersect in more than one point, contradicting (4); the second rule of (1) follows similarly from (4). ■

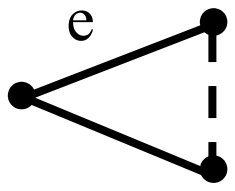

The six-point seven-edge example PL-cograph in Figure 5.3 above, for instance, has a three-point block for each edge A, B, C, or D, and two-point blocks {1,4}, {2,5}, {3,6} for the single-copy edges; and one may see that any two of these intersect in at most one point.

Projective planes are a specialized type of PL-cograph, in which not only do every two points determine a line, but also every two lines determine (i.e., intersect at) a point. These strong properties allow one to assign coordinates to each point P, as follows. Through P one finds the "vertical" line PY passing through a fixed point "at infinity" Y, and the "horizontal" line PX through a fixed point "at infinity" X. One next assigns arbitrary distinct coordinates $(x_i,x_i)$ to the points on the "diagonal" line OI passing through an "origin" O = (0,0) and a fixed point I = (1,1). Then when P is "finite" (i.e. PX ≠ PY), it receives coordinates $(x_i,y_i)$, where PY ∩ OI = $(x_i,x_i)$ and PX ∩ OI = $(y_i,y_i)$, and when P is on the line at infinity XY (i.e. PX = PY), it receives the single coordinate (m) such that (1,m) ∈ OP.

This process cannot be easily carried through for an arbitrary PL-cograph because the requisite line intersections need not exist. It will be shown here now, though, that the edges $\mathcal{C}$(P,X), $\mathcal{C}$(P,Y), and $\mathcal{C}$(P,O) themselves suffice to determine P uniquely. Here, first, is a lemma to establish that nontrivial PL-cographs have sufficiently many edges:

**Lemma 5.1.** *If PL-cograph $\mathcal{C}$ has more than one edge then it has at least three of them, and for any points P and Q there is an R so that $\mathcal{C}$(P,Q), $\mathcal{C}$(P,R), and $\mathcal{C}$(Q,R) are distinct.*

**Proof**: Let $\mathcal{C}$(P,Q) = e. If there is an R with $\mathcal{C}$(P,R) = f ≠ e, then $\mathcal{C}$(Q,R) cannot be e because PL-cograph rule (1) would then force $\mathcal{C}$(P,R) = e ≠ f contrary to assumption; and similarly, $\mathcal{C}$(Q,R) ≠ f.

But such an R always exists. For suppose that instead $\mathcal{C}(P,R_i)$ = e for all $R_i$ ≠ P in $\mathcal{C}$. By assumption $\mathcal{C}$ has another edge f ≠ e, say f = $\mathcal{C}$(R,S), with R and S different from P. But then by hypothesis $\mathcal{C}$(P,R) = e and $\mathcal{C}$(P,S) = e, which by PL-cograph rule (1) forces the contradiction $\mathcal{C}$(R,S) = e. ■

**Corollary.** *The pair of edges $\mathcal{C}(P,Q)$ and $\mathcal{C}(P,R)$ is incident uniquely to P.*

**Proof**: Suppose that $\mathcal{C}(P,Q) = e = \mathcal{C}(P',Q')$ and $\mathcal{C}(P,R) = f = \mathcal{C}(P',R')$ for P',Q',R' distinct from P, Q, and R. Then PL-cograph rule (2) dictates that $\mathcal{C}(P,P')$ must be both e and f, a contradiction. PL-cograph rules (1) and (2) yield similar contradictions in the special cases that, while keeping $P \neq P'$, identify one or more of P', Q', or R' with P, Q, or R. ■

As an example of Lemma 5.1, in the Figure 5.3 six-point seven-edge PL-cograph above $\mathcal{C}(1,2) = A$ is contained in triangle 126 with distinct edges A, B, and C, and only point 1 is incident to both A and B.

Lemma 5.1 now permits one from an *arbitrary* pair of edges X and Y to obtain an "origin" O and coordinatize the entire cograph:

**Theorem 5.1**. *Let X and Y be arbitrary points in a PL-cograph $\mathcal{C}$ having more than one edge. By Lemma 5.1 select a point O so that $\mathcal{C}(X,Y) = e$, $\mathcal{C}(O,X) = f$, and $\mathcal{C}(O,Y) = g$ are distinct. Then any point P different from X and Y is uniquely labeled by:*

| | |
|---|---|
| *the ordered pair ($\mathcal{C}(P,X)$, $\mathcal{C}(P,Y)$)* | *if $\mathcal{C}(P,X) \neq \mathcal{C}(P,Y)$* |
| *$\mathcal{C}(P,O)$* | *if $\mathcal{C}(P,X) = \mathcal{C}(P,Y)$* |

**Proof**: If $\mathcal{C}(P,X) \neq \mathcal{C}(P,Y)$ but there is a point $Q \neq P$ with $\mathcal{C}(P,X) = \mathcal{C}(Q,X)$ and $\mathcal{C}(P,Y) = \mathcal{C}(Q,Y)$, then by PL-cograph rule 1 $\mathcal{C}(P,Q)$ equals both $\mathcal{C}(P,X)$ and $\mathcal{C}(P,Y)$, a contradiction.

Alternatively, suppose $\mathcal{C}(P,X) = \mathcal{C}(P,Y)$ but there is a point $Q \neq P$ also satisfying $\mathcal{C}(Q,X) = \mathcal{C}(Q,Y)$, with $\mathcal{C}(P,O) = \mathcal{C}(Q,O)$. By PL-cograph rule 1, $\mathcal{C}(P,X) = \mathcal{C}(P,Y) \Rightarrow \mathcal{C}(P,X) = \mathcal{C}(P,Y) = \mathcal{C}(X,Y) = e$; and similarly, $\mathcal{C}(Q,X) = \mathcal{C}(Q,Y) = \mathcal{C}(X,Y) = e$. These then imply also that $\mathcal{C}(P,Q) = e$; and then $\mathcal{C}(P,O) = \mathcal{C}(Q,O) \Rightarrow \mathcal{C}(P,O) = \mathcal{C}(Q,O) = \mathcal{C}(P,Q) = e$ as well. But now $\mathcal{C}(P,X) = \mathcal{C}(P,O) = e \Rightarrow \mathcal{C}(O,X) = e \neq f$, a contradiction. ■

| | **Coordinatization** | | |
|---|---|---|---|
| **Point** | **#1** | **#2** | **#3** |
| 1 | ≡ X | ≡ X | ≡ X |
| 2 | ≡ Y | ≡ Y | ≡ O = (A,C) |
| 3 | (G) | (D) | (A,D) |
| 4 | (E,C) | ≡ O = (E,C) | ≡ Y |
| 5 | (B,F) | (B,F) | (B,D) |
| 6 | ≡ O = (B,C) | (B,C) | (B,C) |

**Fig. 5.4.** Example: three coordinatizations of the Fig. 5.3 six-point seven-edge PL-cograph.

In #1 $\mathcal{C}$(X,Y), $\mathcal{C}$(X,O), and $\mathcal{C}$(Y,O) are all multiple-copy; in #2 $\mathcal{C}$(X,O) is single-copy; in #3 $\mathcal{C}$(X,Y) is single-copy. Note that #1 and 2 have an "infinite" (single-coordinate) point, but #3 does not.

Theorem 5.1 and the example thus show that coordinatizing a PL-cograph by *edges* is a natural and straightforward process. Coordinatization via *point* coordinates, in contrast, is complex and problematic, requiring the PL-cograph or its equivalent linear geometry to be embedded in a (possibly infinite) projective plane via some of the deepest and most difficult theorems in the field (Batten & Beutelspacher [3]). *Edge*-coordinatization is even sometimes possible in cases more general than PL-cographs: In Counterexample 1 (but not 2) above, for instance, the sets of incident edges determine each point uniquely.

Figure 5.5 below shows all the PL-cographs on $n \leq 7$ points. Some PL-cographs fall into families--for instance, n-point cographs $S_n$ have all edges single-copy, while n-point cographs $K_n$ have just one edge (thus forming complete *graphs*). Small PL-cographs can combine into larger ones, and the table also illustrates naming cographs by decomposition under two natural operations + and $\vee$:

**Definitions**: Let $\mathcal{C}$ and $\mathcal{D}$ be PL-cographs, where the n points and e edges of $\mathcal{C}$ are distinct from the m points and f edges of $\mathcal{D}$. Then define $\mathcal{C} + \mathcal{D}$ to have point set the union of the point sets of $\mathcal{C}$ and $\mathcal{D}$, with $(\mathcal{C} + \mathcal{D})(P,Q) = \mathcal{C}(P,Q)$ if P and Q are both in $\mathcal{C}$, $\mathcal{D}(P,Q)$ if P and Q are both in $\mathcal{D}$, and a new single-copy edge $e_{\{P,Q\}}$ if $P \in \mathcal{C}$ and $Q \in \mathcal{D}$ or if $P \in \mathcal{D}$ and $Q \in \mathcal{C}$. $\mathcal{C} + \mathcal{D}$ thus has n + m points and e + f + nm edges.

Similarly, if all the points of $\mathcal{C}$ and of $\mathcal{D}$ are equivalent (i.e. if $\mathcal{C}$ and $\mathcal{D}$ have either just one edge, or all single-copy edges), define $\mathcal{C} \vee \mathcal{D}$ to have for its n + m - 1 points the points of $\mathcal{C}$ and of $\mathcal{D}$ joined at one point, with $(\mathcal{C} \vee \mathcal{D})(P,Q) = \mathcal{C}(P,Q)$ if P and Q are both in $\mathcal{C} \sim \mathcal{D}$, $\mathcal{D}(P,Q)$ if P and Q are both in $\mathcal{D} \sim \mathcal{C}$, and a new single-copy edge $e_{\{P,Q\}}$ otherwise. $\mathcal{C} \vee \mathcal{D}$ thus has n + m - 1 points and e + f + (n-1)(m-1) edges.

**Examples**: As shown in Figure 5.5, $K_3 + K_3$ is a 6-point 11-edge "prism," and $K_3 \vee K_3$ is a 5-point 6-edge "butterfly."

The number of PL-cographs on n points is:

| n | 2 | 3 | 4 | 5 | 6 | 7 | 8 | 9 |
|---|---|---|---|---|---|---|---|---|
| PL-cographs | 1 | 2 | 3 | 5 | 10 | 24 | 69 | 384 |

| n | 10 | 11 | 12 |
|---|---|---|---|
| PL-cographs | 5250 | 232929 | 28872973 |

Here the values for n = 2 through 7 are from Figure 5.5, and for 8 through 12 are from sequence A001200 of the Sloane *On-line Encyclopedia of Integer Sequences*, "number of linear geometries on n points," with which the n = 2 through 7 values (non-uniquely) match. The table will thus be justified by Theorem 5.2 below identifying PL-cographs with linear spaces (geometries).

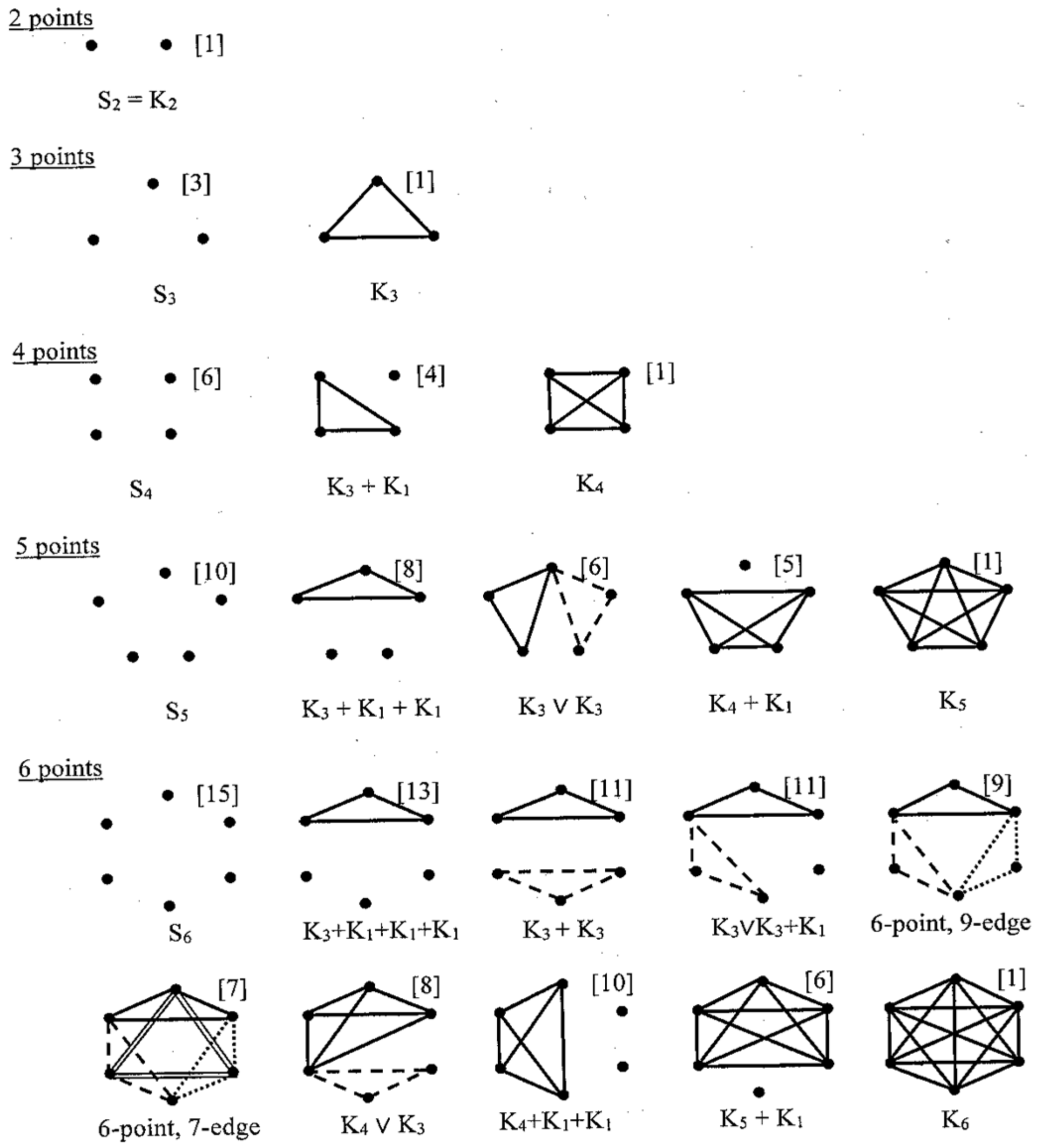

[1] Single-copy edges not shown; number of edges in brackets; 7-point PL cographs sorted by number of blocks with more than two points (i.e. non-single-copy edges).

**Fig. 5.5.** Catalogue of PL-cographs.

(continued next page)

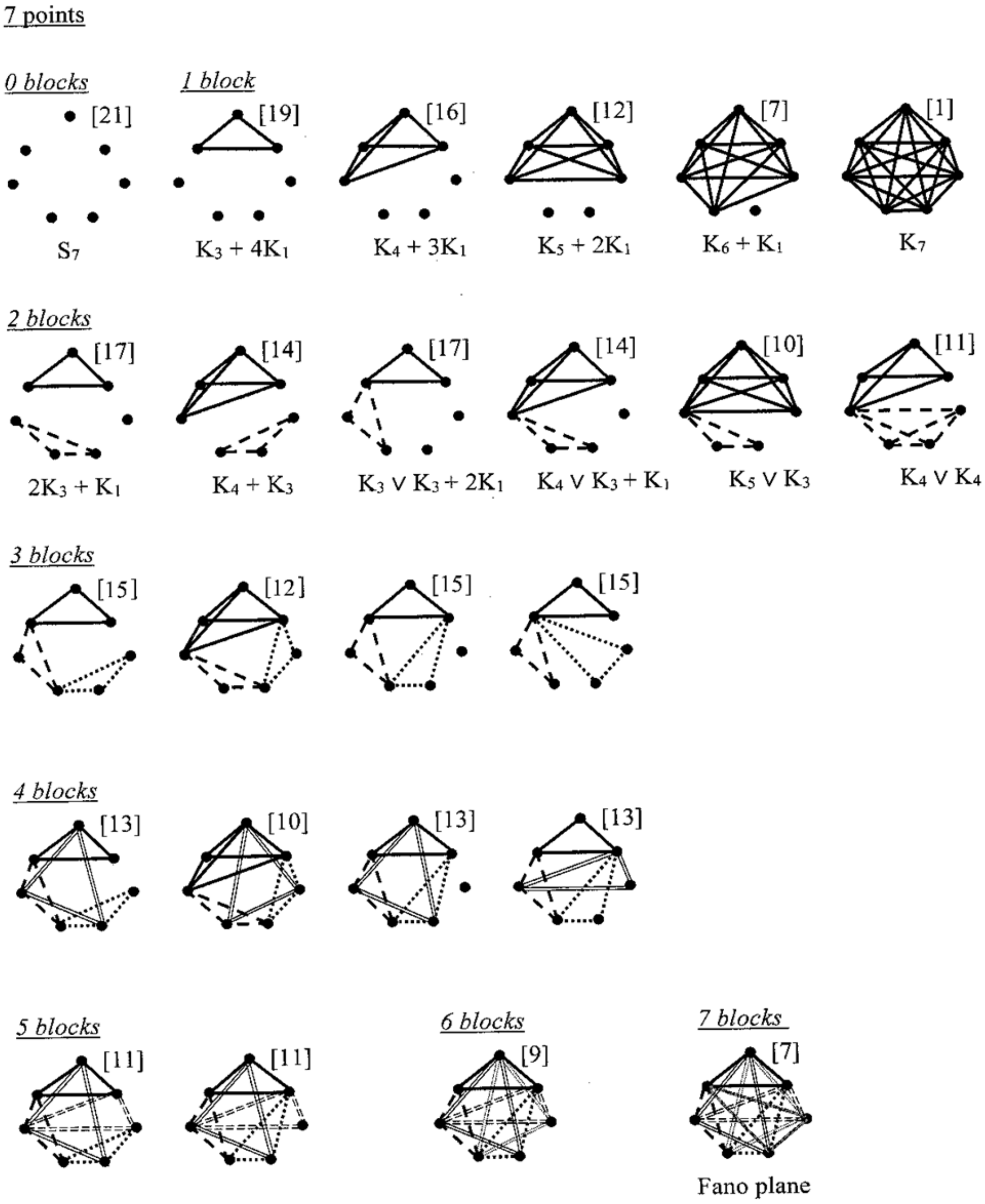

**Fig. 5.5.** Catalogue of PL-cographs.

(continued from previous page)

A linear space $\mathbb{S} = (\wp, \mathcal{L})$ is defined (Batten & Beutelspacher [3]) as a set of "points" $\wp$ and a set of "lines" $\mathcal{L}$ (a line being a set of points) so that:

(1) any two distinct points belong to exactly one line
(2) every line has at least two points.

**Theorem 5.2.** *Linear spaces are equivalent to PL-cographs.*

**Proof**: From a linear space $\mathbb{S} = (\wp, \mathcal{L})$ one may form a cograph $\mathcal{C}$ with points the points $\wp$ of $\mathbb{S}$, edges the lines $\mathcal{L}$ of $\mathbb{S}$, and $\mathcal{C}(P,Q)$ the unique line of $\mathbb{S}$ passing through P and Q. The properties making $\mathcal{C}$ a PL-cograph then follow because the lines of $\mathbb{S}$ are sets.

Conversely, a PL-cograph $\mathcal{C}$ yields a linear space $\mathbb{S} = (\wp, \mathcal{L})$ where $\wp$ is the points of $\mathcal{C}$, and the line through P and Q is the set of all endpoints of the edge $\mathcal{C}(P,Q)$. The PL-cograph rules guarantee that this is well-defined, and that $\mathbb{S}$ is then a linear space.

These two transformations, of linear space to PL-cograph and then back, are easily checked to be inverses, so the two structures are in 1-1 correspondence. ■

PL-cograph blocks correspond under this theorem to linear space lines. Thus, for example, Proposition 5.1(4): two blocks intersect in at most one point, becomes: two lines intersect in at most one point.

For completeness, Figure 5.6 below presents the linear spaces on $n \leq 7$ points corresponding to, and in the same order as, the PL-cographs of Figure 5.5. The point-and-line notation of Figure 5.6 is more efficient for the specialized class of linear spaces than the more general modified cograph notation of Figure 5.5, because it does not "fill in" the blocks. For example, linear spaces in the family $K_n$ are simply those with all their points on one line.

The situation suggests a new classification of the linear spaces, sorting them by the number of *nontrivial* lines (blocks) containing three or more points. Each family is then infinite--e.g., there are one-line spaces $K_n$ for n = 3,4,..., 2-line spaces with n = 5,6,... points, etc. But in a 2-line space the two lines may have either 0 or 1 intersection, in 3-line spaces the lines may be either all parallel, have one intersection, two intersections, or form a triangle, and so on. In other words, one may classify the linear spaces as arising through expansion by adding points (but not nontrivial lines) from "minimal" spaces (having each point on a nontrivial line, and not expanded from any smaller such space). The minimal spaces in Figure 5.6 are marked by ✻. Using the table of linear spaces compiled by Jean Doyen (Batten & Beutelspacher [3]), slightly extended, then yields the following enumeration of minimal linear spaces, in terms of their points or nontrivial lines:

| n | 1 | 2 | 3 | 4 | 5 | 6 | 7 | 8 | 9 |
|---|---|---|---|---|---|---|---|---|---|
| minimal linear spaces: | | | | | | | | | |
| on n points | 0 | 0 | 1 | 0 | 1 | 3 | 8 | 23 | 208 |
| with n non-trivial lines | 1 | 2 | 5 | 16 | $\geq$70 | | | | |

The Sloane *Encyclopedia* does not appear to contain geometry sequences matching either of these, suggesting the minimality approach may be new.

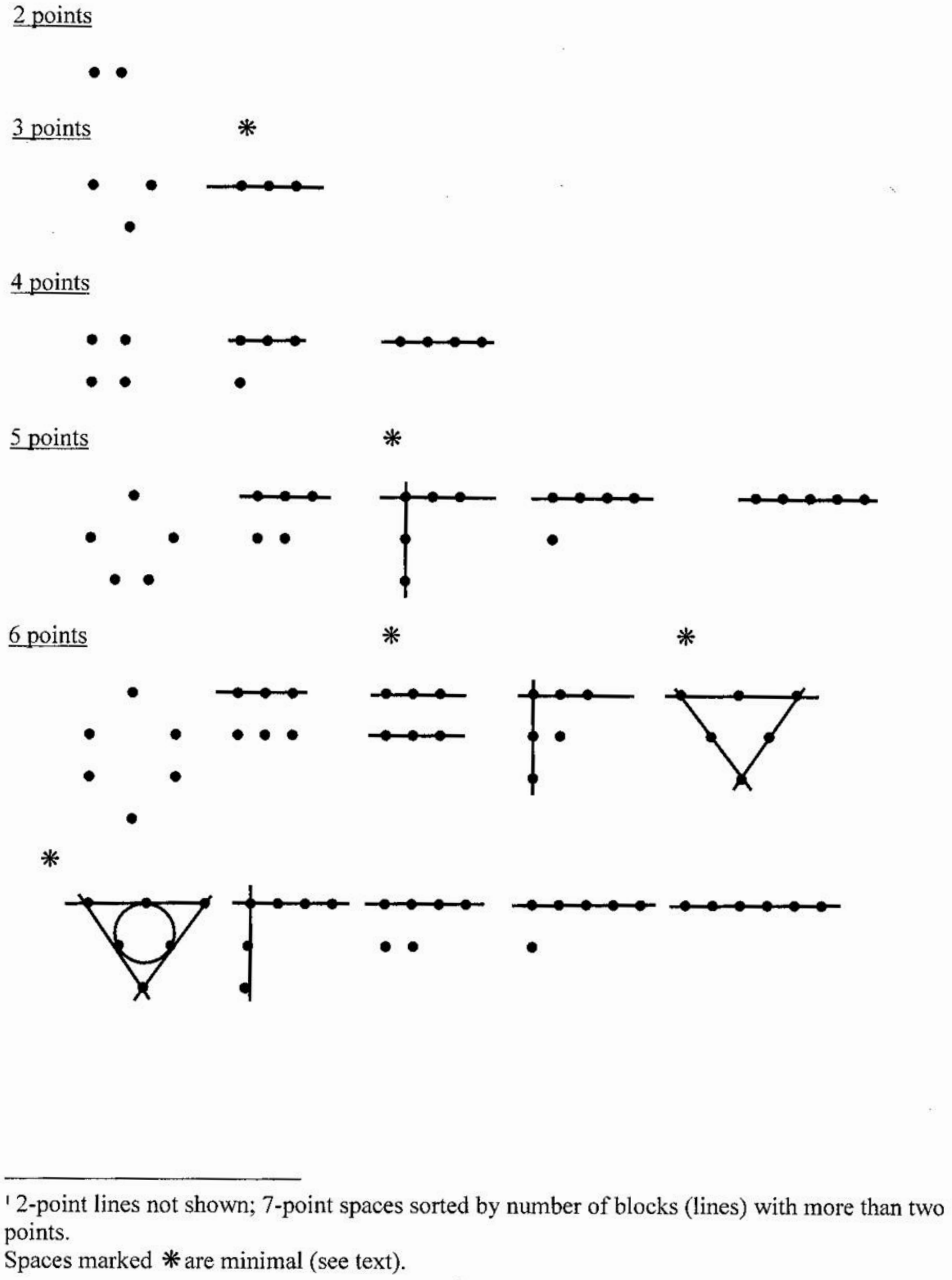

**Fig. 5.6.** Catalogue of linear spaces.

(continued next page)

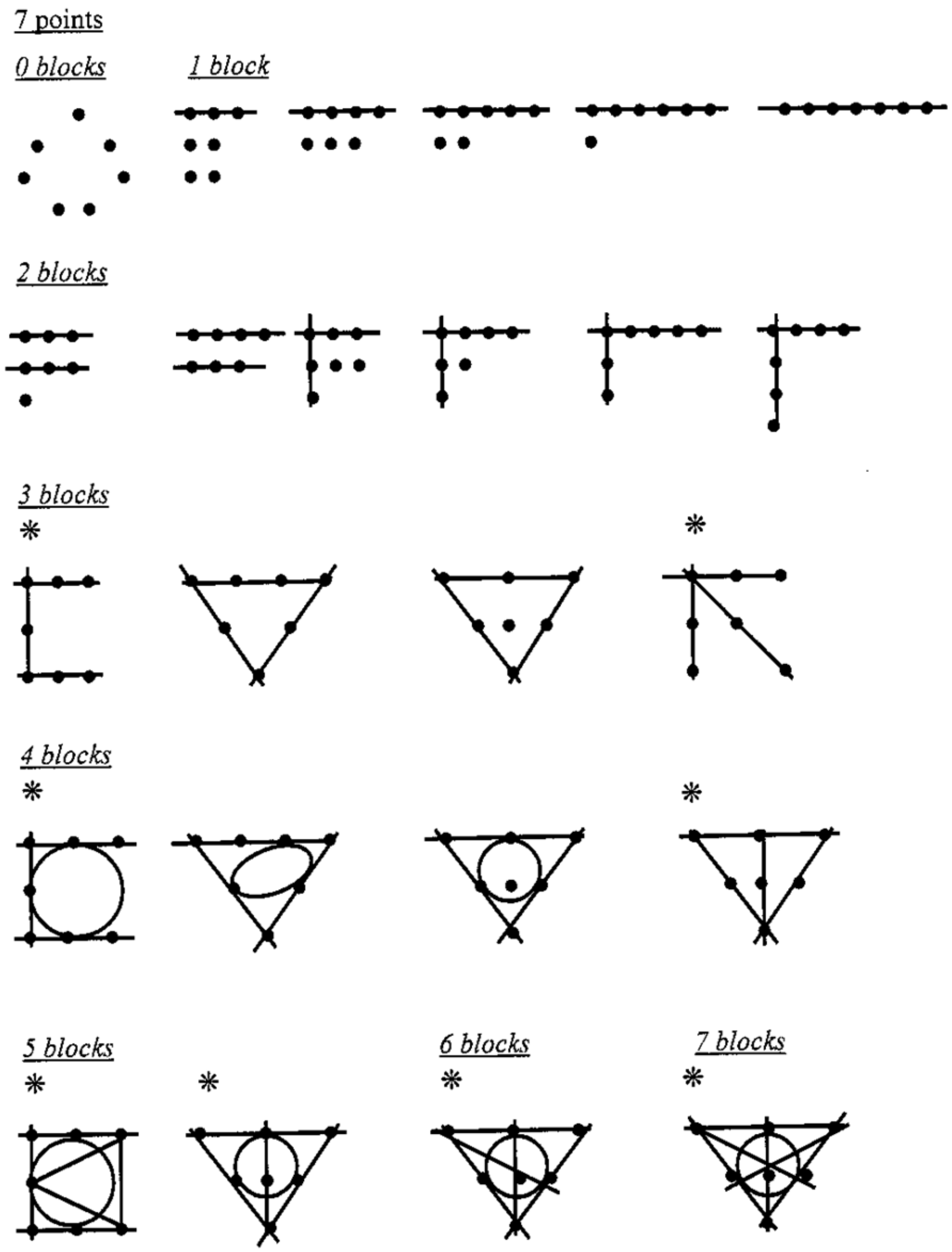

**Fig. 5.6.** Catalogue of linear spaces.

(continued from previous page)

## 6. Group cographs: Ladders, chains, and chain groups

An (additive) abelian group naturally becomes a cograph, as illustrated in Figure 1.4 above, when one defines the edge between each pair of points P and Q to be their sum P+Q. In a general (multiplicative) group, though, there are two possible products PQ and QP, with no evident way to prefer one over the other. To make the group a cograph, one therefore takes the *set* of both:

**Definition**: The *group cograph* $\mathcal{C}$ over a group G has as its points the elements of the group, with $\mathcal{C}(P,Q) = \{PQ,QP\}$.

This section discusses the ladders and chains that characteristically arise from multiple-copy edges in a group cograph, and then solves the interesting group-theoretic problem of characterizing the chain group generated in a finite group from terms that produce a chain.

The first theorem shows how any repeated edge in a group cograph yields whole ladders of repeated edges:

**Theorem 6.1**. *In a group cograph $\mathcal{C}$, if $\mathcal{C}(P,Q) = \mathcal{C}(R,S)$, then there are elements $P_n$ and $Q_n$, for all integers n, so that (possibly switching the names R and S) $P = P_0$, $Q = Q_0$, $R = P_1$, $S = Q_1$, and $\mathcal{C}(P_n,Q_n) = \mathcal{C}(P,Q)$.*

**Proof**. By definition $\{PQ,QP\} = \{RS,SR\}$. Suppose (if necessary switching the names R and S) that $PQ = RS$ and $QP = SR$. Then $R^{-1}PQ = S = QPR^{-1}$, so $PQR = RQP$, implying $QP{\cdot}P^{-1}R = P^{-1}R{\cdot}QP$, whence $X = P^{-1}R$ is in $Z_G(QP)$, the subgroup of G centralizing QP. For all integers *n* define $P_n = PX^n$ and $Q_n = X^{-n}Q$; then $P_nQ_n = PQ$ and $Q_nP_n = QP$, whence $\mathcal{C}(P_n,Q_n) = \mathcal{C}(P,Q)$. ■

**Remarks**: 1) (Diagonal ladders): Note, furthermore, that for all integers *n* and *m*, $P_nQ_m = PQ_{m-n}$ and $Q_mP_n = Q_{m-n}P$, whence $\mathcal{C}(P_n,Q_m) = \mathcal{C}(P,Q_{m-n})$.

2) (Existence): For any pair of points P and Q, a ladder exists containing $\mathcal{C}(P,Q)$: Simply choose X in $Z_G(QP)$ (which always contains at least the subgroup $\langle QP\rangle$), and let $R = PX$, $S = R^{-1}PQ = X^{-1}Q$.

A more specialized situation is when repeated edges join in a chain:

**Definition**. Points $P_0,P_1,\ldots,P_n$ in a group cograph form a *chain* if, for all i, $\mathcal{C}(P_i,P_{i+1}) = \mathcal{C}(P_0,P_1)$. (To rule out trivialities, assume that always $P_i \neq P_{i+1}$ and $P_i \neq P_{i+2}$).

**Theorem 6.2.** *Four or more terms $\{P_i\}$ in a group cograph form a chain if and only if $P = P_0$ and $Q = P_1$ satisfy the three conditions:*

*1) $PQ \neq QP$*

*2) $PQ^2 = Q^2P$*

*3) $P^2Q = QP^2$*

*Then also (defining A = PQ)*

$$P_i = \begin{cases} A^{-i/2}PA^{i/2} & (i \text{ even}) \\ A^{-(i-1)/2}QA^{(i-1)/2} & (i \text{ odd}) \end{cases}$$

*and*

$$(P_i)^2 = \begin{cases} P^2 & (i \text{ even}) \\ Q^2 & (i \text{ odd}) \end{cases}$$

**Proof**. If P,Q,R form a chain then {PQ,QP} = {QR,RQ}. But PQ = RQ would force P = R, which is forbidden. Hence PQ = QR (implying $R = Q^{-1}PQ$) and QP = RQ (implying $R = QPQ^{-1}$, whence $Q^{-1}PQ = R = QPQ^{-1}$, yielding $PQ^2 = Q^2P$. Then also PQ ≠ QP, since PQ = QP would force the forbidden $R = Q^{-1}PQ = P$. Conversely, if P and Q satisfy the given conditions, then defining $R = Q^{-1}PQ$ yields that P,Q,R form a chain. The remainder of the proof is similar, assisted by the observations that $P_2 = R = A^{-1}PA$, and $P_3 = S = A^{-1}QA$, whence also $P^2 = R^2$ and $Q^2 = S^2$, and that (defining B = QP), for any integer i,

1) $PA^i = B^iP$ and $A^iP = PB^i$

2) $QA^i = B^iQ$ and $A^iQ = QB^i$

3) A and B commute

4) $P^2$ and $Q^2$ commute with A and B. ■

**Remarks**. 1) (No branching): Since P and Q determine R uniquely, any chain of three or more terms is uniquely determined by any pair of its adjacent terms. In particular, a chain cannot have forks or alternative branches.

2) (Other associated chains): If $\{P_i\}$ is a chain of four or more terms in a group cograph, then so too are the following:

a) Odd power chains: $\{(P_i)^k\}$ for any fixed odd integer k

b) Odd product chains: $\{\prod_{n=0}^{k-1} P_{i+n}\}$ for any fixed odd integer k

c) j-jump chains: $\{P_{ij+n}\}$ for any fixed integers j and n, provided that $A^j \neq B^j$

A corollary then is that all the $P_i$ must have even order: For if (reindexing) $P = P_0$ had odd order k, $(P_2)^k = (A^{-1}PA)^k = A^{-1}P^kA = 1 = (P_0)^k$, contradicting assertion (a) that $\{(P_i)^k\}$ is a chain.

3) (Cycles): If a chain ever repeats a term, by uniqueness (Remark 1) the entire chain must then form a closed cycle (i.e. $P_{k+i} = P_i$ for all i, and some fixed $k \geq 3$). This will happen if and only if $P = P_0$ and $Q = P_1$, besides the three conditions of the Theorem, also satisfy:

$$\underbrace{PQP \dots}_{k \text{ terms}} = \underbrace{QPQ \dots}_{k \text{ terms}}$$

Thus for any chain $\{P_i\}$ and integer j, either $A^j \neq B^j$ and the j-jump chains $\{P_{ij+n}\}$ exist (Remark 2c), or else $A^j = B^j$ and $\{P_i\}$ forms a 2j-cycle (Remark 3).

4) (Ladders from a chain): The terms of a chain also form ladders as analyzed in the previous theorem. When P, Q from a chain form one rung, and $A^{-1}PA$, $A^{-1}QA$ the next, then the ladder term $X = P^{-1}R$ has a natural expression by commutators: $X = P^{-1}A^{-1}PA = [P,A] = [P,Q]$.

The remainder of this section uses the method of Coxeter & Moser to show that the "chain group" generated by P and Q satisfying the three conditions of the preceding theorem (when P, Q, and A = PQ have finite order) is just a glorified dihedral group.

**Theorem 6.** *In the chain group*

$$C_{p,q,n} = \langle P,Q \mid P^p = Q^q = (PQ)^n = 1,\ PQ \neq QP,\ P^2Q = QP^2,\ PQ^2 = Q^2P \rangle,$$

*p and q are even, and* $P^{2n}Q^{2n} = 1$*. Furthermore,* $C_{p,q,n} = C_{p',q',n}$*, where* $p' = (p,[q,2n])$ *and* $q' = (q,[p,2n])$*, and then (the case* $p' = q' \equiv 2 \bmod 4$*, n = 2 not being possible),* $|C_{p,q,n}| = p'(q',2n)n/2$*, and* $C_{p,q,n}/S \cong D_n$*, where the abelian group* $S = \langle P^2,Q^2 \rangle$ *is central of order* $|S| = p'(q',2n)/4$*.*

**Proof**. $p$ is even, because if $p = 2k+1$ were odd, then $Q = 1 \cdot Q = P^{2k+1}Q = PP^{2k}Q = PQP^{2k} \neq QPP^{2k} = QP^{2k+1} = Q$, a contradiction; and similarly for $q$.

To prove $P^{2n}Q^{2n} = 1$, note that $(PQ)^n = 1$ implies that also $(QP)^n = 1$. Then $1 = 1 \cdot 1 = (PQ)^n(QP)^n = (PQ)^{n-1}PQQP(QP)^{n-1} = (PQ)^{n-1}(QP)^{n-1}P^2Q^2 = \dots = P^{2n}Q^{2n}$.

$C_{p,q,n}$ then reduces to $C_{p',q',n} = C_{(p,[q,2n]),(q,[p,2n]),n}$ (where ( , ) denotes the greatest common divisor, and [ , ] the least common multiple; e.g. $C_{26,16,6} = C_{2,4,6}$). For to prove that $P^p = 1 \Leftrightarrow P^{p'} = 1$ (and similarly with Q), observe that the Euclidean algorithm guarantees integers $a$ and $b$ making $p' = ap + b[q,2n]$, whence $P^{2n}Q^{2n} = 1$ yields $P^{p'} = 1$; conversely $P^{p'} = 1 \Rightarrow P^p = 1$ since $p'$ is a factor of $p$.

It then follows that $|S| = p'(q',2n)/4$: The powers of P in S are $1, P^2, P^4, \dots, P^{p-2}$, and since $P^{2n}Q^{2n} = 1$, one may write any $Q^{2k}$ as $P^{2i}Q^{2j}$ with $0 < 2j < (q',2n)$.

The quotient group $C_{p,q,n}/S \cong D_{n'}$ is then dihedral, where $n'$ is a factor of $n$. For, since P and Q generate $C_{p,q,n}$, their images generate the quotient. But these images are involutions, and the finite group generated by two involutions is always dihedral $D_{n'}$ (since conjugating the product by either one inverts it; Gorenstein [9] 301). Here $n'$ is the order of the image of PQ, which is a factor of $|PQ| = n$.

To prove $C_{p,q,n}/S \cong D_n$ (i.e. $n' = n$, that is, that no nontrivial relation of the form $(PQ)^i = P^{2u}Q^{2v}$ can ever occur), one applies the coset enumeration method of Todd and Coxeter (Coxeter & Moser [6] 12-18). It is simplest to introduce the method by concrete examples, $C_{4,6,6}$ (Table 6.1), and then $C_{4,12,3}$ for odd $n$ and "collapsing" (Table 6.2). The left-most column defines the cosets, labeled by numbers and divided into $b$ "blocks" of length $q' = |Q|$, with the general formula:

$$kq' + i = \langle P \rangle (QP)^k Q^{i-1} \qquad\qquad 0 < k < b\text{-}2,\ 1 < i < q'$$

The numbers in the table represent multiplication of the cosets by the group terms *between* columns at the top (e.g. $2 \cdot P = 7$). Sub-tables I, II, and V encode the relations $P^4 = Q^6 = (PQ)^6 = 1$ by the fact that their left and right columns are the same; sub-tables III and IV encode $P^2Q = QP^2$ and $PQ^2 = Q^2P$ by having their rightmost column match the column to the left of the dashed line. Rows not containing "essentially different" products are suppressed, e.g. rows 2-6 of table II (since $2 \cdot Q$, $3 \cdot Q$, etc. are already in row 1).

The method of Coxeter & Moser consists of progressively filling in the table from the coset definitions and other deduced relations; the method succeeds if the table "closes up," as it does after the 18th row of Table 1, leaving no gaps, inconsistencies, or undefined terms. The circled 7 2 in the first row of table III illustrates the process of deducing relations: Entries 1 1 1 2 1 2 7 follow directly from the coset definitions $1 \cdot P = 1$, $1 \cdot Q = 2$, and $2 \cdot P = 7$. But the final circled entry must be 2, since $P^2Q = QP^2$, thus implying the new relation $7 \cdot P = 2$. This, together with the definition $2 \cdot P = 7$, then allows one to complete row 2 of table I. Once $14 \cdot P = 14$ in table V halts the generation of cosets, the circled deduced entries in tables III and IV allow one to complete the entire table.

The completed coset table gives the order of the group: there are 18 cosets of 4 elements each, so $|C_{4,6,6}| = 72$.

In $C_{4,6,6}$ the relation $P^{2n}Q^{2n} = 1$ deduced at the beginning of the Theorem is vacuous since $P^{12} = Q^{12} = 1$. But otherwise the coset table collapses to a much shorter length, because each block has not $q'$ cosets, but only $(q',2n)$, repeating thereafter. Table 2 illustrates the shortened coset table of $C_{4,12,3}$ (retaining the coset numbering from the "uncollapsed" case). $C_{4,12,3}$ also illustrates the special features of the top block in a coset table with odd $n$: table I is empty, the top block has only half as many cosets as the rest, and its table II consequently cycles through its cosets twice as often.

In the general case, periodicities apparent in the examples simplify the analysis considerably. Observe, first, that $P^2 = 1$ on cosets. For by definition $1 \cdot P = 1$, so also $1 \cdot P^2 = 1$; and a succeeding coset $k$+1 has the form either $k \cdot P$ (inductively unaffected by $P^2$), or $k \cdot Q$ (inductively unaffected by $P^2$ since $QP^2 = P^2Q$). Hence the columns of table I (showing coset multiplication by P) repeat in pairs; and table III holds trivially and will always "close up" if table I does.

Furthermore, each row beyond the first two in a block is the $Q^2$-translate of the row two above it. This is trivially true for the leftmost column in each sub-table, which is just the cosets in ascending order. Columns thereafter arise through multiplication by P or Q at the top of the table, and $Q^2P = PQ^2$ (by definition), and $Q^2Q = QQ^2$ (trivially).

In consequence, the entries in block $k$+1 ($1<k<b$-1) are those of block $k$ shifted forward by $q$. For table V has entries only in the first block, while table II, with entries only in the first row of each block, has the required shift by definition. The other sub-tables, by the above discussion, need be checked only in the first two rows of each block. In table I, by the discussion, one need check only the first two columns. Only even rows appear in these blocks (cf. sample Table 6.1).

The assertion is trivial for column 1, which is just the even numbered cosets; and the row 2 column 2 entry in the $k$th block is by definition always the first coset in the $k$+1st block, yielding the required pattern. Tables III and IV follow because multiplication by P has the desired shift (via table I) and changes the block by at most 1, and $P^2 = 1$, while multiplication by Q has the desired shift (via table II) and leaves the block unchanged.

Since applying PQ repetitively moves a coset steadily first up the blocks, then back down (table V), the constraint $(PQ)^n = 1$ dictates that the number of blocks is $b = n/2$. (When $n$ is odd, the top block is only half-length, as in Table 6.2.) In more detail, the values $1\cdot(PQ)^i$ in the left half of the first row of table V are $2+(i-1)q'$ for $i = 1,2,...$, while in the right half the values $1\cdot(PQ)^{-i}P$ are $2+(i-1)q' - 2i$ for $i = 1,2,...$ (Here $2+(i-1)q'$ is always the second coset in block $i$, and $-2i$ denotes an offset taken modulo $(q',2n)$ within that block.) Hence when $n$ is even, the table V row will close up at its central meeting point $i = b = n/2$, at the product $[2+(b-1)q']\cdot P = 2+(b-1)q' - 2b = 2+(b-1)q' - n$. The offset $-n$, congruent to $n$ mod $(q,2n)$ [which is 0 or $(q'/2,n)$, depending as the 2-factor of $q'$ is < or > than the 2-factor of $n$], is consistent with $P^2 = 1$, so allowing table I to close. When $n$ is odd, the central meeting point is $[1 + (i-1)q']\cdot Q = (i-1)q' + (q',2n) -2i+3$ at $i = (n+1)/2$, an offset mod $(q',2n)$ of $1-n$. But Q, by definition in table II, offsets by 1. Hence, to make $-n \equiv 0$, the top block must collapse by half, to a length of $(q',n)$ rather than $(q',2n)$. In summary, the requirement that table V "close up" dictates: when $n$ is even, a P-product allowing table I to close up in block $b$; when $n$ is odd, a Q-product allowing table II to close up provided that block $b$ collapse to half size as in the example Table 6.2.

The discussed periodicities reduce the analysis to checking just the first two rows of just the first, second, and last ($b$th) blocks of the table; simple direct calculations now show that these behave properly, so the coset table is consistent and "closes up." In particular, $C_{p,q,n}$ has $n/2$ blocks each containing $(q',2n)$ cosets of $p'$ elements, proving that $|C_{p,q,n}| = p'(q',2n)n/2$.

Since the coset table proof given so far was actually for the group defined by generators and relations without the constraint $PQ \neq QP$, it remains necessary, finally, to derive $PQ \neq QP$ from the rest. P and Q cannot commute when $n > 2$ because they generate the nonabelian quotient $D_n$. When $n = 2$, $QP = P^{-1}Q^{-1}$ (trivially) and $1 = P^{2n}Q^{2n} = P^4Q^4$, whence $QP = P^{-1}(P^4Q^4)Q^{-1} = PQP^2Q^2$, proving $PQ \neq QP$ unless $P^2Q^2 = 1$. The latter does actually hold in the excluded case $n = 2$, $p' = q' = 4k + 2$, since then $1 = 1\cdot 1 = P^{p'}Q^{q'} = (P^4Q^4)^kP^2Q^2 = P^2Q^2$. But otherwise the coset table for $C_{p',q',2}$ closes up successfully after four rows, on account of first row coset products being either $2\cdot P = 4$, $3\cdot P = 3$, and $4\cdot P = 2$ in subtables V, IV, and III, respectively, or else $2\cdot P = 2$, $3\cdot P = 3$, and $4\cdot P = 4$, the generation of cosets $1 = \langle P\rangle$, $2 = 1\cdot Q$, $3 = 2\cdot Q$, $4 = 3\cdot Q$ halting at 4 because $Q^4 = P^{-4} \in 1$. The relation $P^2Q^2 = 1$ (which would halt the cosets already at 2) thus is not true; hence $PQ \neq QP$. ■

| Coset Definitions | I | II | III | | IV | | V |
|---|---|---|---|---|---|---|---|
| | P P P P | Q QQ QQQ | P P Q | Q P P | P Q Q | Q Q P | P QP QP Q P Q P Q P Q |
| 1 = <P> | 1 1 1 1 1 | 1 2 3 4 5 6 1 | 1 1 1 2 | 1 2 7 2 | 1 1 2 3 | 1 2 3 3 | 1 1 2 7 8 13 14 14 15 10 11 6 1 |
| 2 = 1·Q | 2 7 2 7 2 | | 2 7 2 3 | 2 3 3 3 | 2 7 8 9 | 2 3 4 9 | |
| 3 = 2·Q | 3 3 3 3 3 | | 3 3 3 4 | 3 4 9 4 | 3 3 4 5 | 3 4 5 5 | 3 3 4 9 10 15 16 16 17 12 7 2 3 |
| 4 = 3·Q | 4 9 4 9 4 | | 4 9 4 5 | 4 5 5 5 | 4 9 10 11 | 4 5 6 11 | |
| 5 = 4·Q | 5 5 5 5 5 | | 5 5 5 6 | 5 6 11 6 | 5 5 6 1 | 5 6 1 1 | 5 5 6 11 12 17 18 18 13 8 9 4 5 |
| 6 = 5·Q | 6 11 6 11 6 | | 6 11 6 1 | 6 1 1 1 | 6 11 | 12 7 6 1 | 2 7 |
| 7 = 2·P | | 7 8 9 10 11 12 7 | 7 2 7 8 | 7 8 13 8 | 7 2 3 4 | 7 8 9 4 | |
| 8 = 7·Q | 8 13 8 13 8 | | 8 13 8 9 | 8 9 4 9 | 8 13 14 15 | 8 9 10 15 | |
| 9 = 8·Q | | | 9 4 9 10 | 9 10 15 10 | 9 4 5 6 | 9 10 11 6 | |
| 10 = 9·Q | 10 15 10 15 10 | | 10 15 10 11 | 10 11 6 11 | 10 15 16 17 | 10 11 12 17 | |
| 11 = 10·Q | | | 11 6 11 12 | 11 12 17 12 | 11 6 1 2 | 11 12 7 2 | |
| 12 = 11·Q | 12 17 12 17 12 | | 12 17 12 7 | 12 7 2 7 | 12 17 18 13 | 12 7 8 13 | |
| 13 = 8·P | | 13 14 15 16 17 18 13 | 13 8 13 14 | 13 14 14 14 | 13 8 9 10 | 13 14 15 10 | |
| 14 = 13·Q | 14 14 14 14 14 | | 14 14 14 15 | 14 15 10 15 | 14 14 15 16 | 14 15 16 16 | |
| 15 = 14·Q | | | 15 10 15 16 | 15 16 16 16 | 15 10 11 12 | 15 16 17 12 | |
| 16 = 15·Q | 16 16 16 16 16 | | 16 16 16 17 | 16 17 12 17 | 16 16 17 18 | 16 17 18 18 | |
| 17 = 16·Q | | | 17 12 17 18 | 17 18 18 18 | 17 12 7 8 | 17 18 13 8 | |
| 18 = 17·Q | 18 18 18 18 18 | | 18 18 18 13 | 18 13 8 13 | 18 18 13 14 | 18 13 14 14 | |

**Table 6.1.** $C_{4,6,6}$.

| Coset Definitions | I | II | III | | IV | | V |
|---|---|---|---|---|---|---|---|
| | P P P P | Q Q Q Q Q Q Q Q Q Q Q Q | P P Q | Q P P | P Q Q | Q Q P | P Q P Q P Q |
| 1 = <P> | 1 1 1 1 1 | 1 2 3 4 5 6 1 2 3 4 5 6 1 | 1 1 1 2 | 1 2 13 2 | 1 1 2 3 | 1 2 3 3 | 1 1 2 13 14 6 1 |
| 2 = 1·Q | 2 13 2 13 2 | | 2 13 2 3 | 2 3 3 3 | 2 13 14 15 | 2 3 4 15 | |
| 3 = 2·Q | 3 3 3 3 3 | | 3 3 3 4 | 3 4 15 4 | 3 3 4 5 | 3 4 5 5 | 3 3 4 15 13 2 3 |
| 4 = 3·Q | 4 15 4 15 4 | | 4 15 4 5 | 4 5 5 5 | 4 15 13 14 | 4 5 6 14 | |
| 5 = 4·Q | 5 5 5 5 5 | | 5 5 5 6 | 5 6 14 6 | 5 5 6 1 | 5 6 1 1 | 5 5 6 14 15 4 5 |
| 6 = 5·Q | 6 14 6 14 6 | | 6 14 6 1 | 6 1 1 1 | 6 14 15 13 | 6 1 2 13 | |
| 13 = 2·P | | 13 14 15 13 14 15 13 14 15 13 14 15 13 | 13 2 13 14 | 13 14 6 14 | 13 2 3 4 | 13 14 15 4 | |
| 14 = 13·Q | | | 14 6 14 15 | 14 15 4 15 | 14 6 1 2 | 14 15 13 2 | |
| 15 = 14·Q | | | 15 4 15 13 | 15 13 2 13 | 15 4 5 6 | 15 13 14 6 | |

**Table 6.2.** $C_{4,12,3}$.

## Acknowledgments

I am grateful for the encouragement of my late advisor, Prof. Peter Hilton, at the beginning of this project, and then for the opportunity to present my results in a series of talks at John Carroll University in 2000, and at section meetings of the Mathematical Association of America held at John Carroll University in 2004, and Cleveland State University in 2013. Technology trainer Ann MacNamara, and other Cleveland Heights - University Heights public library staff, assisted my struggles to PC-computer format my manuscript.

## References

[1] Adams, Colin C. *The Knot Book: An Elementary Introduction to the Mathematical Theory of Knots,* New York: W. H. Freeman, 1994.
[2] Bach, Anna Magdalena. Minuet in g, in: *Notebook of Anna Magdalena Bach.* www.youtube.com/user/The Great Repertoire, 2015.
[3] Batten, Lynn Margaret, and Beutelspacher, Albrecht. *The Theory of Finite Linear Spaces: Combinatorics of Points and Lines*, Cambridge: Cambridge University Press, 1993.
[4] Birkhoff, George D. *Aesthetic Measure*, Cambridge, MA: Harvard University Press, 1933.
[5] Borsuk, K. *Multidimensional Analytic Geometry,* Polish Scientific Publishers, 1969.
[6] Coxeter, H. S. M., and Moser, W. O. J. *Generators and Relations for Discrete Groups*, Berlin: Springer-Verlag, 1957.
[7] Dickinson, Emily. "Safe in their Alabaster Chambers," In: *The Complete Poems of Emily Dickinson* (Thomas H. Johnson, ed.), Boston: Little, Brown, and Co., 1960, No. 216, p. 100.
[8] Forte, Allen, and Gilbert, Steven E. *Introduction to Schenkerian Analysis*, New York: W. W. Norton, 1982.
[9] Gorenstein, Daniel. *Finite Groups*, New York: Harper & Row, 1968.
[10] Haas, Robert. "Three-colorings of finite groups, or An algebra of nonequalities," *Mathematics Magazine 63*(4) (1990) 211-225.
[11] Hall, Marshall, Jr. *The Theory of Groups*, 2nd ed., New York: Chelsea, 1976.
[12] Harary, Frank, and Palmer, Edgar M. *Graphical Enumeration*, New York: Academic Press, 1973.
[13] Johnson, Thomas H. *Emily Dickinson: An Interpretive Biography*, Cambridge, Mass.: The Belknap Press of Harvard Univ. Press, 1955, pp. 104-108.
[14] Kauffman, Louis H. "New invariants in the theory of knots," *American Math. Monthly 95* (1988) 195-242.
[15] Koshy, Thomas. *Fibonacci and Lucas Numbers with Applications* (Wiley-Interscience series, Pure and Applied Mathematics), New York: John Wiley & Sons, 2001.
[16] Mu Qi. Six Persimmons, in: Sherman Lee, *A History of Far Eastern Art* (4th edition), New York: Harry N. Abrams, Inc., 1982, p. 362.
[17] Sloane, Neil J. A. *On-line Encyclopedia of Integer Sequences*, http://www.research.att.com /cgi-bin/access.cgi/as/njas/sequences/eismum.cgi.
[18] Waley, Arthur. *An Introduction to the Study of Chinese Painting*, London: Benn, 1923, p. 231. Cited in: Sherman Lee, *A History of Far Eastern Art* (4th edition), New York: Harry N. Abrams, Inc., 1982, p. 362.
[19] Weyl, Hermann. *Symmetry*, Princeton: Princeton Univ. Press, 1952.

[i] Remark on nomenclature: The cographs in this article are quite distinct from the combinatorial objects termed cographs. The term used in this article serves to recall 1) the familiar category-theoretic duality, "reversing the arrows," in Definition 1.1, and 2) the complete colored graphs in Definition 1.2. I have used the term privately since 1980, in lectures in 2000, and in talks to local section meetings of the MAA, in 2004 and 2013. Some evidence that cographs have not been extensively studied under some other name is that the sequence 1, 3, 25, 1299, the number of cographs on $n$ points (Table 1.1), is not currently (2019) in the *On-Line Encyclopedia of Integer Sequences* [17].